\documentclass[10pt]{extarticle}

\usepackage[english]{babel}
\usepackage{graphicx}
\usepackage{framed}
\usepackage[normalem]{ulem}
\usepackage{indentfirst}
\usepackage{amsmath,amsthm,amssymb,amsfonts}
\usepackage{bbm}
\usepackage[italicdiff]{physics}
\usepackage[T1]{fontenc}
\usepackage{wrapfig} 
\usepackage{lmodern,mathrsfs}
\usepackage[inline,shortlabels]{enumitem}
\setlist{topsep=2pt,itemsep=2pt,parsep=0pt,partopsep=0pt}
\usepackage[dvipsnames]{xcolor}
\usepackage[utf8]{inputenc}
\usepackage[a4paper,top=0.7in,bottom=0.6in,left=0.7in,right=0.7in]{geometry} 
\usepackage[most]{tcolorbox}
\usepackage{tikz,tikz-3dplot,tikz-cd,tkz-tab,tkz-euclide,pgf,pgfplots}
\pgfplotsset{compat=newest}
\usepackage{multicol}
\usepackage[bottom,multiple]{footmisc} 
\usepackage[backend=bibtex,style=numeric]{biblatex}
\usepackage{xurl} 
\usepackage{hyperref}
\usepackage[nameinlink]{cleveref}

\usepackage{fancyhdr}
\fancypagestyle{plain}{
\fancyhf{}

\fancyfoot[C]{\sffamily\thepage}}

\newcommand{\hide}[1]{} 

\definecolor{contcol1}{HTML}{48C6EF}
\definecolor{contcol2}{HTML}{6F86D6}
\definecolor{convcol1}{HTML}{EA8D8D}
\definecolor{convcol2}{HTML}{A890FE}
\definecolor{abscol1}{HTML}{E89403}
\definecolor{abscol2}{HTML}{F23912}
\definecolor{cdcol1}{HTML}{02AAB0}
\definecolor{cdcol2}{HTML}{00CDAC}
\definecolor{conclcol1}{HTML}{24D8D2}
\definecolor{conclcol2}{HTML}{73E8C9}

\newcommand{\ep}{\varepsilon}

\newcommand{\lam}{\lambda}

\renewcommand{\ip}[1]{\ensuremath{\left\langle#1\right\rangle}}

\newcommand{\mmat}[4]{\ensuremath{\begin{pmatrix} #1 & #2 \\ #3 & #4 \end{pmatrix}}}

\newcommand{\C}{\mathbb{C}}

\newcommand{\F}{\mathbb{F}}
\newcommand{\K}{\mathbb{K}}
\newcommand{\N}{\mathbb{N}}

\newcommand{\R}{\mathbb{R}}

\newcommand{\Z}{\mathbb{Z}}
\newcommand{\As}{\mathcal{A}}
\newcommand{\Bs}{\mathcal{B}}

\newcommand{\Hs}{\mathcal{H}}
\newcommand{\Is}{\mathcal{I}}
\newcommand{\Js}{\mathcal{J}}
\newcommand{\Ks}{\mathcal{K}}
\newcommand{\Ls}{\mathcal{L}}
\newcommand{\Ms}{\mathcal{M}}
\newcommand{\Ns}{\mathcal{N}}

\newcommand{\Ps}{\mathcal{P}}

\newcommand{\Rs}{\mathcal{R}}
\newcommand{\Ss}{\mathcal{S}}

\newcommand{\Vs}{\mathcal{V}}

\newcommand{\xb}{\textbf{x}}

\newcommand{\zbar}{\overline{z}}

\newcommand{\Dbar}{\overline{D}}

\newcommand{\1}{\mathbbm{1}} 

\newcommand{\limm}{\lim_{m\to\infty}}
\newcommand{\limn}{\lim_{n\to\infty}}

\newcommand{\liminfn}{\liminf_{n\to\infty}}
\newcommand{\limsupn}{\limsup_{n\to\infty}}
\newcommand{\sumkn}{\sum_{k=1}^n}

\newcommand{\sumn}[1][1]{\sum_{n=#1}^\infty}
\newcommand{\emp}{\varnothing}
\newcommand{\exc}{\backslash}
\newcommand{\sub}{\subseteq}
\newcommand{\sups}{\supseteq}
\newcommand{\capp}{\bigcap}

\newcommand{\dx}{\,dx}
\newcommand{\dy}{\,dy}

\newcommand{\ahat}{\widehat{a}}
\newcommand{\bhat}{\widehat{b}}
\newcommand{\fh}{\widehat{f}}
\newcommand{\ph}{\widehat{p}}
\newcommand{\xh}{\widehat{x}}

\newcommand{\inv}{^\times}

\renewcommand{\iff}{\Leftrightarrow}
\DeclareMathOperator{\im}{\text{im}}
\let\Re\relax\let\Im\relax

\DeclareMathOperator{\Re}{\text{Re}}
\DeclareMathOperator{\Im}{\text{Im}}

\newtheoremstyle{mystyle}{}{}{}{}{\sffamily\bfseries}{.}{ }{}
\newtheoremstyle{cstyle}{}{}{}{}{\sffamily\bfseries}{.}{ }{\thmnote{#3}}
\makeatletter
\renewenvironment{proof}[1][\proofname] {\par\pushQED{\qed}{\normalfont\sffamily\bfseries\topsep6\p@\@plus6\p@\relax #1\@addpunct{.} }}{\popQED\endtrivlist\@endpefalse}
\makeatother
\theoremstyle{mystyle}{\newtheorem{definition}{Definition}[section]}
\theoremstyle{mystyle}{\newtheorem{proposition}[definition]{Proposition}}
\theoremstyle{mystyle}{\newtheorem{theorem}[definition]{Theorem}}
\theoremstyle{mystyle}{\newtheorem{lemma}[definition]{Lemma}}
\theoremstyle{mystyle}{\newtheorem{corollary}[definition]{Corollary}}
\theoremstyle{mystyle}{\newtheorem*{remark}{Remark}}
\theoremstyle{mystyle}{\newtheorem*{remarks}{Remarks}}
\theoremstyle{mystyle}{\newtheorem*{example}{Example}}
\theoremstyle{mystyle}{\newtheorem*{examples}{Examples}}
\theoremstyle{definition}{}
\theoremstyle{cstyle}{\newtheorem*{cthm}{}}

\newtheoremstyle{warn}{}{}{}{}{\normalfont}{}{ }{}
\theoremstyle{warn}
\newtheorem*{warning}{\warningsign{0.2}\relax}

\newcommand{\warningsign}[1]{\tikz[scale=#1,every node/.style={transform shape}]{\draw[-,line width={#1*0.8mm},red,fill=yellow,rounded corners={#1*2.5mm}] (0,0)--(1,{-sqrt(3)})--(-1,{-sqrt(3)})--cycle;
\node at (0,-1) {\fontsize{48}{60}\selectfont\bfseries!};}}

\tcolorboxenvironment{definition}{boxrule=0pt,boxsep=0pt,colback={red!10},left=8pt,right=8pt,enhanced jigsaw, borderline west={2pt}{0pt}{red},sharp corners,before skip=10pt,after skip=10pt,breakable}
\tcolorboxenvironment{proposition}{boxrule=0pt,boxsep=0pt,colback={Orange!10},left=8pt,right=8pt,enhanced jigsaw, borderline west={2pt}{0pt}{Orange},sharp corners,before skip=10pt,after skip=10pt,breakable}
\tcolorboxenvironment{theorem}{boxrule=0pt,boxsep=0pt,colback={blue!10},left=8pt,right=8pt,enhanced jigsaw, borderline west={2pt}{0pt}{blue},sharp corners,before skip=10pt,after skip=10pt,breakable}
\tcolorboxenvironment{lemma}{boxrule=0pt,boxsep=0pt,colback={Cyan!10},left=8pt,right=8pt,enhanced jigsaw, borderline west={2pt}{0pt}{Cyan},sharp corners,before skip=10pt,after skip=10pt,breakable}
\tcolorboxenvironment{corollary}{boxrule=0pt,boxsep=0pt,colback={violet!10},left=8pt,right=8pt,enhanced jigsaw, borderline west={2pt}{0pt}{violet},sharp corners,before skip=10pt,after skip=10pt,breakable}
\tcolorboxenvironment{proof}{boxrule=0pt,boxsep=0pt,blanker,borderline west={2pt}{0pt}{CadetBlue!80!white},left=8pt,right=8pt,sharp corners,before skip=10pt,after skip=10pt,breakable}
\tcolorboxenvironment{remark}{boxrule=0pt,boxsep=0pt,blanker,borderline west={2pt}{0pt}{Green},left=8pt,right=8pt,before skip=10pt,after skip=10pt,breakable}
\tcolorboxenvironment{remarks}{boxrule=0pt,boxsep=0pt,blanker,borderline west={2pt}{0pt}{Green},left=8pt,right=8pt,before skip=10pt,after skip=10pt,breakable}
\tcolorboxenvironment{example}{boxrule=0pt,boxsep=0pt,blanker,borderline west={2pt}{0pt}{Black},left=8pt,right=8pt,sharp corners,before skip=10pt,after skip=10pt,breakable}
\tcolorboxenvironment{examples}{boxrule=0pt,boxsep=0pt,blanker,borderline west={2pt}{0pt}{Black},left=8pt,right=8pt,sharp corners,before skip=10pt,after skip=10pt,breakable}
\tcolorboxenvironment{cthm}{boxrule=0pt,boxsep=0pt,colback={gray!10},left=8pt,right=8pt,enhanced jigsaw, borderline west={2pt}{0pt}{gray},sharp corners,before skip=10pt,after skip=10pt,breakable}

\newtcbox{\personbox}[2]{
title={#2},
flip title={interior hidden}, 
fonttitle=\sffamily\selectfont,center title,
colframe=#1,coltitle=black,
hbox,boxsep=0mm, 
top=0mm,bottom=0mm,left=0mm,right=0mm,toptitle=1mm,bottomtitle=1mm, 
clip upper, 
enhanced}

\newtcbox{\conclbox}{
interior style={left color=conclcol1!60,right color=conclcol2!60},
frame style={left color=conclcol1!80!black,right color=conclcol2!80!black},
halign upper=center,
fontupper=\sffamily\selectfont,
hbox,center, 
top=1mm,bottom=1mm,left=1mm,right=1mm,arc=3.5mm,
before skip=50pt,after skip=10pt, 
enhanced}

\newenvironment{talign*}{\let\displaystyle\textstyle\csname align*\endcsname}{\endalign}

\usepackage[explicit]{titlesec}
\titleformat{\section}{\fontsize{24}{30}\sffamily\bfseries}{\thesection}{20pt}{#1}
\titleformat{\subsection}{\fontsize{16}{18}\sffamily\bfseries}{\thesubsection}{12pt}{#1}
\titleformat{\subsubsection}{\fontsize{10}{12}\sffamily\large\bfseries}{\thesubsubsection}{8pt}{#1}

\titlespacing*{\section}{0pt}{5pt}{5pt}
\titlespacing*{\subsection}{0pt}{5pt}{5pt}
\titlespacing*{\subsubsection}{0pt}{5pt}{5pt}

\newcommand{\xnewpage}{\newpage} 
  
\newcommand{\Disp}{\displaystyle}

\DeclareMathAlphabet\mathbfcal{OMS}{cmsy}{b}{n}
\newcommand{\ind}{\hspace{0.2in}} 

\newcommand{\cst}{$C^*$}
\newcommand{\csta}{$C^*$-algebra}
\newcommand{\cstas}{$C^*$-algebras}

\makeatletter
\g@addto@macro\normalsize{
\setlength\abovedisplayskip{3pt}
\setlength\belowdisplayskip{3pt}
\setlength\abovedisplayshortskip{0pt}
\setlength\belowdisplayshortskip{0pt}}
\makeatother

\makeatletter
\def\supervisor#1{\gdef\@supervisor{#1}}
\def\@supervisor{\@latex@warning@no@line{No \noexpand\supervisor given}}
\def\module#1{\gdef\@module{#1}}
\def\@module{\@latex@warning@no@line{No \noexpand\module given}}
\renewcommand\maketitle{
\begin{center}
{\fontsize{22}{26}\sffamily\bfseries\selectfont\@title}\\
\vspace{8mm}
{\fontsize{24}{28}\sffamily\selectfont\@author}\\
\vspace{12mm}
{\fontsize{20}{24}\sffamily\selectfont Project Supervisor: \@supervisor}\\
\vspace{6mm}
{\fontsize{18}{22}\sffamily\selectfont Module: \@module}\\
\vspace{6mm}
{\fontsize{18}{22}\sffamily\selectfont\@date}
\end{center}}
\makeatother

\makeatletter
\renewenvironment{abstract}{
\if@twocolumn
\section*{\abstractname}
\else
\begin{center}
{\large\sffamily\bfseries\abstractname\vspace{\z@}}
\end{center}
\quotation
\fi}
{\if@twocolumn\else\endquotation\fi}
\renewenvironment{quotation}{
\list{}{
\rightmargin\leftmargin
\parsep\z@\@plus\p@}
\item\relax}
\makeatother

\title{\cstas: The (Quantum) Path Less Traveled}
\author{Senan Sekhon}
\supervisor{Benjamin Doyon}
\module{6CCM345A}
\date{2020/21}

\begin{document}

\begin{titlepage}
\vspace*{25mm}

\begin{center}
\begin{tikzpicture}[complexnode/.pic={
\begin{axis}[scale=1,axis lines=none,axis on top,xtick=\empty,ytick=\empty,ztick=\empty,xrange=-4:4,yrange=-4:4,unit vector ratio=1 1 1]
\addplot3[domain=-4:4,y domain=-4:4,colormap/viridis,surf,samples=51] {((x^2+y^2)/2-(x^4+y^4)/12+(x^6+y^6)/288)};
\end{axis}
}]
\newcommand{\wavepart}[8][0.2]{
\begin{scope}[xshift=#6cm,yshift=#7cm,rotate=#8,transform shape]
\shade[left color=#2,right color=white,shading angle=#8+90] (0,#5) sin (#3/8,#4+#5/2) cos (2*#3/8,#5) sin (3*#3/8,-#4+#5/2) cos (4*#3/8,#5) sin (5*#3/8,#4+#5/2) cos (6*#3/8,#5) sin (7*#3/8,-#4+#5/2) cos (8*#3/8,#5) -- (8*#3/8,-#5) sin (7*#3/8,-#4-#5/2) cos (6*#3/8,-#5) sin (5*#3/8,#4-#5/2) cos (4*#3/8,-#5) sin (3*#3/8,-#4-#5/2) cos (2*#3/8,-#5) sin (#3/8,#4-#5/2) cos (0*#3/8,-#5) -- cycle;
\path[shading=ball,ball color=#2] (0,0) circle (#1);
\end{scope}
}
\clip (-1,-0.5) rectangle (14.7,11.5);
\draw (0,0) pic[scale=2] {complexnode};
\wavepart[0.2]{red}{5}{0.5}{0.1}{7}{7}{20}
\wavepart[0.2]{orange}{4.2}{0.5}{0.1}{6}{11}{-10}
\wavepart[0.2]{yellow!80!black}{3.4}{0.5}{0.1}{0}{2}{30}
\wavepart[0.2]{green!50!black}{3.4}{0.5}{0.1}{3}{10}{-160}
\wavepart[0.2]{blue}{3}{0.5}{0.1}{4}{7}{150}
\wavepart[0.2]{violet}{2.4}{0.5}{0.1}{13}{3}{-140}
\end{tikzpicture}
\end{center}

\vspace*{20mm}
\maketitle

\end{titlepage}
\thispagestyle{empty} 
\addtocounter{page}{-1} 

\begin{tcolorbox}[title=Contents, fonttitle=\huge\sffamily\bfseries\selectfont,interior style={left color=contcol1!40!white,right color=contcol2!40!white},frame style={left color=contcol1!80!white,right color=contcol2!80!white},coltitle=black,top=2mm,bottom=2mm,left=2mm,right=2mm,drop fuzzy shadow,enhanced]
\makeatletter
\sffamily\selectfont\@starttoc{toc}
\makeatother
\end{tcolorbox}

\vspace*{10mm}

\begin{tcolorbox}[title=Conventions, fonttitle=\large\sffamily\bfseries\selectfont,interior style={left color=convcol1!40!white,right color=convcol2!40!white},frame style={left color=convcol1!80!white,right color=convcol2!80!white},coltitle=black,top=2mm,bottom=2mm,left=2mm,right=2mm,drop fuzzy shadow,enhanced]
\begin{itemize}[leftmargin=0.15in]
    \item $\F$ denotes either $\R$ or $\C$.
    \item $\N$ denotes the set $\{1,2,3,...\}$ of natural numbers (excluding $0$).
    \item $\N_0$ denotes the set $\{0,1,2,3,...\}$ of natural numbers (including $0$).
    \item $D(z_0,r)$ denotes the open disk $\{z\in\C\mid\abs{z-z_0}<r\}$ in the complex plane.
    \item $\Dbar(z_0,r)$ denotes the closed disk $\{z\in\C\mid\abs{z-z_0}\le r\}$ in the complex plane.
    \item Inner products are assumed to be linear in the first argument and conjugate linear in the second.
    \item $\Hs$ denotes a (real or complex) Hilbert space.
    \item All algebras are assumed to be associative.
\end{itemize}
\end{tcolorbox}

\vspace*{20mm}

\begin{tcolorbox}[interior style={left color=abscol1!40!white,right color=abscol2!40!white},frame style={left color=abscol1!80!white,right color=abscol2!80!white},width=6in,rounded corners,center,top=2pt,bottom=2pt,left=2pt,right=2pt,drop fuzzy shadow,enhanced]
\begin{abstract}
    In this project, we will develop the theory of Banach algebras and prove two celebrated theorems, the Gelfand representation theorem and the GKZ theorem. We will then proceed to develop the theory of \cstas\ and prove the Gelfand-Naimark theorem and the GNS theorem. Finally, we will show how they can be used to give an alternative formulation of quantum mechanics.
\end{abstract}
\end{tcolorbox}

\vfill

\begin{tcolorbox}[interior style={left color=cdcol1!40,right color=cdcol2!40},frame style={left color=cdcol1!80!black,right color=cdcol2!80!black},fontupper=\sffamily,top=2mm,bottom=2mm,left=2mm,right=2mm,drop fuzzy shadow,enhanced,breakable]
The cover image was created by the author using {\rmfamily Ti\textit{k}Z} and {\rmfamily\scshape PGF}. It depicts a set of quantum particles (or equivalently, wave packets) moving around in an oscillating potential, interfering with it as they vibrate.
\end{tcolorbox}

\section{Introduction}\label[chapter]{introduction}
\ind When quantum mechanics was first developed in the 1920s by \href{https://mathshistory.st-andrews.ac.uk/Biographies/Schrodinger/}{Erwin Schrödinger} (1887--1961), \href{https://mathshistory.st-andrews.ac.uk/Biographies/Heisenberg/}{Werner Heisenberg} (1901--1976), \href{https://mathshistory.st-andrews.ac.uk/Biographies/Pauli/}{Wolfgang Pauli} (1900--1958) and \href{https://mathshistory.st-andrews.ac.uk/Biographies/Born/}{Max Born}\footnote{Outside of the physics community, Born is more famous as the grandfather of singer and actress Olivia Newton-John.} (1882--1970) (among many others), it was formulated in terms of matrix mechanics, differential equations and calculus of variations. While these are sufficient in practice for modeling simple quantum systems (e.g. a particle in a potential well), they are not mathematically rigorous, and this leads to problems when trying to generalize them to more complicated (e.g. infinite-dimensional) systems.

\begin{center}
\begin{minipage}{24em}\begin{center}
\personbox{Periwinkle!50}{\href{https://mathshistory.st-andrews.ac.uk/Biographies/Dirac/}{Paul Dirac}\\(1902--1984)}{\includegraphics[width=1.3in]{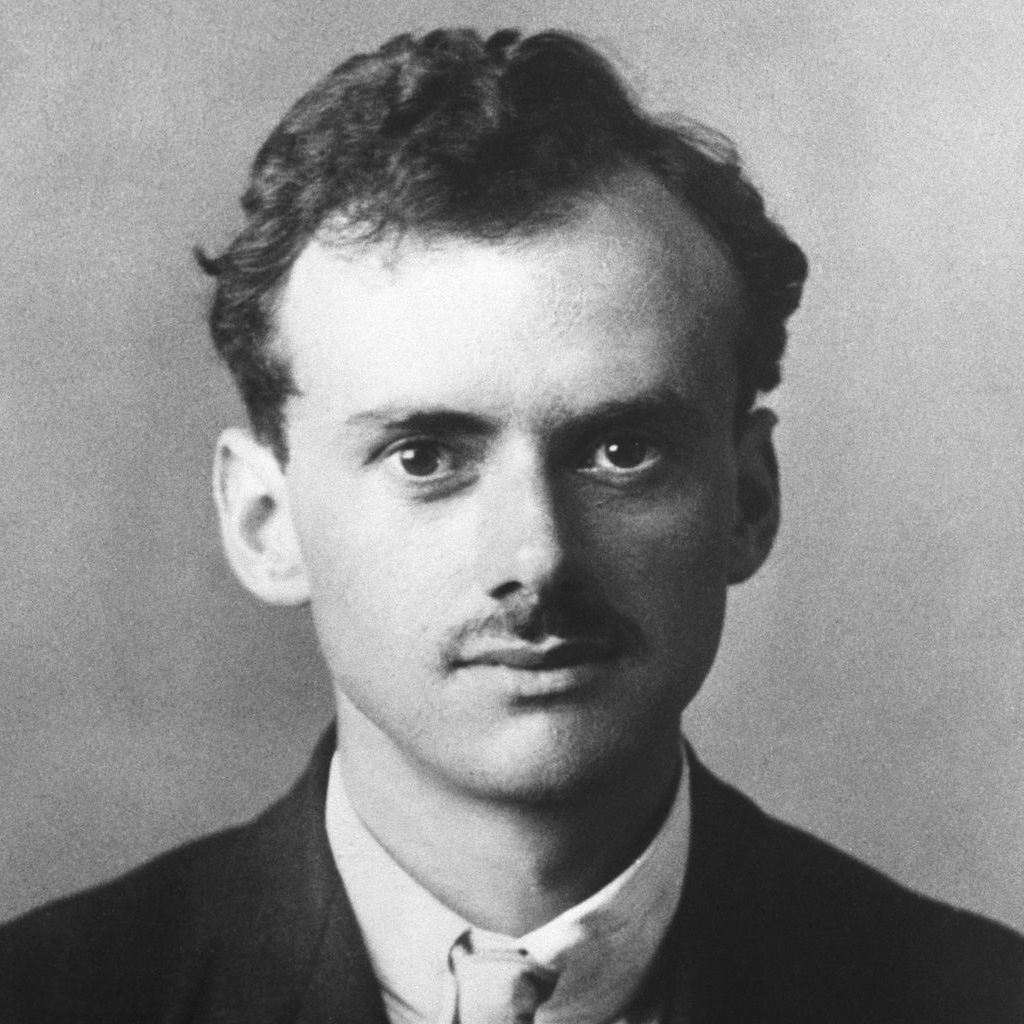}}\label{xdirac} 
\end{center}\end{minipage}
\begin{minipage}{24em}\begin{center}
\personbox{Mulberry!50}{\href{https://mathshistory.st-andrews.ac.uk/Biographies/Von_Neumann/}{John von Neumann}\\(1903--1957)}{\includegraphics[width=1.3in]{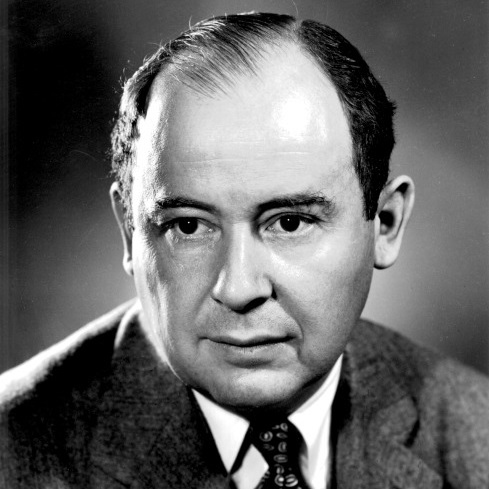}}\label{xvonneumann} 
\end{center}\end{minipage}
\end{center}

\ind The first rigorous mathematical formulation of quantum mechanics was given in 1930 by Dirac \cite{dirac} and in 1932 by von Neumann\footnote{Not to be confused with \href{https://mathshistory.st-andrews.ac.uk/Biographies/Neumann_Carl/}{Carl Neumann} (1832--1925), known for \hyperref[neumann]{Neumann series} and Neumann boundary conditions, or \href{https://en.wikipedia.org/wiki/Ezra_T._Newman}{Ezra T. Newman} (1929--), known for the Newman-Penrose (NP) formalism and the Kerr-Newman metric.} \cite{neumann}. This is primarily based on functional analysis in Hilbert spaces, and boils down to a set of three axioms, now known as the \emph{Dirac-von Neumann axioms}:

\begin{tcolorbox}[title=Dirac-von Neumann Axioms (Hilbert space version), fonttitle=\sffamily\bfseries\selectfont,colback=Green!10,colframe=Green!40,coltitle=black,top=2mm,bottom=2mm,left=2mm,right=2mm,drop fuzzy shadow,before skip=10pt,after skip=10pt,enhanced]\label{diracvn}
Every quantum system can be modeled using a complex Hilbert space $\Hs$.
\begin{enumerate}
    \item The observables of the system are the self-adjoint operators $A:\Hs\to\Hs$.
    \item The states of the system are the unit vectors in $\Hs$.
    \item The expected value of an observable $A$ when the system is in a state $\psi$ is\footnotemark\ $\ip{A\psi,\psi}$.
\end{enumerate}
\end{tcolorbox}
\footnotetext{Most texts on quantum mechanics define this as $\ip{\psi,A\psi}$ instead. This is because they assume inner products to be linear in the second argument, rather than the first as we do here.}

\ind Dirac and von Neumann's formulation of quantum mechanics starts with a Hilbert space of states and builds the observables from there (using operators). But we can do this the other way around. This leads to an alternative description of the Dirac-von Neumann axioms, using what are now known as \emph{\cstas}:

\begin{tcolorbox}[title=Dirac-von Neumann Axioms (\csta\ version), fonttitle=\sffamily\bfseries\selectfont,colback=Green!10,colframe=Green!40,coltitle=black,top=1mm,bottom=1mm,left=1mm,right=1mm,drop fuzzy shadow,before skip=10pt,after skip=10pt,enhanced]
Every quantum system can be modeled using a complex \csta\ $\As$.
\begin{enumerate}
    \item The observables of the system are the self-adjoint elements of $\As$.
    \item The states of the system are the unit positive linear functionals $\omega:\As\to\C$.
    \item The expected value of an observable $a$ when the system is in a state $\omega$ is $\omega(a)$.
\end{enumerate}
\end{tcolorbox}

\ind In this formulation, we start with a \csta\ of observables and build the states from there (using linear functionals). This time, it is the \emph{states} that act on the \emph{observables}, rather than the \emph{observables} acting on the \emph{states}. The fact that this is compatible with the Hilbert space formulation is a consequence of the \hyperref[gns]{GNS theorem}, which we will prove in \Cref{cstas}.\\

\ind The move to consolidate the mathematics of quantum mechanics into an abstract theory was pioneered in 1933 by \href{https://mathshistory.st-andrews.ac.uk/Biographies/Jordan_Pascual/}{Pascual Jordan}\footnote{Not to be confused with \href{https://en.wikipedia.org/wiki/Wilhelm_Jordan_(geodesist)}{Wilhelm Jordan} (1842--1899), known for Gauss-Jordan elimination, or \href{https://mathshistory.st-andrews.ac.uk/Biographies/Jordan/}{Camille Jordan} (1838--1922), known for the Jordan curve theorem and the Jordan normal form.} (1902--1980), though the abstract theory of \cstas\ did not come until 1946. Meanwhile, other formulations of quantum mechanics were also under way, such as the phase space formulation (developed in 1946 by \href{https://en.wikipedia.org/wiki/Hilbrand_J._Groenewold}{Hilbrand Groenewold} (1910--1996)) and the path integral formulation (developed in 1948 by \href{https://mathshistory.st-andrews.ac.uk/Biographies/Feynman/}{Richard Feynman} (1918--1988)).

\newpage 

\subsection{The Name of the Game}
We first introduce the two most important and fundamental examples, the Batman and Robin\footnote{or Bonnie and Clyde, or Thelma and Louise, or whichever dynamic duo you prefer.}  of \cstas:

{
\definecolor{col}{HTML}{872143}
\begin{tcolorbox}[title=$C_0(X)$: Continuous functions vanishing at infinity on a locally compact Hausdorff space, fonttitle=\sffamily\bfseries\selectfont,colback=col!10,colframe=col!40,coltitle=black,top=1mm,bottom=1mm,left=1mm,right=1mm,drop fuzzy shadow,before skip=10pt,after skip=10pt,enhanced]
Suppose $X$ is a locally compact Hausdorff space. Then $C_0(X)$ is the space of all functions $f:X\to\F$ such that:
\begin{itemize}[align=left,leftmargin=0.5in,font=\sffamily\itshape\selectfont]
    \item[Continuous:] For any open set $U\sub\F$, the pre-image $f^{-1}(U)=\{x\in X\mid f(x)\in U\}$ is open in $X$.
    \item[Vanishing at infinity:] For any $\ep>0$, there is a compact set $K\sub X$ such that $\abs{f(x)}<\ep$ for all $x\in X\exc K$.
\end{itemize}
Elements of $C_0(X)$ can be added, scaled and multiplied in the usual way:
\begin{align*}
    (f+g)(x)=f(x)+g(x) && (\alpha f)(x)=\alpha f(x) && (fg)(x)=f(x)g(x)
\end{align*}
We also define the following operations on $C_0(X)$:
\vspace*{-8pt} 
\begin{multicols}{2}
\begin{itemize}[align=left,leftmargin=0.5in,font=\sffamily\itshape\selectfont]
    \item[Adjoint:] $f^*:X\to\F$, $f^*(x)=\overline{f(x)}$
\end{itemize}
\begin{itemize}[align=left,leftmargin=0.5in,font=\sffamily\itshape\selectfont]
    \item[Norm:] $\Disp\norm{f}=\sup_{x\in X} \abs{f(x)}$
\end{itemize}
\end{multicols}
\vspace*{-10pt} 
\end{tcolorbox}
}

\begin{example}
If $X=\R^n$, then $C_0(X)$ is the set of all continuous functions on $\R^n$ that vanish as $\abs{\xb}\to\infty$, i.e. those that decay to zero (hence the name \emph{vanishing at infinity}).
\end{example}

For any locally compact Hausdorff space $X$, we can define three other related spaces:
\begin{itemize}
    \item $C(X)$, the space of \emph{all} continuous functions $f:X\to\F$.
    \item $C_b(X)$, the space of all \emph{bounded} continuous functions $f:X\to\F$.
    \item $C_c(X)$, the space of all continuous functions $f:X\to\F$ with \emph{compact support}, i.e. there is a compact set $K\sub X$ such that $f(x)=0$ for all $x\in X\exc K$.
\end{itemize}
These spaces are nested as follows:
\begin{equation*}
    C_c(X)\sub C_0(X)\sub C_b(X)\sub C(X)
\end{equation*}
All of these inclusions are strict unless $X$ is compact, in which case all four spaces are identical.

{
\definecolor{col}{HTML}{1232A4}
\begin{tcolorbox}[title=$\Bs(\Hs)$: Bounded linear operators on a Hilbert space, fonttitle=\sffamily\bfseries\selectfont,colback=col!10,colframe=col!40,coltitle=black,top=1mm,bottom=1mm,left=1mm,right=1mm,drop fuzzy shadow,before skip=10pt,after skip=10pt,enhanced]
Suppose $\Hs$ is a Hilbert space. Then $\Bs(\Hs)$ is the space of all functions $T:\Hs\to\Hs$ such that:
\begin{itemize}[align=left,leftmargin=0.5in,font=\sffamily\itshape\selectfont]
    \item[Linear:] For any $x,y\in\Hs$ and any $\alpha\in\F$, we have $T(x+y)=T(x)+T(y)$ and $T(\alpha x)=\alpha T(x)$.
    \item[Bounded:] There exists $M\ge0$ such that $\norm{Tx}\le M\norm{x}$ for all $x\in\Hs$.
\end{itemize}
Elements of $\Bs(\Hs)$ can be added, scaled and multiplied in the usual way:
\begin{align*}
    (S+T)(x)=S(x)+T(x) && (\alpha T)(x)=\alpha T(x) && (ST)(x)=S(T(x))
\end{align*}
We also define the following operations on $\Bs(\Hs)$:
\vspace*{-8pt} 
\begin{multicols}{2}
\begin{itemize}[align=left,leftmargin=0.5in,font=\sffamily\itshape\selectfont]
    \item[Adjoint:] $T^*:\Hs\to\Hs$, $\ip{T^*(x),y}=\ip{x,T(y)}$
\end{itemize}
\begin{itemize}[align=left,leftmargin=0.5in,font=\sffamily\itshape\selectfont]
    \item[Norm:] $\Disp\norm{T}=\sup_{x\in\Hs,\norm{x}_\Hs\le1} \norm{T(x)}_\Hs$
\end{itemize}
\end{multicols}
\vspace*{-10pt} 
\end{tcolorbox}
}

\begin{examples}\leavevmode
\begin{enumerate}
    \item If $\Hs=\F^n$, then $\Bs(\Hs)=M_n(\F)$, the set of all $n\times n$ matrices with entries in $\F$. If we use the standard (Euclidean) inner product on $\F^n$, then the adjoint of any matrix $A\in M_n(\F)$ is simply its conjugate transpose (if $\F=\R$, this is even more simply its transpose). This \csta\ can be used to describe a single quantum particle with a finite number of eigenstates, e.g. a stationary particle with spin.
    \item If $\Hs=\ell^2=\left\{(x_n)_{n=1}^\infty\,\middle|\,\sumn\abs{x_n}^2<\infty\right\}$ (the set of all square-summable sequences in $\F$), then $\Bs(\Hs)$ is the set of all bounded operators $T:\ell^2\to\ell^2$. This \csta\ can be used to describe a single quantum particle with a countably infinite number of eigenstates, e.g. a particle in a box.
\end{enumerate}
\end{examples}

\ind The key difference between the examples $C_0(X)$ and $\Bs(\Hs)$ is that the former is \emph{abelian}, i.e. $f(x)g(x)=g(x)f(x)$, while the latter is \emph{non-abelian}, i.e. $S\circ T\ne T\circ S$ (in general). As we will see later, this will lead us to formulate two different theorems, one for abelian \cstas\ (the \hyperref[gelfandnaimark]{Gelfand-Naimark theorem}) and the other for non-abelian \cstas\ (the \hyperref[gns]{GNS theorem}\footnote{Technically, the GNS theorem is also applicable to abelian \cstas, but in this case, the Gelfand-Naimark theorem is simpler.}). Together, these will allow us to characterize \emph{all} \cstas.

\newpage 

\subsection{What This Project is About}
\ind With the end (quantum mechanics) in mind, let us outline the journey to get there. We will combine tools from algebra and analysis with concepts from theoretical physics to develop a theory that is mathematically rigorous, yet also suitable to describe real-world physical systems. Before the 1930s, this was a delicate balance that was rarely emphasized, let alone achieved\footnote{Perhaps the closest approximation of this is when Einstein first formulated his theory of general relativity in 1915, using tools from tensor calculus and differential geometry of manifolds, developed merely 20 years earlier by Ricci and Levi-Civita.}. As the 20\textsuperscript{th} century progressed, there were several major developments that helped strike this balance, such as topological manifolds (used in string theory and supergravity), compact Lie algebras (used in Yang-Mills theory) and, as we will discuss here, \cstas.
\begin{center}
\begin{tikzpicture}[scale=0.7,every node/.style={transform shape,align=center,font=\Large\sffamily}]
\draw[fill=red,fill opacity=0.3] (-150:{sqrt(3)}) circle (3);
\draw[fill=green,fill opacity=0.3] (90:{sqrt(3)}) circle (3);
\draw[fill=blue,fill opacity=0.3] (-30:{sqrt(3)}) circle (3);
\node[scale=1.1] at (-150:{2*sqrt(3)}) {Algebra};
\node[scale=1.1] at (90:{2*sqrt(3)}) {Analysis};
\node at (-30:{2*sqrt(3)}) {Theoretical\\Physics};
\node[xshift=-2mm,yshift=1mm] at (150:2) {Lie\\theory};
\node[yshift=-1mm] at (-90:2) {Particle\\physics};
\node[xshift=1mm,yshift=3mm] at (30:2) {Relativity};
\node[scale=1.3] at (0:0) {\bfseries This\\\bfseries project}; 
\end{tikzpicture}\\
This diagram is not exhaustive, there are \emph{many} other topics in the overlapping regions.
\end{center}

\vspace{10mm}

\ind In \Cref{banachalgs}, we will build the theory of Banach algebras and use the combined power of algebra and analysis to prove some powerful results in spectral theory and representation theory. In particular, we will prove the \hyperref[gelfandrep]{Gelfand representation theorem} and the \hyperref[gkz]{GKZ theorem}.\\

\ind In \Cref{cstas}, we will strengthen these ideas and construct the theory of \cstas, and show how just one axiom combining the algebraic and analytic sides of the game leads to such a rich and vivid structure. We will then use these to prove the \hyperref[gelfandnaimark]{Gelfand-Naimark theorem} and the \hyperref[gns]{GNS theorem}, thereby explicitly classifying all \cstas.\\

\ind Finally, in \Cref{rebuildingqm}, we will recast the Dirac-von Neumann axioms in terms of \cstas\ and discuss how they can be used to model quantum mechanical systems.\\

\ind Now that we have an idea of the journey ahead, let's take a trip down the (quantum) path less traveled...

\vfill

\subsection*{Acknowledgments}
\ind I would like to thank my project supervisor, Prof. Benjamin Doyon, for his invaluable feedback and comments on my project. I would also like to thank Dr. Paul Cook for his unwavering support in my learning endeavors and his encouragement for me to embark on mathematical research from early on, some of which has led to this project. Finally, I would like to thank my peers for their support throughout my years at King's College London. In particular, I would like to thank Tamanna Sehgal and John Campbell for their help and feedback on the format and aesthetics of this project.

\vspace{8mm}

\begin{flushright} 
Senan Sekhon\\
March 31, 2021
\end{flushright}

\section{Banach Algebras}\label[chapter]{banachalgs}
\ind Before we can talk about \cstas, we first need to introduce the concept of \emph{Banach algebras}. Banach algebras are where algebra meets analysis. They have a norm structure (which allows us to talk about things like distances, convergence and continuity), as well as an algebra structure (which allows us to talk about things like polynomials, ideals and homomorphisms). Combining these two structures opens up more concepts for us, such as power series and spectra, which are fundamental to the theory of functional calculus.\\

\ind While these concepts will play a vital role in our study of \cstas\ and quantum mechanics, they are also interesting on their own. As such, we will spend some time investigating the general theory of Banach algebras before we specialize to \cstas.\\

\ind For further reading on Banach algebras, see \cite{muscat}, \cite{ward}, \cite{kaniuth}, \cite{douglas}, \cite{conway2}, \cite{grudin} and \cite{rickart}.

\subsection{Banach Algebras}
\begin{definition}\label{algebradefn}
An \textbf{algebra} over a field $\F$ is a vector space $\As$ over $\F$, together with a binary operation $\cdot:\As\times\As\to\As$ (known as \textbf{multiplication} or \textit{vector multiplication}), that satisfies the following:
\begin{enumerate}
    \item If $x,y\in\As$, then $x\cdot y\in\As$. \hfill (Closure)
    \item If $x,y,z\in\As$, then $x\cdot(y+z)=x\cdot y+x\cdot z$. \hfill (Left distributivity)
    \item If $x,y,z\in\As$, then $(x+y)\cdot z=x\cdot z+y\cdot z$. \hfill (Right distributivity)
    \item If $x,y\in\As$ and $\alpha\in\F$, then $(\alpha x)\cdot y=\alpha(x\cdot y)=x\cdot(\alpha y)$. \hfill (Compatibility with scalar multiplication)
    \item If $x,y,z\in\As$, then $x\cdot(y\cdot z)=(x\cdot y)\cdot z$. \hfill (Associativity)
\end{enumerate}
\end{definition}

\begin{examples}\leavevmode
\begin{enumerate}
    \item Every field $\F$ is an algebra over itself, with the vector multiplication operation given by the field multiplication.
    \item $\F^n$ is an algebra over $\F$, with componentwise addition, scalar multiplication and vector multiplication, e.g. $(a_1,a_2,...,a_n)\cdot(b_1,b_2,...,b_n)=(a_1b_1,a_2b_2,...,a_nb_n)$.
    \item Suppose $X$ is a set. Then the set $\F^X$ of all functions $f:X\to\F$ is an algebra, with pointwise addition, scalar multiplication and vector multiplication, e.g. $(f\cdot g)(x)=f(x)g(x)$.
    \item Suppose $X$ is a vector space. Then the set $\Ls(X)$ of all linear maps $T:X\to X$ is an algebra, with the multiplication operation given by composition.
    \item The set $M_n(\F)$ of all $n\times n$ matrices is an algebra, with matrix addition, scalar multiplication and matrix multiplication.
\end{enumerate}
\end{examples}

\begin{remark}
Some sources do not require algebras to be associative. Indeed, there are non-associative algebras, such as $\R^3$ with the cross product. Here, we will require all algebras to be associative. In particular, this makes them rings under the addition and multiplication operations. We will also usually omit the symbol for vector multiplication, i.e. write $ab$ instead of $a\cdot b$.
\end{remark}

\begin{definition}
An algebra $\As$ is \textbf{abelian} (or \textit{commutative}) if $a\cdot b=b\cdot a$ for all $a,b\in\As$.
\end{definition}
\begin{examples}\leavevmode
\begin{enumerate}
    \item $\F^n$ (with componentwise multiplication) is abelian.
    \item $\F^X$ (with pointwise multiplication) is abelian.
    \item $M_n(\F)$ is non-abelian if $n\ge2$.
    \item If $X$ is a vector space and $\dim(X)\ge2$, then $\Ls(X)$ is non-abelian.
\end{enumerate}
\end{examples}

\begin{definition}
Suppose $\As$ is an algebra. A \textbf{subalgebra} of $\As$ is a subset $\Bs\sub\As$ that is also an algebra (with the same operations as $\As$).
\end{definition}
\begin{remark}
Equivalently, a subalgebra of $\As$ is a (vector) subspace $\Bs$ of $\As$ that is closed under vector multiplication, i.e. for all $a,b\in\Bs$, we have $a\cdot b\in\Bs$.
\end{remark}

\ind The intersection of any collection of subalgebras of $\As$ is a subalgebra of $\As$. The proof of this is straightforward from the definition. An important consequence is that every subset of $\As$ ``generates'' a subalgebra of $\As$, the intersection of all subalgebras of $\As$ containing it. This is analogous to the concept of span in linear algebra.\\

\ind We can easily construct larger algebras from smaller ones. If $\As$ and $\Bs$ are algebras over $\F$, then the Cartesian product $\As\times\Bs$ is an algebra over $\F$, with coordinatewise addition, scalar multiplication and vector multiplication:
\begin{align*}
    (a_1,b_1)+(a_2,b_2)=(a_1+b_1,a_2+b_2) && \lam(a,b)=(\lam a,\lam b) && (a_1,b_1)\cdot(a_2,b_2)=(a_1\cdot b_1,a_2\cdot b_2)
\end{align*}

\begin{definition}
Suppose $\As$ is an algebra. An \textbf{identity} (or \textit{identity element} or \textit{unit}) of $\As$ is an element $\1\in\As$, $\1\ne0$ such that $\1a=a\1=a$ for all $a\in\As$.\\
An algebra $A$ is \textbf{unital} if it has an identity.
\end{definition}
\begin{warning}
This is NOT the same as a unit in a ring (where it refers to any element with a multiplicative inverse). For example, both $1$ and $-1$ are units in $\Z$, but only $1$ is an identity element.
\end{warning}
\begin{remark}
The assumption $\1\ne0$ is necessary to exclude the case where $\As=\{0\}$ (although $0\cdot a=a\cdot 0=0$ for all $a\in\{0\}$, we do not consider $0$ an identity).
\end{remark}

\begin{examples}\leavevmode
\begin{enumerate}
    \item $\F^n$ is unital, its identity is $(1,1,...,1)$.
    \item $M_n(\F)$ is unital, its identity is $I_n$ (the $n\times n$ identity matrix).
    \item The set $\As=\left\{[a_{ij}]\in M_n(\F)\mid a_{ij}=0\text{ for all }i>j\right\}$ of all upper triangular $n\times n$ matrices is a subalgebra of $M_n(\F)$. It is unital as $I_n\in\As$.
    \item The set $\As=\left\{[a_{ij}]\in M_n(\F)\mid a_{ij}=0\text{ for all }i\ge j\right\}$ of all strictly upper triangular $n\times n$ matrices is a subalgebra of $M_n(\F)$. It is not unital as $I_n\notin\As$.
\end{enumerate}
\end{examples}

\begin{proposition}
The identity of any algebra is unique.
\end{proposition}
\begin{proof}
Suppose $\1_1$ and $\1_2$ are two identities of $\As$. Since $\1_2$ is an identity, we have $\1_1\1_2=\1_1$, and since $\1_1$ is an identity, we have $\1_1\1_2=\1_2$. Thus $\1_1=\1_2$.
\end{proof}

\begin{definition}
Suppose $\As$ is a unital algebra. An element $a\in\As$ is \textbf{invertible} if there is an element $b\in\As$ such that $ab=ba=\1$. If so, we call $b$ the \textbf{inverse} of $a$ and denote it by $a^{-1}$.\\
The set of all invertible elements of $\As$ is denoted by $\As\inv$.
\end{definition}
\begin{remark}
It is enough that there are elements $b,c\in\As$ such that $ab=\1$ and $ca=\1$, as this implies that $b=c$ (since $b=\1b=(ca)b=c(ab)=c\1=c$).
\end{remark}
\begin{warning}
It is NOT enough that there is an element $b\in\As$ such that $ab=\1$, see the example after the next lemma. 
\end{warning}

\begin{examples}\leavevmode
\begin{enumerate}
    \item Suppose $\As=\F^n$ with componentwise multiplication. Then $\As\inv$ is the set of all elements of $\F^n$ with all components nonzero.
    \item Suppose $\As=M_n(\F)$. Then $\As\inv=GL_n(\F)$, the set of all $n\times n$ matrices with nonzero determinant.
    \item Suppose $X$ is a locally compact Hausdorff space and $f\in C_0(X)$. Then $f$ is invertible if and only if it does not vanish, i.e. $f(x)\ne0$ for all $x\in X$.
\end{enumerate}
\end{examples}

\begin{lemma}
Suppose $\As$ is a unital algebra and $a,b\in\As$. If $ab$ and $ba$ are invertible, then $a$ and $b$ are also invertible.
\end{lemma}
\begin{proof}
Since $ab(ab)^{-1}=\1$ and $(ba)^{-1}ba=\1$, it follows that $a$ is invertible and $a^{-1}=b(ab)^{-1}=(ba)^{-1}b$. Switching $a$ and $b$ shows that $b$ is also invertible.
\end{proof}
\begin{warning}
It is NOT sufficient that $ab$ is invertible!
\end{warning}
\begin{example}
Suppose $\As=\Bs(\ell^\infty)$ and define the \textit{left shift} and \textit{right shift} operators $L,R\in\As$ respectively by $L((x_1,x_2,x_3,...))=(x_2,x_3,...)$ and $R((x_1,x_2,x_3,...))=(0,x_1,x_2,...)$. Then $LR=\1$ is invertible, but $L$, $R$ and $RL$ are not. This is because $L$ is surjective but not injective, and vice versa for $R$.
\end{example}

\begin{definition}
An \textbf{algebra norm} (or \textit{norm}) on an algebra $\As$ is a vector space norm $\norm{\cdot}:\As\to\R$ such that for all $x,y\in\As$, we have $\norm{x\cdot y}\le\norm{x}\norm{y}$.\\
An algebra with a norm is known as a \textbf{normed algebra}.
\end{definition}
\begin{remark}
Since every normed algebra is (by definition) a normed vector space, it is also a metric space with respect to the induced metric $d(x,y)=\norm{x-y}$.
\end{remark}

\begin{examples}\leavevmode
\begin{enumerate}
    \item $\F^n$ is a normed algebra with the $p$-norm $\norm{\xb}_p=(\sumkn\abs{x_k}^p)^{1/p}$ for any $1\le p<\infty$ (if $p=2$, this is simply the Euclidean norm $\norm{\xb}_2=\sqrt{\sumkn\abs{x_k}^2}$). It is also a normed algebra with the maximum norm $\norm{\xb}_\infty=\max_{1\le k\le n} \abs{x_k}$.
    \item $M_n(\F)$ is a normed algebra with the spectral norm $\norm{A}=\sup_{\xb\in\F^n,\norm{\xb}_2\le1} \norm{A\xb}_2$. It is also a normed algebra with the Frobenius norm $\norm{A}_F=\sqrt{\sumkn\sum_{l=1}^n \abs{a_{kl}}^2}$.
    \item Suppose $X$ is a normed vector space. Then the set $\Bs(X)$ of all bounded linear operators $T:X\to X$ is a normed algebra (with the operator norm).
    \item Suppose $X$ is a locally compact Hausdorff space. Then $C_0(X)$ is a normed algebra with the \emph{supremum norm} $\norm{f}_\infty=\sup_{x\in X} \abs{f(x)}$ and the pointwise product $(fg)(x)=f(x)g(x)$.
    \item Suppose $(X,\As,\mu)$ is a measure space. Then the Lebesgue space $L^\infty(X,\As,\mu)$ is a normed algebra with the \emph{essential supremum} norm:
    \begin{equation*}
        \norm{f}_\infty=\inf\left\{M\ge0\mid\mu(\{x\in X\mid\abs{f(x)}>M\})=0\right\}
    \end{equation*}
    And the pointwise product. In particular, we have the following important special cases:
    \begin{enumerate}
        \item If $X=\N$, $\As=\Ps(\N)$ and $\mu$ is the counting measure, we get the sequence space $\ell^\infty$.
        \item If $X=\R$, $\As$ is the Lebesgue $\sigma$-algebra on $\R$ and $\mu$ is the Lebesgue measure, we get the Lebesgue space $L^\infty(\R)$.
    \end{enumerate}
\end{enumerate}
\end{examples}

Two notable examples of normed algebras arise from the concept of \emph{convolution}:

\begin{example}
The sequence space $\ell^1(\Z)$ is a normed algebra with the pointwise product $(a\cdot b)_n=a_nb_n$. It is also a (different) normed algebra with the convolution product:
\begin{center}
\begin{tikzpicture}[scale=0.5]
\draw[<->,thick] (-9.25,0)--(9.25,0) node[right]{$n$}; 
\foreach \a in {-2,-1,0,1,2}
\fill[red] (\a,1) circle (3pt);
\foreach \a in {-9,-8,...,-3,3,4,...,9}
\fill[red] (\a,0) circle (3pt);
\foreach \n in {-9,-8,...,9}
\fill[blue] (\n,{2^(-abs(\n))}) circle (3pt);
\foreach \na in {-9,-8,...,-2,2,3,...,9}
\fill[red!50!blue] (\na,{31*2^(-abs(\na)-2)}) circle (3pt);
\foreach \pt in {(-1,19/8),(0,5/2),(1,19/8)}
\fill[red!50!blue] \pt circle (3pt);
\node[red,above] at (-1,1) {$a_n$};
\node[blue,left] at (1,0.3) {$b_n$};
\node[red!50!blue,above right] at (3,1) {$a_n*b_n$};
\node at (-16,1) {$\Disp(a*b)_n=\sum_{m=-\infty}^\infty a_mb_{n-m}$};
\end{tikzpicture}
\end{center}
More generally, suppose $G$ is a group. Then the space $\ell^1(G)$ of all absolutely summable functions $f:G\to\F$ is a normed algebra with the convolution product:
\begin{equation*}
    (f*g)(x)=\sum_{y\in G} f(y)g(y^{-1}x)
\end{equation*}
This is known as the \emph{discrete group algebra} of $G$.
\end{example}

\begin{example}
The Lebesgue space $L^1(\R)$ is NOT a normed algebra with the pointwise product $(fg)(x)=f(x)g(x)$, as the product of two functions in $L^1(\R)$ may not be in $L^1(\R)$. However, it is a normed algebra with the convolution product:
\begin{center}
\begin{tikzpicture}[scale=0.5]
\draw[<->,thick] (-9.15,0)--(9.15,0) node[right]{$x$}; 
\draw[-,very thick,red] (-9,0)--(-2,0) (-2,1)--(2,1) (2,0)--(9,0);
\draw[very thick,blue,smooth,domain=-9:9,samples=200] plot(\x,{1/((\x)^2+1)});
\draw[very thick,red!50!blue,smooth,domain=-9:9,samples=200] plot(\x,{pi/180*(atan(\x+2)-atan(\x-2))});
\node[red,above] at (-1,1) {$f$};
\node[blue,left] at (0.9,0.3) {$g$};
\node[red!50!blue,above right] at (2.7,0.6) {$f*g$};
\node at (-16,1) {$\Disp(f*g)(x)=\int_{-\infty}^\infty f(y)g(x-y)\dy$};
\end{tikzpicture}
\end{center}
More generally, suppose $G$ is a topological group, $\Bs(G)$ is the Borel $\sigma$-algebra on $G$ and $\mu$ is a left-invariant Haar measure on $(G,\Bs(G))$. Then the space $L^1(G,\Bs(G),\mu)$ of all absolutely integrable functions $f:G\to\F$ is a normed algebra with the convolution product:
\begin{equation*}
    (f*g)(x)=\int_G f(y)g(y^{-1}x)\,d\mu
\end{equation*}
This is known as the \emph{Haar measure algebra} of $(G,\Bs(G),\mu)$.
\end{example}

\begin{definition}\label{unitalnormedalgebra}
A normed algebra $\As$ is \textbf{unital} if it has an identity $\1\in\As$ such that $\norm{\1}=1$.
\end{definition}
\begin{remark}
A unital normed algebra is NOT simply a normed algebra with an identity, it has the additional requirement that $\norm{\1}=1$. This is not vacuous, as the next example shows.
\end{remark}

\begin{example}
$\R$ is a unital normed algebra with the norm $\norm{x}=\abs{x}$. However, it is not a unital normed algebra with the norm $\norm{x}=2\abs{x}$ since $\norm{1}=2\ne1$. (though it is still a normed algebra and a unital algebra).
\end{example}

\ind If a normed algebra $\As$ has an identity $\1$, it follows automatically that $\norm{\1}\ge1$, since $\norm{\1}=\norm{\1\cdot\1}\le\norm{\1}\norm{\1}$. In practice, we usually assume that $\norm{\1}=1$, as we can always replace the norm with an equivalent norm to make this true. See \cite[Proposition 1.1.1, Page 2]{kaniuth} for a proof.

\begin{example}
Suppose $X$ is a locally compact (but not compact) Hausdorff space. Then $C_0(X)$ is a non-unital normed algebra, since the identity element (the constant function $1$) is not in $C_0(X)$.
\end{example}

Some normed algebras do not have an identity, but they have an \emph{approximate identity}, i.e. a sequence $(e_n)_{n=1}^\infty$ such that $\limn\norm{e_na-a}=0$ for all $a\in\As$.

\begin{example}
Suppose $e$ is the identity element in $L^1(\R)$. Then for all $f\in L^1(\R)$, we must have $e*f=f$, i.e.
\begin{equation*}
    \int_{-\infty}^\infty e(x-y)f(y)\dy=f(x)
\end{equation*}
This is satisfied if and only if $e(x)=\delta(x)$, the Dirac delta function. This function is not in $L^1(\R)$, and so $L^1(\R)$ has no identity. However, we can define a sequence $(e_n)_{n=1}^\infty$ of functions in $L^1(\R)$ such that for any $f\in L^1(\R)$, we have $\norm{e_n*f-f}_1\to0$ as $n\to\infty$. Below are some examples of such a sequence:
\begin{center}
\begin{tikzpicture}
\clip (-3.5,-0.7) rectangle (11,3.2);
\begin{scope}[xshift=0cm] 
\draw[-] (-3,0)--(3,0) node[right]{$x$};
\draw[-] (0,0)--(0,3.2);
\draw[thick,red,smooth,domain=-3:3,samples=200] plot(\x,{1/2*exp(-1/2*(\x)^2)});
\draw[thick,Orange,smooth,domain=-3:3,samples=200] plot(\x,{2/2*exp(-2*(\x)^2)});
\draw[thick,Dandelion,smooth,domain=-3:3,samples=200] plot(\x,{3/2*exp(-9/2*(\x)^2)});
\draw[thick,Green,smooth,domain=-3:3,samples=200] plot(\x,{4/2*exp(-8*(\x)^2)});
\draw[thick,blue,smooth,domain=-3:3,samples=200] plot(\x,{5/2*exp(-25/2*(\x)^2)});
\draw[thick,violet,smooth,domain=-3:3,samples=200] plot(\x,{6/2*exp(-18*(\x)^2)});
\node[red] at (1.8,0.3) {$e_1$};
\node[Orange] at (0.7,0.7) {$e_2$};
\node[Dandelion] at (0.5,1.3) {$e_3$};
\node[Green] at (0.45,1.8) {$e_4$};
\node[blue] at (0.4,2.3) {$e_5$};
\node[violet] at (0.35,2.8) {$e_6$};
\node at (-2,2.4) {$\Disp e_n(x)=\frac{n}{\sqrt{\pi}}e^{-n^2x^2}$};
\end{scope}
\begin{scope}[xshift=8cm] 
\draw[-] (-3,0)--(3,0) node[right]{$x$};
\draw[-] (0,0)--(0,3.2);
\draw[thick,red,smooth,domain=-3:3,samples=200] plot(\x,{sin(deg(\x))/(2*\x)});
\draw[thick,Orange,smooth,domain=-3:3,samples=200] plot(\x,{sin(deg(2*\x))/(2*\x)});
\draw[thick,Dandelion,smooth,domain=-3:3,samples=200] plot(\x,{sin(deg(3*\x))/(2*\x)});
\draw[thick,Green,smooth,domain=-3:3,samples=200] plot(\x,{sin(deg(4*\x))/(2*\x)});
\draw[thick,blue,smooth,domain=-3:3,samples=200] plot(\x,{sin(deg(5*\x))/(2*\x)});
\draw[thick,violet,smooth,domain=-3:3,samples=200] plot(\x,{sin(deg(6*\x))/(2*\x)});
\node[red] at (1.8,0.45) {$e_1$};
\node[Orange] at (1,0.7) {$e_2$};
\node[Dandelion] at (0.6,1.3) {$e_3$};
\node[Green] at (0.5,1.8) {$e_4$};
\node[blue] at (0.45,2.3) {$e_5$};
\node[violet] at (0.4,2.8) {$e_6$};
\node at (-2,2.4) {$\Disp e_n(x)=\frac{\sin(nx)}{x}$};
\end{scope}
\end{tikzpicture}
\end{center}
Approximate identities of $L^1(\R)$, such as those above, are also known as \emph{summability kernels}. They are vital in Fourier analysis, harmonic analysis, PDEs and distribution theory.
\end{example}

Some normed algebras do not even have an approximate identity.

\begin{example}
Suppose $\As\ne\{0\}$ is a normed vector space and define the multiplication operation by $ab=0$ for all $a,b\in\As$. Then $\As$ has no approximate identity, since for all $a,b\in\As$, $a\ne0$, we have $\norm{ba-a}=\norm{0-a}=\norm{a}$, so this cannot approach $0$.
\end{example}

\begin{proposition}\label{multcontinuous}
Suppose $\As$ is a normed algebra. Then the addition, scalar multiplication and vector multiplication operations are continuous.
\end{proposition}
\begin{proof}
Suppose $a_0,b_0\in\As$, $\lam\in\F$ and $\ep>0$.
\begin{itemize}[align=left,leftmargin=0pt,itemindent=!,itemsep=5pt,font=\sffamily\selectfont]
    \item[Addition:] Set $\delta=\frac{\ep}{2}$. Then for all $a,b\in\As$ such that $\norm{a-a_0},\norm{b-b_0}<\delta$, we have:
    \begin{equation*}
        \norm{(a+b)-(a_0+b_0)}=\norm{(a-a_0)+(b-b_0)}\le\norm{a-a_0}+\norm{b-b_0}<\delta+\delta=2\delta=\ep
    \end{equation*}
    \item[Scalar multiplication:] Set $\delta=\frac{\ep}{\abs{\lam}+1}$. Then for all $a\in\As$ such that $\norm{a-a_0}<\delta$, we have:
    \begin{talign*}
        \norm{\lam a-\lam a_0}=\norm{\lam(a-a_0)}=\abs{\lam}\norm{a-a_0}<\abs{\lam}\delta=\abs{\lam}\frac{\ep}{\abs{\lam}+1}<\ep
    \end{talign*}
    \item[Vector multiplication:] Set $m=\frac{\norm{a_0}+\norm{b_0}}{2}$ and set $\delta=\sqrt{m^2+\ep}-m$. Note that $\delta>0$ and $\delta^2+2m\delta=\ep$. Then for all $a,b\in\As$ such that $\norm{a-a_0},\norm{b-b_0}<\delta$, we have:
    \begin{align*}
        \norm{ab-a_0b_0}&=\norm{ab-a_0b+a_0b-a_0b_0}=\norm{(a-a_0)b+a_0(b-b_0)}\le\norm{a-a_0}\norm{b}+\norm{a_0}\norm{b-b_0}\\
        &<\delta\norm{b}+\delta\norm{a_0}=\delta(\norm{b}+\norm{a_0})=\delta(\norm{b-b_0+b_0}+\norm{a_0})\le\delta(\norm{b-b_0}+\norm{b_0}+\norm{a_0})\\
        &\le\delta(\delta+2m)=\delta^2+2m\delta=\ep \qedhere
    \end{align*}
\end{itemize}
\end{proof}

\begin{definition}
A \textbf{Banach algebra} is a normed algebra $\As$ that is complete, i.e. every Cauchy sequence in $\As$ converges in $\As$.
\end{definition}

This simply means that, as a normed vector space, it is a Banach space, hence the name `Banach algebra'.

\begin{examples}\leavevmode
\begin{enumerate}
    \item $\F^n$ is a Banach algebra with any of the $p$-norms, as well as with the maximum norm. More generally, every finite-dimensional normed algebra is a Banach algebra.
    \item Suppose $X$ is a Banach space. Then $\Bs(X)$ is a Banach algebra.
    \item Suppose $X$ is a locally compact Hausdorff space. Then $C_b(X)$ and $C_0(X)$ are Banach algebras. However, $C_c(X)$ is not a Banach algebra (unless $X$ is compact), as it is incomplete.
    \item Suppose $(X,\As,\mu)$ is a measure space. Then $L^\infty(X,\As,\mu)$ is a Banach algebra (the fact that it is complete follows from the Riesz-Fischer theorem).
    \item The set $C^1[a,b]$ of all continuously differentiable functions $f:[a,b]\to\F$ is NOT a Banach algebra with the supremum norm $\norm{f}=\sup_{x\in[a,b]} \abs{f(x)}$ (as it is not complete). However, it is a Banach algebra with the norm $\norm{f}=\sup_{x\in[a,b]} \abs{f(x)}+\sup_{x\in[a,b]} \abs{f'(x)}$.
\end{enumerate}
\end{examples}

\begin{definition}
Suppose $\As$ is an algebra and $\Is$ is a subspace of $\As$.
\begin{itemize}
    \item $\Is$ is a \textbf{left ideal} of $\As$ if $ab\in\Is$ for all $a\in\As$ and all $b\in\Is$.
    \item $\Is$ is a \textbf{right ideal} of $\As$ if $ba\in\Is$ for all $a\in\As$ and all $b\in\Is$.
    \item $\Is$ is an \textbf{ideal} (or \textit{two-sided ideal}) of $\As$ if it is both a left ideal and a right ideal of $\As$.
\end{itemize}
\end{definition}

\begin{examples}\leavevmode
\begin{enumerate}
    \item Every algebra $\As$ has $\{0\}$ and $\As$ as ideals (these are the \emph{trivial} ideals of $\As$).
    \item For any $a\in\As$, the set $\As a=\{ba\mid b\in\As\}$ is a left ideal of $\As$. Similarly, the set $a\As=\{ab\mid b\in\As\}$ is a right ideal of $\As$.
    \item Suppose $X$ is a normed vector space. Then the set $\Ks(X)$ of all compact operators in $\Bs(X)$ is a (two-sided) ideal of $\Bs(X)$.
    \item Suppose $X$ is a locally compact Hausdorff space. Then $C_0(X)$ is a (two-sided) ideal of $C_b(X)$.
    \item $\left\{\mmat{a}{b}{0}{0}\,\middle|\, a,b\in\F\right\}$ is a right ideal of $M_2(\F)$, but not a left ideal. Similarly, $\left\{\mmat{a}{0}{b}{0}\,\middle|\, a,b\in\F\right\}$ is a left ideal of $M_2(\F)$, but not a right ideal.
\end{enumerate}
\end{examples}

Every left/right/two-sided ideal of $\As$ is a subalgebra of $\As$. The intersection of any collection of left (resp. right, two-sided) ideals of $\As$ is a left (resp. right, two-sided) ideal of $\As$. The proofs of these are straightforward from the definitions.

\begin{definition}
Suppose $\As$ and $\Bs$ are algebras. A \textbf{homomorphism} (or \textit{algebra homomorphism}) from $\As$ to $\Bs$ is a linear map $\phi:\As\to\Bs$ such that $\phi(ab)=\phi(a)\phi(b)$ for all $a,b\in\As$.\\
An \textbf{isomorphism} (or \textit{algebra isomorphism}) is a bijective homomorphism.
\end{definition}

Basically, a homomorphism is a linear map that is also multiplicative.

\begin{examples}\leavevmode
\begin{enumerate}
    \item The map $\phi:\F\to\F$, $\phi(z)=\zbar$ is an isomorphism (if $\F=\R$, then $\phi$ is the identity map).
    \item The map $\phi:M_n(\F)\to\F$, $\phi(A)=\det(A)$ is NOT a homomorphism, as it is not linear.
\end{enumerate}
\end{examples}

\begin{proposition}\label{imskeri}
Suppose $\As$ and $\Bs$ are algebras and $\phi:\As\to\Bs$ is a homomorphism. Then:
\begin{enumerate}
    \item $\im(\phi)$ is a subalgebra of $\Bs$.
    \item $\ker(\phi)$ is an ideal of $\As$.
\end{enumerate}
\end{proposition}
\begin{proof}
Since $\phi$ is a homomorphism, it is linear, so $\im(\phi)$ is a subspace of $\Bs$, while $\ker(\phi)$ is a subspace of $\As$.
\begin{enumerate}
    \item Suppose $x,y\in\im(\phi)$. Then we have $x=\phi(a)$ and $y=\phi(b)$ for some $a,b\in\As$. This yields $xy=\phi(a)\phi(b)=\phi(ab)$, so $xy\in\im(\phi)$. Thus $\im(\phi)$ is a subalgebra of $\Bs$.
    \item Suppose $a\in\As$ and $b\in\ker(\phi)$. Then $\phi(b)=0$. This yields $\phi(ab)=\phi(a)\phi(b)=\phi(a)0=0$ and $\phi(ba)=\phi(b)\phi(a)=0\phi(a)=0$, so $ab\in\ker(\phi)$. Thus $\ker(\phi)$ is an ideal of $\As$. \qedhere
\end{enumerate}
\end{proof}

\begin{definition}
Suppose $\As$ and $\Bs$ are unital algebras. A \textbf{unital homomorphism} (or \textit{unital algebra homomorphism}) from $\As$ to $\Bs$ is an algebra homomorphism $\phi:\As\to\Bs$ such that $\phi(\1_\As)=\1_\Bs$.\\
In other words, a unital homomorphism is a homomorphism that maps the identity to the identity.
\end{definition}

\begin{examples}\leavevmode
\begin{enumerate}
    \item If $\As$ is a subalgebra of $\Bs$, then the inclusion map $\iota:\As\to\Bs$, $\iota(a)=a$ is a unital homomorphism.
    \item The map $\phi:\F\to M_2(\F)$, $\phi(x)=\mmat{x}{0}{0}{x}$ is a unital homomorphism.
    \item The map $\phi:\F\to M_2(\F)$, $\phi(x)=\mmat{x}{0}{0}{0}$ is a non-unital homomorphism.
\end{enumerate}
\end{examples}

\xnewpage

\begin{proposition}\label{unithominv}
Suppose $\As$ and $\Bs$ are unital algebras and $\phi:\As\to\Bs$ is a unital homomorphism. Then for all $a\in\As\inv$, we have $\phi(a)\in\Bs\inv$ and $(\phi(a))^{-1}=\phi(a^{-1})$.
\end{proposition}
\begin{proof}
Suppose $a\in\As\inv$. Then $\1_\Bs=\phi(\1_\As)=\phi(aa^{-1})=\phi(a)\phi(a^{-1})$ and $\1_\Bs=\phi(\1_\As)=\phi(a^{-1}a)=\phi(a^{-1})\phi(a)$. Thus $\phi(a)\in\Bs\inv$ and $(\phi(a))^{-1}=\phi(a^{-1})$.
\end{proof}

\begin{definition}
Suppose $\As$ is an algebra. An ideal $\Is$ of $\As$ is \textbf{proper} if $\Is\ne\As$. It is \textbf{maximal} if it is a proper ideal of $\As$ and it is not contained in any other proper ideal of $\As$.
\end{definition}

If $\As$ is unital, then every proper ideal of $\As$ is contained in a maximal ideal of $\As$. This is a special case of Krull's theorem, which is equivalent to the axiom of choice.

\begin{example}
Suppose $X$ is a locally compact Hausdorff space and $a\in X$. Then $\Is=\{f\in C_0(X)\mid f(a)=0\}$ is a maximal ideal of $C_0(X)$. To see this, suppose $\Js\supset\Is$ is an ideal of $C_0(X)$ that strictly contains $\Is$. Then there exists $f\in\Js$ such that $f(a)=c\ne0$. Suppose $g\in C_0(X)$ and define $h=\frac{g(a)}{c}f-g$. Then $h\in C_0(X)$ as it is a linear combination of two functions in $C_0(X)$. We also have $h(a)=\frac{g(a)}{c}c-g(a)=0$, so $h\in\Is\subset\Js$. Thus $g=\frac{g(a)}{c}f-h$ is a linear combination of two functions in $\Js$, and so $g\in\Js$. Thus $\Js=C_0(X)$, and so $\Is$ is a maximal ideal of $C_0(X)$.
\end{example}
\begin{example}
Suppose $\eta\in\R$. Then $\Is=\left\{f\in L^1(\R)\,\middle|\,\fh(\eta)=0\right\}$ is a maximal ideal of $L^1(\R)$, where $\fh(\xi)=\int_{-\infty}^\infty e^{-2\pi i\xi x}f(x)\dx$ is the Fourier transform of $f$. To see this, we can use Fourier transforms, i.e. pass from the time domain to the frequency domain. This is an algebra homomorphism from $L^1(\R)$ into $C_0(\R)$ (as it is linear and turns convolution into pointwise multiplication), and so the previous example shows that $\Is$ is maximal.
\end{example}

\begin{lemma}
Suppose $\As$ is a unital algebra and $\Is$ is a proper ideal of $\As$. Then $\As\inv\cap\Is=\emp$, i.e. $\Is$ does not contain any invertible elements of $\As$.
\end{lemma}
\begin{proof}
Suppose $a\in\As\inv\cap\Is$. Since $a\in\As\inv$, it has an inverse $a^{-1}$. Since $a\in\Is$ and $\Is$ is an ideal of $\As$, we have $aa^{-1}=\1\in\Is$. Thus for any $b\in\As$, we have $b=b\1\in\Is$, and so $\Is=\As$. Therefore if $\Is\ne\As$, then $\As\inv\cap\Is=\emp$.
\end{proof}

\begin{proposition}
Suppose $\As$ is an algebra and $\Is$ is an ideal of $\As$. Then $\Is$ is maximal if and only if $\As/\Is$ is simple, i.e. it does not have a non-trivial ideal.
\end{proposition}
\begin{proof}
($\Rightarrow$) Suppose $\Is$ is maximal but $\Bs=\As/\Is$ has a non-trivial ideal $\Js\ne\{0\},\Bs$. Then the quotient maps $\phi:\As\to\Bs$ and $\psi:\Bs\to\Bs/\Js$ are nonzero homomorphisms. Thus $\psi\circ\phi$ is also a nonzero homomorphism, and so $\ker(\psi\circ\phi)$ is an ideal of $\As$ containing $\Is$. Since $\Js\ne\{0\}$, there exists $a\in\As$ such that $[a]\ne[0]$, i.e. $a+\Js\ne\Js$. Thus $a\in\ker(\psi\circ\phi)$ but $a\notin\Js$, which contradicts the maximality of $\Js$.\\

($\Leftarrow$) Suppose $\Is$ is not maximal. Then there is a proper ideal $\Ks\supset\Is$. Define $\Bs=\As/\Is$ and $\Js_\Bs=\{x+\Is\mid x\in\Ks\exc\Is\}$. Since $\Ks\exc\Is\ne\emp$, we have $\Js_\Bs\ne\{0\}$. Since $\Ks\ne\As$, we have $\Js_\Bs\ne\Bs$. Also, for all $[a]\in\Bs$ and all $[j]\in\Js_\Bs$, we have $[a][j]=[aj]=[j]$ and $[j][a]=[ja]=[j]$, since $ja\in\Ks$ and $aj\in\Ks$ as $\Ks$ is an ideal of $\As$. Thus $\Js_\Bs$ is a non-trivial ideal of $\Bs$.
\end{proof}

\ind We can easily construct larger normed algebras from smaller ones. If $\As$ and $\Bs$ are normed algebras, then $\As\times\Bs$ is a normed algebra, with coordinatewise operations and the norm given by:
\begin{equation*}
    \norm{(a,b)}=\max\left\{\norm{a}_\As,\norm{b}_\Bs\right\}
\end{equation*}
\ind Of course, there are other possible norms on $\As\times\Bs$, such as $\norm{(a,b)}_1=\norm{a}_\As+\norm{b}_\Bs$.\\

\ind We can also do this for infinite collections of algebras, though this requires more care. See \cite[\S26, Definition 10, Page 136]{bonsall} for details.

\begin{definition}
Suppose $\As$ is an algebra and $\Is$ is an ideal of $\As$. The \textbf{quotient algebra} of $\As$ \textbf{modulo} $\Is$, denoted by $\As/\Is$, is given by:
\begin{equation*}
    \As/\Is=\left\{[a]\mid a\in\As\right\}
\end{equation*}
Where $[a]=a+\Is=\left\{a+b\mid b\in\Is\right\}$. The addition, scalar multiplication and vector multiplication operations in $\As/\Is$ are given by:
\begin{align*}
    [a]+[b]=[a+b] && \lam[a]=[\lam a] && [a][b]=[ab]
\end{align*}
The \textbf{quotient map} from $\As$ to $\As/\Is$ is given by $q:\As\to\As/\Is$, $q(a)=[a]=a+\Is$.
\end{definition}

\begin{proposition}\label{quotientalg}
The quotient algebra $\As/\Is$ is an algebra, and the quotient map $q$ is an algebra homomorphism.
\end{proposition}
\begin{proof}
Clearly $\As/\Is$ is a vector space and $q$ is a linear map. We now show that $\As/\Is$ is an algebra:
\begin{enumerate}
    \item Suppose $a+\Is,b+\Is\in\As/\Is$. Then $a,b\in\As$, so $ab\in\As$, and so $(a+\Is)(b+\Is)=ab+\Is\in\As/\Is$.
    \item $[a]([b]+[c])=[a][b+c]=[a(b+c)]=[ab+ac]=[ab]+[ac]=[a][b]+[a][c]$
    \item $([a]+[b])[c]=[a+b][c]=[(a+b)][c]=[ac+bc]=[ac]+[bc]=[a][c]+[b][c]$
    \item $(\lam[a])[b]=[\lam a][b]=[(\lam a)b]=[\lam(ab)]=\lam[ab]$ and $[a](\lam[b])=[a][\lam b]=[a(\lam b)]=[\lam(ab)]=\lam[ab]$
    \item $[a]([b][c])=[a][bc]=[a(bc)]=[(ab)c]=[ab][c]=([a][b])[c]$
\end{enumerate}
Thus $\As/\Is$ is an algebra. We also have $q(ab)=[ab]=[a][b]=q(a)q(b)$, so $q$ is an algebra homomorphism.
\end{proof}

\begin{theorem}[Neumann Series]\label{neumann}
Suppose $\As$ is a unital Banach algebra, $a\in\As$ and $\norm{a}<1$. Then $\1-a\in\As\inv$ and $(\1-a)^{-1}=\sumn[0] a^n$.
\end{theorem}
\begin{remark}
This is reminiscent of the formula for a geometric series of numbers: $\sumn[0] a^n=\frac{1}{1-a}$.
\end{remark}
\begin{proof}
Since $\norm{a}<1$, we have $\sumn[0]\norm{a^n}\le\sumn[0]\norm{a}^n<\infty$, so $\sumn[0] a^n$ is absolutely convergent (and thus convergent, since $\As$ is a Banach space). Define $b=\sumn[0] a^n$. For each $N\in\N_0$, we have $(\1-a)(\sum_{n=0}^N a^n)=\sum_{n=0}^N a^n-\sum_{n=0}^N a^{n+1}=1-a^{N+1}$. Taking the limit as $N\to\infty$ yields $(\1-a)b=\1$, since $\norm{a^{N+1}}\le\norm{a}^{N+1}\to0$ as $N\to\infty$. A similar argument shows that $b(\1-a)=\1$. Thus $b=(\1-a)^{-1}$.
\end{proof}

\begin{proposition}\label{invopen}
$\As\inv$ is an open set in $\As$.
\end{proposition}
\begin{proof}
Suppose $a\in\As\inv$ and $b\in\As$ such that $\norm{a-b}<\frac{1}{\norm{a^{-1}}}$. Then we have:
\begin{equation*}
    \norm{1-a^{-1}b}=\norm{a^{-1}a-a^{-1}b}=\norm{a^{-1}(a-b)}\le\norm{a^{-1}}\norm{a-b}<1
\end{equation*}
By the previous theorem, $a^{-1}b\in\As\inv$, and so $b=a(a^{-1}b)\in\As\inv$. Thus $\As\inv$ is open in $\As$.
\end{proof}

\begin{proposition}\label{invhomeo}
The inversion map $f:\As\inv\to\As\inv$, $f(a)=a^{-1}$ is a homeomorphism.
\end{proposition}
\begin{proof}
Suppose $a\in\As\inv$ and $b\in\As$ such that $\norm{a-b}\le\frac{1}{2\norm{a^{-1}}}$. Then $\norm{a-b}<\frac{1}{\norm{a^{-1}}}$, so by the previous proposition, we have $b\in\As\inv$. We also have:
\begin{align*}
    \norm{b^{-1}}-\norm{a^{-1}}&\le\abs{\norm{b^{-1}}-\norm{a^{-1}}}\le\norm{b^{-1}-a^{-1}}=\norm{b^{-1}aa^{-1}-b^{-1}ba^{-1}}=\norm{b^{-1}(a-b)a^{-1}}\\
    &\le\norm{b^{-1}}\norm{a-b}\norm{a^{-1}}\le\norm{b^{-1}}\frac{1}{2\norm{a^{-1}}}\norm{a^{-1}}=\frac{\norm{b^{-1}}}{2}
\end{align*}
Rearranging yields $\frac{\norm{b^{-1}}}{2}\le\norm{a^{-1}}$, and so $\norm{b^{-1}}\le2\norm{a^{-1}}$. We now have:
\begin{equation*}
    \norm{b^{-1}-a^{-1}}\le\norm{b^{-1}}\norm{a-b}\norm{a^{-1}}\le2\norm{a^{-1}}\norm{a-b}\norm{a^{-1}}=2\norm{a^{-1}}^2\norm{a-b}
\end{equation*}

Now suppose $a\in\As\inv$ and $\ep>0$. Set $\delta=\min\left\{\frac{1}{2\norm{a^{-1}}},\frac{\ep}{2\norm{a^{-1}}^2}\right\}$. Then for all $b\in\As$, $\norm{a-b}<\delta$, we have $\norm{a-b}<\frac{1}{2\norm{a^{-1}}}$, so $b\in\As\inv$ and $\norm{b^{-1}-a^{-1}}\le2\norm{a^{-1}}^2\norm{a-b}<2\norm{a^{-1}}^2\frac{\ep}{2\norm{a^{-1}}^2}=\ep$. Thus $f$ is continuous. Since $f$ maps each element of $\As\inv$ to its inverse, we have $f^{-1}=f$, so $f^{-1}$ is also continuous. Thus $f$ is a homeomorphism.
\end{proof}

\Cref{multcontinuous} and \Cref{invhomeo} show that $\As\inv$ is a topological group (with the subspace topology from the norm topology on $\As$). As a subset of a normed vector space, it is also a metric space (and thus Hausdorff).

\begin{lemma}
Suppose $\As$ is a normed algebra and $\Bs$ is a subalgebra (resp. ideal) of $\As$. Then $\overline{\Bs}$ (the closure of $\Bs$) is also a subalgebra (resp. ideal) of $\As$.
\end{lemma}
\begin{proof}
Since $\Bs$ is a subspace of $\As$, so is $\overline{\Bs}$. Suppose $a,b\in\overline{\Bs}$. Then there are sequences $(a_n),(b_n)$ in $\Bs$ such that $a_n\to a$ and $b_n\to b$. Since $\Bs$ is a subalgebra of $\As$, we have $a_nb_n\in\Bs$ for all $n\in\N$. Since $a_n\to a$ and $b_n\to b$, we have $a_nb_n\to ab$. Thus $ab\in\overline{\Bs}$, and so $\overline{\Bs}$ is a subalgebra of $\As$.\\

Now suppose $\Bs$ is an ideal of $\As$, $a\in\As$ and $b\in\overline{\Bs}$. Then there is a sequence $(b_n)$ in $\Bs$ such that $b_n\to b$. Since $\Bs$ is an ideal of $\As$, we have $ab_n,b_na\in\Bs$ for all $n\in\N$. Since $b_n\to b$, we have $ab_n\to ab$ and $b_na\to ba$. Thus $ab,ba\in\overline{\Bs}$, and so $\overline{\Bs}$ is an ideal of $\As$.
\end{proof}

\begin{proposition}
Suppose $\As$ is a unital Banach algebra and $\Is$ is a proper ideal of $\As$. Then $\overline{\Is}$ (the closure of $\Is$) is also a proper ideal of $\As$.
\end{proposition}
\begin{proof}
Since $\Is$ is a proper ideal of $\As$, it cannot contain any invertible elements, i.e. $\Is\cap\As\inv=\emp$, so $\Is\sub\As\exc\As\inv$. Since $\As\inv$ is open, its complement $\As\exc\As\inv$ is closed, and so $\overline{\Is}\sub\As\exc\As\inv$. In particular, $\overline{\Is}\ne\As$. Since $\Is$ is an ideal of $\As$, so is $\overline{\Is}$, by the previous lemma. Thus $\overline{\Is}$ is a proper ideal of $\As$.
\end{proof}

\begin{theorem}\label{maxclosed}
Every maximal ideal of a unital Banach algebra is closed.
\end{theorem}
\begin{proof}
Suppose $\As$ is a unital Banach algebra and $\Ms$ is a maximal ideal of $\As$. By the previous proposition, $\overline{\Ms}$ is a proper ideal of $\As$. Since $\Ms$ is maximal, we must have $\overline{\Ms}=\Ms$, and so $\Ms$ is closed.
\end{proof}

The following example shows that a maximal ideal of a \textit{non-unital} Banach algebra may not be closed:

\begin{example}
Suppose $X$ is an infinite-dimensional Banach space and define vector multiplication on $X$ by $ab=0$ for all $a,b\in X$. Then $X$ is a non-unital abelian Banach algebra. Now suppose $f$ is an unbounded linear functional\footnotemark\ on $X$. Then $\ker(f)$ is a dense subspace of $X$ with codimension 1. Since $ab=ba=0\in\ker(f)$, it is also an ideal of $X$, and since it has codimension 1, it must be maximal. However, it is dense in $X$ and not equal to $X$, so it cannot be closed.
\end{example}
\footnotetext{This is always possible, as every infinite-dimensional Banach space has unbounded linear functionals. In general, however, these functionals cannot be constructed explicitly, and proving they exist requires the axiom of choice. See \cite[Example 4.2, Page 126]{heil} for more details.}

We will now extend \Cref{quotientalg} to Banach algebras:

\begin{theorem}\label{quotientbalg}
Suppose $\As$ is a Banach algebra and $\Is$ is a closed ideal of $\As$. Define the norm $\norm{\cdot}$ on $\As/\Is$ by $\norm{a+\Is}=\inf_{b\in\Is} \norm{a+b}$. Then $\As/\Is$ is a Banach algebra, and the quotient map $q:\As\to\As/\Is$, $q(a)=a+\Is$ is a homorphism and $\norm{q}\le1$. If $\As$ is unital, then $\As/\Is$ is unital with identity $\1_{\As/\Is}=\1_\As+\Is$.
\end{theorem}
\begin{proof}
Since $\Is$ is an ideal of $\As$, by \Cref{quotientalg}, the quotient $\As/\Is$ is an algebra. Since $\Is$ is a closed ideal (and thus a closed subspace) of $\As$, $\As/\Is$ is also a Banach space. We now show that the quotient norm is an algebra norm:
\begin{align*}
    \norm{[a][b]}_{\As/\Is}
    &=\norm{[ab]}_{\As/\Is}=\inf_{c\in\Is} \norm{ab+c}=\inf_{c_1,c_2\in\Is} \norm{ab+ac_2+bc_1+c_1c_2}=\inf_{c_1,c_2\in\Is} \norm{(a+c_1)(b+c_2)} \\
    &\le\inf_{c_1,c_2\in\Is} \norm{a+c_1}\norm{b+c_2}\le\inf_{c_1\in\Is} \norm{a+c_1}\inf_{c_2\in\Is} \norm{b+c_2}=\norm{[a]}_{\As/\Is}\norm{[b]}_{\As/\Is}
\end{align*}
Thus the quotient norm is an algebra norm, and so $\As/\Is$ is a Banach algebra.\\

We also have $\norm{q(a)}=\norm{a+\Is}=\inf_{b\in\Is} \norm{a+b}\le\norm{a+0}=\norm{a}$. Thus $q$ is bounded and $\norm{q}\le1$.\\

Finally, suppose $\As$ is unital. Then for all $a\in\As$, we have:
\begin{align*}
    (a+\Is)(\1_\As+\Is)=a\1_\As+\Is=a+\Is && (\1_\As+\Is)(a+\Is)=\1_\As a+\Is=a+\Is
\end{align*}
Thus $\As/\Is$ is unital with identity $1_{\As/\Is}=1_\As+\Is$.
\end{proof}
\begin{remark}
One can also show that if $\Is\ne\As$, then $\norm{q}=1$. This is slightly more involved, and requires Riesz's lemma (see \cite{rieszlemma}).
\end{remark}


We will see in \Cref{characters} what happens in the case that $\Is$ is a \emph{maximal} ideal of $\As$.\\

\ind We conclude this section with the definition and some basic results about isometries, which we will revisit later in \Cref{cstas}.

\begin{definition}
Suppose $X$ and $Y$ are normed vector spaces. An \textbf{isometry} (or \textit{linear isometry}) from $X$ to $Y$ is a linear map $T:X\to Y$ such that $\norm{Tx}=\norm{x}$ for all $x\in X$.\\
A \textbf{isometric isomorphism} (or \textit{global linear isometry}) is an isometry that is also surjective.
\end{definition}
Basically, an isometry is a homomorphism of normed vector spaces, i.e. a map that preserves the vector space structure as well as the norm structure. Likewise, an isometric isomorphism is an isomorphism of normed vector spaces.

\xnewpage

\begin{lemma}\label{isomcontinj}
Every isometry is continuous and injective.
\end{lemma}
\begin{proof}
Since $T:X\to Y$ is an isometry, we have $\norm{Tx}=\norm{x}$ for all $x\in X$, so $T$ is bounded (and thus continuous) and $\norm{T}\le1$. Now suppose $T(x)=T(y)$ for some $x,y\in X$. Then we have $\norm{x-y}=\norm{T(x)-T(y)}=\norm{0}=0$, so $x=y$. Thus $T$ is injective.
\end{proof}

\begin{proposition}\label{isomimclosed}
Suppose $X$ is a Banach space, $Y$ is a normed vector space and $T:X\to Y$ is an isometry. Then $\im(T)$ is closed in $Y$.
\end{proposition}
\begin{proof}
Suppose $(y_n)$ is a sequence in $\im(T)$ that converges to some $y\in Y$. We want to show that $y\in\im(T)$. For each $n\in\N$, since $y_n\in\im(T)$, there exists $x_n\in X$ such that $T(x_n)=y_n$. Since the sequence $(y_n)$ converges in $Y$, it is Cauchy, and since $T$ is an isometry, we have $\norm{y_m-y_n}=\norm{x_m-x_n}$ for all $m,n\in\N$. Thus the sequence $(x_n)$ is also Cauchy. Since $X$ is complete, $(x_n)$ converges to some $x\in X$. By the previous lemma, $T$ is continuous, so $(y_n)=(Tx_n)$ converges to $Tx$. Thus $y=Tx\in\im(T)$, and so $\im(T)$ is closed in $Y$.
\end{proof}
\begin{remark}
A slight modification of this proof shows that $\im(T)$ is also complete, i.e. it is a Banach space (even if $Y$ is not).
\end{remark}

\ind Since every isometry $T:X\to Y$ is injective, it has an inverse $T^{-1}:\im(T)\to X$. Also, since $\norm{Tx}=\norm{x}$ for all $x\in X$, it follows that $T^{-1}$ is also an isometry (and thus continuous). In particular, every isometry is a homeomorphism onto its image. Of course, the converse is not true (you can multiply an isometry by any nonzero constant and it would still be a homeomorphism, but not an isometry).

\subsection{The Spectrum}
\ind One of the most important concepts in linear algebra is that of eigenvalues and eigenvectors. Recall that if $V$ is a vector space and $T:V\to V$ is a linear operator, we say that $v\in V\exc\{0\}$ is an \emph{eigenvector} of $T$ with \emph{eigenvalue} $\lam\in\F$ if $Tv=\lam v$. If $\F=\C$ and $V$ is finite-dimensional, every linear operator $T:V\to V$ must have at least one eigenvector (this is a consequence of the fundamental theorem of algebra). However, for infinite-dimensional spaces, this fails horribly, even in the simplest cases:

\begin{example}
Suppose $T:C[0,1]\to C[0,1]$, $(Tf)(x)=xf(x)$. Then $T$ has NO eigenvalues or eigenvectors. To see this, suppose $f$ is an eigenvector of $T$ with eigenvalue $\lam$, i.e. $Tf=\lam f$. Then we have $xf(x)=\lam f(x)$, i.e. $(x-\lam)f(x)=0$. This must be true for all $x\in[0,1]$, which implies that $f\equiv 0$, a contradiction. Thus $T$ has no eigenvalues or eigenvectors.
\end{example}

\ind To develop a meaningful analogous concept for infinite-dimensional spaces, we need to introduce a different definition: $\lam\in\F$ is in the \emph{spectrum} of $T$ if $T-\lam I$ is \emph{not} invertible (where $I$ is the identity operator). For finite-dimensional spaces, this is equivalent to $\lam$ being an eigenvalue of $T$, but for infinite-dimensional spaces, this definition is much more general (in particular, if $\F=\C$, this eliminates the problem of such values not existing). It also has the advantage of being applicable to all unital algebras, not just those that arise from linear operators on vector spaces.\\

\ind The theory of eigenvalues and eigenvectors, now known as \emph{spectral theory}, has a very rich history, which you can read about in \cite{steen}.

\begin{definition}
Suppose $\As$ is a unital algebra and $a\in\As$. The \textbf{spectrum} of $a$, denoted by $\sigma(a)$, is given by:
\begin{equation*}
    \sigma(a)=\left\{\lam\in\F\mid a-\lam\1\notin\As\inv\right\}
\end{equation*}
The \textbf{resolvent set} of $a$, denoted by $\rho(a)$, is given by $\rho(a)=\F\exc\sigma(a)$.\\
In other words, the spectrum of $a$ is the set of all numbers $\lam\in\F$ such that $a-\lam\1$ is \textit{not} invertible.
\end{definition}

\begin{center}
\begin{tikzpicture}
\begin{scope}[xshift=0cm] 
\colorlet{scol}{red}
\draw (-1.7,0)--(1.7,0) node[right]{$\Re(z)$};
\draw (0,-1.7)--(0,1.7) node[above]{$\Im(z)$};
\foreach \a in {(1.2,1),(-0.5,-0.9),(-1.2,0.4),(0.9,-0.5),(-0.4,1.1),(0.5,0.2)}
\fill[scol!50] \a circle (2pt);
\end{scope}
\begin{scope}[xshift=5cm] 
\colorlet{scol}{Green}
\draw (-1.7,0)--(1.7,0) node[right]{$\Re(z)$};
\draw (0,-1.7)--(0,1.7) node[above]{$\Im(z)$};
\foreach \c in {(-0.7,0.8),(0.7,0.8)}
\draw[thick,scol,draw opacity=0.5,fill=scol,fill opacity=0.5] \c circle (0.2);
\draw[thick,scol,draw opacity=0.5,fill=scol,fill opacity=0.5] (-1.1,-0.2) arc(-180:0:1.1) --(0.9,-0.2) arc(0:-180:0.9) --cycle;
\end{scope}
\begin{scope}[xshift=10cm] 
\colorlet{scol}{blue}
\draw (-1.7,0)--(1.7,0) node[right]{$\Re(z)$};
\draw (0,-1.7)--(0,1.7) node[above]{$\Im(z)$};
\draw[thick,scol,draw opacity=0.5,fill=scol,fill opacity=0.5] (-0.2,0.3)--(1,0.3)--(1,0.8)--(-0.2,0.8)--cycle;
\draw[thick,scol,draw opacity=0.5,fill=scol,fill opacity=0.5] (0.9,-0.7) circle (0.3);
\draw[very thick,scol!50] (-0.6,1.2)--(-1.2,0.2)--(-0.4,-0.8);
\draw[very thick,scol!50] (0.7,1.4) arc(120:0:0.5);
\fill[scol!50] (-0.9,-1.1) circle (2pt);
\end{scope}
\end{tikzpicture}\\
Some possible spectra of elements.
\end{center}

\begin{examples}\leavevmode
\begin{enumerate}
    \item Suppose $V$ is a finite-dimensional vector space and $T:V\to V$ is a linear operator. Then $\lam\in\sigma(T)$ if and only if $\lam$ is an eigenvalue of $T$.
    \item Suppose $X$ is a locally compact Hausdorff space and $f\in C_0(X)$. Then $\sigma(f)=\im(f)$ (the image of $f$).
    \item Suppose $(X,\As,\mu)$ is a measure space and $f\in L^\infty(X,\As,\mu)$. Then:
    \begin{equation*}
        \sigma(f)=\left\{\lam\in\F\mid\mu(\{x\in X\mid\abs{f(x)-\lam}<\ep\})>0\text{ for all }\ep>0\right\}
    \end{equation*}
    This is sometimes known as the \emph{essential image} of $f$.
\end{enumerate}
\end{examples}

\ind The definition of the spectrum requires $\As$ to be unital, as it directly uses the identity and invertibility. We will restrict our discussion of the spectrum to unital Banach algebras (and usually, unital \emph{complex} Banach algebras) as far as possible. When required, we define the spectrum of $a\in\As$, where $\As$ is a \emph{non-unital} Banach algebra, to be the spectrum of $a$ as an element of the unitization $\As_1$ (see \Cref{embedding}).
 
\begin{definition}
Suppose $U\sub\C$ is open and $\As$ is a complex Banach algebra. A function $f:U\to\As$ is \textbf{holomorphic} on $U$ if for all $z_0\in U$, the limit $\lim_{z\to z_0} \frac{f(z)-f(z_0)}{z-z_0}$ exists in $\As$.
\end{definition}
\begin{remark}
This notion is sometimes known as \textit{strongly holomorphic}, as opposed to \textit{weakly holomorphic}, which means that for all $\phi\in\As^*$, the function $\phi\circ f:U\to\C$ is holomorphic on $U$. In fact, these are equivalent (see \cite[Theorem 8.20]{ward} for a proof), so there is no need to distinguish them, and so we will call them both ``holomorphic''.
\end{remark}

\begin{theorem}\label{reshol}
Suppose $\As$ is a unital complex Banach algebra and $a\in\As$. Then the function $f:\rho(a)\to\As$, $f(z)=(a-z\1)^{-1}$ is holomorphic on $\rho(a)$.
\end{theorem}
\begin{remark}
The function $f$ is sometimes known as the \textit{resolvent} of $a$. By the definition of $\rho(a)$, it is well-defined. Also, as we will show in \Cref{ftba}, $\rho(a)$ is open, so it makes sense to say that $f$ is holomorphic on $\rho(a)$.
\end{remark}
\begin{proof}
Suppose $z,z_0\in\rho(a)$. Since $a-z\1$ and $a-z_0\1$ are linear functions of $a$, they commute, and so their inverses also commute ($*$). This yields:
\begin{align*}
    f(z)-f(z_0)&=(a-z\1)^{-1}-(a-z_0\1)^{-1} \\
    &=\textcolor{red!70!black}{(a-z_0\1)(a-z_0\1)^{-1}}(a-z\1)^{-1}-\textcolor{red!70!black}{(a-z\1)(a-z\1)^{-1}}(a-z_0\1)^{-1} \\
    &=\textcolor{Green!80!black}{(a-z_0\1)}(a-z\1)^{-1}(a-z_0\1)^{-1}-\textcolor{Green!80!black}{(a-z\1)}(a-z\1)^{-1}(a-z_0\1)^{-1} \tag{$*$} \\
    &=(\textcolor{Green!80!black}{(a-z_0\1)-(a-z\1)})(a-z\1)^{-1}(a-z_0\1)^{-1} \\
    &=(\textcolor{Green!80!black}{(z-z_0)\1})(a-z\1)^{-1}(a-z_0\1)^{-1} \\
    &=(z-z_0)(a-z\1)^{-1}(a-z_0\1)^{-1}
\end{align*}
Dividing by $z-z_0$, we get $\frac{f(z)-f(z_0)}{z-z_0}=(a-z\1)^{-1}(a-z_0\1)^{-1}$. Taking the limit as $z\to z_0$ yields:
\begin{equation*}
    \lim_{z\to z_0} \frac{f(z)-f(z_0)}{z-z_0}=\lim_{z\to z_0} (a-z\1)^{-1}(a-z_0\1)^{-1}=(a-z_0\1)^{-2}\in\As
\end{equation*}
Thus $f$ is holomorphic on $\rho(a)$.
\end{proof}

\ind The next lemma is a generalization of Liouville's theorem, which states that every bounded entire function $f:\C\to\C$ is constant\footnote{An \emph{entire} function is a function that is holomorphic on all of $\C$.} (see \cite[Theorem 3.3.1, Page 89]{taylor} for a proof). This is the starting point of the discussion we will soon embark on that is specific to complex normed algebras. Indeed, almost every result in functional analysis that holds for complex normed algebras but not for real ones can be traced back to Liouville's theorem.

\begin{lemma}[Louville's Theorem for Normed Vector Spaces]\label{liouvillenvs}
Suppose $X$ is a complex normed vector space. Then every bounded entire function $f:\C\to X$ is constant.
\end{lemma}
\begin{proof}
Suppose $f$ is not constant. Then there exist $z,w\in\C$ such that $f(z)\ne f(w)$, i.e. $f(z)-f(w)\ne0$. By the Hahn-Banach theorem, there exists $\phi\in X^*$ such that $\phi(f(z)-f(w))\ne0$, i.e. $\phi(f(z))\ne\phi(f(w))$. Thus $\phi\circ f$ is not constant. Since $f$ is bounded and entire, the composition $\phi\circ f:\C\to\C$ is also bounded and entire, so by Liouville's theorem, it must be constant, a contradiction. Thus $f$ is constant.
\end{proof}

\begin{theorem}[Fundamental Theorem of Banach Algebras]\label{ftba}
Suppose $\As$ is a unital complex Banach algebra and $a\in\As$. Then $\sigma(a)$ is non-empty and compact.
\end{theorem}
\begin{proof}
Suppose $\lam\in\C$ and $\abs{\lam}>\norm{a}$. We want to show that $a-\lam\1$ is invertible. Since $\abs{\lam}>\norm{a}$, we have $\norm{\frac{1}{\lam}a}<1$, so by \Cref{neumann}, $\1-\frac{1}{\lam}a$ is invertible, and so is $a-\lam\1=-\lam\left(\1-\frac{1}{\lam}a\right)$. Thus every $\lam\in\C$, $\abs{\lam}>\norm{a}$ is in $\rho(a)$, and so $\sigma(a)\sub\Dbar(0,\norm{a})$. Thus $\sigma(a)$ is bounded.\\

Define $g:\C\to\As$, $g(z)=a-z\1$. Then for all $z,w\in\C$, we have $\norm{g(z)-g(w)}=\norm{(a-z\1)-(a-w\1)}=\norm{-(z-w)\1}=\abs{z-w}$. Thus $g$ is continuous. By \Cref{invopen}, the set $\As\inv$ is open in $\As$, so its pre-image $g^{-1}(\As\inv)=\rho(a)$ is open in $\C$. Thus $\sigma(a)=\C\exc\rho(a)$ is closed. Since $\sigma(a)$ is closed and bounded, by the Heine-Borel theorem, it is compact.\\

Suppose $\sigma(a)=\emp$. Then $\rho(a)=\C$. Thus the function $f$ from \Cref{reshol} is holomorphic on $\C$ (and thus entire). We also have $\norm{f(z)}=\norm{(a-z\1)^{-1}}\ge\frac{1}{\norm{a-z\1}}\ge\frac{1}{\norm{a}+\abs{z}}\to0$ as $z\to\infty$. Thus $f$ is bounded on $\C$. By \Cref{liouvillenvs}, $f$ is constant, and since it vanishes as $z\to\infty$, it must be identically zero. This is a contradiction, as $f(z)=(a-z\1)^{-1}\in\As\inv$. Thus $\sigma(a)\ne\emp$.
\end{proof}
\begin{remark}
If $\As$ is a unital \textit{real} Banach algebra, $\sigma(a)$ is still compact, but it may be empty, as the next example shows.
\end{remark}

\begin{example}
Suppose $A\in M_2(\R)$, $A=\mmat{0}{-1}{1}{0}$. Then $\sigma(A)=\emp$ as its characteristic polynomial is $\lam^2+1$, which has no real roots. If we instead view $A$ as an element of $M_2(\C)$, then $\sigma(A)=\{i,-i\}$.
\end{example}

\ind In general, there are no other restrictions on the spectrum of an element, besides being non-empty and compact. This can be seen as follows: Suppose $K\sub\C$ is non-empty and compact. Then we can find a sequence $(a_n)_{n=1}^\infty$ in $K$ that is dense in $K$. Define $T:\ell^\infty\to\ell^\infty$ by $T(x_1,x_2,x_3,...)=(a_1x_1,a_2x_2,a_3x_3,...)$. Then $T\in\Bs(\ell^\infty)$ and $\sigma(T)=K$.

\begin{definition}
Suppose $\As$ is a unital Banach algebra and $a\in\As$. The \textbf{spectral radius} of $\As$, denoted by $r(a)$, is given by:
\begin{equation*}
    r(a)=\sup_{\lam\in\sigma(a)} \abs{\lam}
\end{equation*}
In other words, the spectral radius of $a$ is the largest absolute value of all the numbers in $\sigma(a)$.
\end{definition}
\begin{remarks}\leavevmode
\begin{enumerate}
    \item Since $\sigma(a)$ is non-empty and compact, this supremum is always attained and finite, i.e. it is a maximum.
    \item By the proof of the previous theorem, we have $\sigma(a)\sub\Dbar(0,\norm{a})$. Thus $r(a)\le\norm{a}$.
\end{enumerate}
\end{remarks}

\begin{center}
\begin{tikzpicture}
\colorlet{scol}{Mulberry}
\draw (-1.7,0)--(1.7,0);
\draw (0,-1.7)--(0,1.7);
\draw[dashed] (0,0) circle (1.5);
\draw[thick,scol,draw opacity=0.5,fill=scol,fill opacity=0.5,bend angle=30] (1.2,0.4)--(3/4,{3/4*sqrt(3)})--(-0.6,0.6)--(-0.2,-0.4)--(-0.4,-1.2) to[bend right] cycle;
\draw[very thick,scol!50] (-1.1,0.5)--(-0.4,1.2);
\foreach \a in {(0.2,1.2),(1,-0.8),(-0.8,0.3),(-0.6,-1)}
\fill[scol!50] \a circle (2pt);
\draw[<->,thick] (0,0)--({-3/4*sqrt(2)},{-3/4*sqrt(2)}) node[pos=0.4,left,xshift=-1pt,yshift=1pt]{$r(a)$};
\end{tikzpicture}\\
The spectral radius of an element.
\end{center}

\begin{example}
Suppose $A=\mmat{0}{1}{0}{0}$. Then $A-\lam I$ is invertible for all $\lam\in\C$ except $\lam=0$. Thus $\sigma(A)=\{0\}$, and so $r(A)=0$. However, since $A\ne0$, we have $\norm{A}>0$ (in fact, $\norm{A}=1$).
\end{example}

\ind The above example shows that $r(a)$ is not necessarily equal to $\norm{a}$. We now present two classes of examples where equality always holds. These are special cases of a more general phenomenon we will encounter in \Cref{cstas}, when we discuss abelian \cstas.

\begin{examples}\leavevmode
\begin{enumerate}
    \item Suppose $X$ is a locally compact Hausdorff space and $f\in C_0(X)$. Then $r(f)=\sup_{\lam\in\sigma(f)} \abs{\lam}=\sup_{x\in X} \abs{f(x)}$. This is the supremum of $\abs{f}$, which by definition is $\norm{f}$.
    \item Suppose $(X,\As,\mu)$ is a measure space and $f\in L^\infty(X,\As,\mu)$. Then:
    \begin{equation*}
        r(f)=\sup_{\lam\in\sigma(f)} \abs{\lam}=\inf\left\{M\ge0\mid\mu(\{x\in X\mid\abs{f(x)}>M\})=0\right\}
    \end{equation*}
    This is the essential supremum of $\abs{f}$, which by definition is $\norm{f}$.
\end{enumerate}
\end{examples}

\xnewpage

\ind The next theorem, which first appeared in \cite{gelfand}, is a remarkable result that links the algebraic and analytic aspects of Banach algebras. To prove it, we will need a special case of Laurent's theorem, which says that a series of the form $\sumn[0] c_nz^{-n}$ converges for all $\abs{z}>R$, where $R=\limsupn\abs{c_n}^{1/n}$.

\begin{theorem}[Spectral Radius Formula or Beurling-Gelfand Formula]\label{specradius}
Suppose $\As$ is a unital complex Banach algebra and $a\in\As$. Then:
\begin{equation*}
    r(a)=\limn\norm{a^n}^{1/n}=\inf_{n\in\N}\norm{a^n}^{1/n}
\end{equation*}
\end{theorem}
\begin{remarks}\leavevmode
\begin{enumerate}
    \item The fact that the above limit exists is highly nontrivial (and thus surprising in itself), and is part of the theorem.
    \item The spectrum (and thus the spectral radius) does not depend on the norm used. As such, the above formula holds for \textit{every} norm on $\As$ that makes it a Banach algebra.
\end{enumerate}
\end{remarks}
\begin{proof}
We first show that $r(a)\le\inf_{n\in\N}\norm{a^n}^{1/n}$. Suppose $\lam\in\C$ and $n\in\N$. Then we have:
\begin{equation*}
    a^n-\lam^n\1=(a-\lam\1)(a^{n-1}+\lam a^{n-2}+\cdots+\lam^{n-2}a+\lam^{n-1}\1)
\end{equation*}
Also, these factors commute as they are polynomials in $a$. Thus if $a-\lam\1$ is not invertible, neither is $a^n-\lam^n\1$. In other words, if $\lam\in\sigma(a)$, then $\lam^n\in\sigma(a^n)$. This yields $\abs{\lam}^n=\abs{\lam^n}\le\norm{a^n}$, and so $\abs{\lam}\le\norm{a^n}^{1/n}$. Taking the supremum over all $\lam\in\sigma(a)$, we get $r(a)\le\norm{a^n}^{1/n}$, and taking the infimum over all $n\in\N$, we get $r(a)\le\inf_{n\in\N}\norm{a^n}^{1/n}$.\\

We now show that $r(a)\ge\limsupn\norm{a^n}^{1/n}$. Define $f:\rho(a)\to\As$ by $f(z)=(a-z\1)^{-1}$. By definition, this is well-defined on $\rho(A)$, and for $\abs{z}>\norm{a}$, we can use \Cref{neumann} to get:
\begin{align*}
    (a-z\1)^{-1}=-\frac{1}{z}\left(\1-\frac{1}{z}a\right)^{-1}=-\frac{1}{z}\sumn[0] \left(\frac{1}{z}a\right)^n=-\sumn[0] z^{-n-1}a^n
\end{align*}
This is a Laurent series in $z$ with no terms in positive powers of $z$. By Laurent's theorem, it converges absolutely for all $z\in\C$, $\abs{z}>R$, where $R=\limsupn\norm{a^n}^{1/n}$. Thus every $z\in\C$, $\abs{z}>R$ is in $\rho(a)$, and so $r(a)\le R$.\\

We now have $\inf_{n\in\N}\norm{a^n}^{1/n}\ge r(a)\ge\limsupn\norm{a^n}^{1/n}$. Since $\limsupn\beta_n\ge\liminfn\beta_n\ge\inf_{n\in\N}\beta_n$ for any sequence $(\beta_n)$ in $\R$, it follows that:
\begin{equation*}
    \inf_{n\in\N}\norm{a^n}^{1/n}\ge r(a)\ge\limsupn\norm{a^n}^{1/n}\ge\liminfn\norm{a^n}^{1/n}\ge\inf_{n\in\N}\norm{a^n}^{1/n}
\end{equation*}
Thus $\limn\norm{a^n}^{1/n}$ exists and $\limn\norm{a^n}^{1/n}=\inf_{n\in\N}\norm{a^n}^{1/n}=r(a)$.
\end{proof}

\begin{example}
Suppose $A=\mmat{1}{1}{0}{2}$. Then the eigenvalues of $A$ are $1$ and $2$, so $\rho(A)=2$. We also have $A^n=\mmat{1}{2^n-1}{0}{2^n}$, and (after \emph{a lot} of computation), $\norm{A^n}=\sqrt{2^{2n}-2^n+1+(2^n-1)\sqrt{2^{2n}+1}}$. As $n\to\infty$, this is asymptotic to $2^{n+\tfrac{1}{2}}$, so $\norm{A^n}^{1/n}\sim 2^{1+\tfrac{1}{2n}}$. Thus $\limn\norm{A^n}^{1/n}=\limn 2^{1+\tfrac{1}{2n}}=2^1=2$, which agrees with the spectral radius formula (for comparison, $\norm{A}=\sqrt{3+\sqrt{5}}\approx2.28825$, while $\norm{A^{100}}^{1/100}\approx2.00694$).
\end{example}

\ind While the spectral radius is perfectly well-defined for elements of \emph{real} unital Banach algebras, it is customary to use the \emph{complex} spectral radius even if $\As$ is real. One reason is intuition: It makes more sense to say that the spectral radius of $\mmat{0}{-1}{1}{0}$ is $1$, even though as an element of $M_2(\R)$, the definition gives it a spectral radius of $0$ (since it has no real eigenvalues). Another reason is to ensure various important theorems still hold, such as \Cref{specradius}. In particular, this convention will find its way into the proofs of the \hyperref[gelfandrep]{Gelfand representation theorem} and the \hyperref[gelfandnaimark]{Gelfand-Naimark theorem}, which will allow us to conclude that they still hold if $\As$ is real.

\begin{corollary}\label{specradiuscor}
Suppose $\As$ is a unital complex Banach algebra and $a\in\As$. Then $\limn\norm{a^n}=0$ if and only if $r(a)<1$.
\end{corollary}
\begin{proof}
Note that $r(a)=\limn\norm{a^n}^{1/n}$ implies $\limn (r(a))^n=\limn\norm{a^n}$. Thus:
\begin{equation*}
    \limn\norm{a^n}=0 \iff \limn (r(a))^n=0 \iff r(a)<1 \qedhere
\end{equation*}
\end{proof}

\ind The last corollary is especially important in stability theory, where dynamical systems are modeled by matrices and one would like to know if such a system stabilizes over time, blows up to infinity or does something else.\\

\ind While the proofs of \Cref{ftba} and \Cref{specradius} here relied heavily on methods from complex analysis, it is actually possible to prove them without any complex analysis (of course, still assuming the algebra is complex). See \cite[Theorem 1.2.8, Pages 10-13]{kaniuth}.

\begin{theorem}[Spectral Mapping Theorem]
Suppose $\As$ is a unital complex Banach algebra, $a\in\As$ and $p(z)=c_nz^n+c_{n-1}z^{n-1}+\cdots+c_1z+c_0$ is a polynomial. Then $\sigma(p(a))=p(\sigma(a))$.\\
Here, $p(a)\in\As$ denotes the element $c_na^n+c_{n-1}a^{n-1}+\cdots+c_1a+c_0\1$, while $p(\sigma(a))\subset\C$ denotes the set $\{p(\zeta)\mid\zeta\in\sigma(a)\}$, i.e. the image of $\sigma(a)$ under the polynomial $p$.
\end{theorem}
Basically, when an element of $\As$ is mapped by a polynomial, its spectrum is also mapped by the same polynomial. In other words, the following diagram commutes:
\begin{center}
\begin{tikzcd}[column sep=4em]
a \arrow[r, "p"] \arrow[d, "\sigma"'] & p(a) \arrow[d, "\sigma"] \\
\sigma(a) \arrow[r, "p"'] & \sigma(p(a))
\end{tikzcd}
\end{center}
\begin{proof}
Suppose $\mu\in\C$. Then the equation $p(z)=\mu$ is a polynomial equation of degree $n$, so by the fundamental theorem of algebra, it has $n$ solutions $\lam_1,\lam_2,...,\lam_n\in\C$ (including multiplicity). This yields $p(z)-\mu=c_n(z-\lam_1)(z-\lam_2)\cdots(z-\lam_n)$, and similarly, $p(a)-\mu\1=c_n(a-\lam_1\1)(a-\lam_2\1)\cdots(a-\lam_n\1)$. Note that all the factors $a-\lam_k\1$ commute as they are polynomials in $a$. By definition, $\mu\notin\sigma(p(a))$ if and only if $p(a)-\mu\1\notin\As\inv$. Since all the factors $a-\lam_k\1$ commute, this is equivalent to all of them being invertible. In other words, $\lam_1,\lam_2,...,\lam_n\notin\sigma(a)$. Since $\lam_1,\lam_2,...,\lam_n$ are the zeros of $p(z)-\mu$, this is equivalent to $\mu\notin p(\sigma(a))$. Thus $\sigma(p(a))=p(\sigma(a))$.
\end{proof}
\begin{remark}
Using functional calculus, this theorem can be generalized to larger classes of functions than just polynomials, see \cite[Theorem 10.28, Page 263]{grudin}, \cite[Chapter VIII, Theorem 2.7, Page 239]{conway2} and \cite[Theorem 2.1.14, Page 43]{murphy}.
\end{remark}

\begin{example}
Suppose $A=\mmat{3}{2}{1}{4}$ and $p(x)=x^3+10x^2+8x+5$. Then $\sigma(A)=\{2,5\}$, so $p(\sigma(A))=\{p(2),p(5)\}=\{69,420\}$. On the other hand, $p(A)=\mmat{186}{234}{117}{303}$, so $\sigma(p(A))=\{69,420\}$ as expected.
\end{example}

\ind In general, Banach algebras have a lot of non-invertible elements. This is in contrast to $\R$ and $\C$, where the only non-invertible element is $0$. Algebras where every nonzero element has a multiplicative inverse are known as \emph{division algebras}. As we will see now, there are only three Banach division algebras (up to isomorphism).

\begin{theorem}[Gelfand-Mazur Theorem]\label{gelfmazur}
Suppose $\As$ is a complex Banach division algebra. Then $\As\cong\C$, i.e. $\As=\{\lam\1\mid\lam\in\C\}$.
\end{theorem}
In other words, $\C$ is the \emph{only} complex Banach division algebra (up to isomorphism).
\begin{proof}
For any $a\in\As$, we have $\sigma(a)\ne\emp$, so there exists $\lam\in\sigma(a)$. Equivalently, $a-\lam\1$ has no inverse, so by assumption, it must be zero. Thus $a=\lam\1$.
\end{proof}

\begin{example}
The quaternions\footnotemark\ $\mathbb{H}$ form a division algebra that is not isomorphic to $\C$. Consequently, there is \emph{no} norm on $\mathbb{H}$ that makes it a complex Banach algebra (at least not with the usual multiplication of quaternions).
\end{example}
\footnotetext{The set $\mathbb{H}$ of quaternions is given by $\mathbb{H}=\{a+bi+cj+dk\mid a,b,c,d\in\R\}$, with the relations $i^2=j^2=k^2=ijk=-1$.}

\ind The fundamental theorem of algebra follows easily from the previous theorem and some basic results from abstract algebra, see \cite{fta}.\\

\ind Clearly the Gelfand-Mazur theorem cannot hold for \emph{real} Banach algebras, as $\R$, $\C$ and $\mathbb{H}$ are all real Banach division algebras. Nonetheless, we have the following result:

\begin{theorem}[Gelfand-Mazur-Kaplansky Theorem]
Suppose $\As$ is a real Banach division algebra. Then $\As\cong\R$, $\C$ or $\mathbb{H}$.
\end{theorem}
See \cite[\S14, Theorem 7, Pages 73-74]{bonsall} for a proof.\\

\ind It is a well-known result in abstract algebra that the only finite-dimensional division algebra over an algebraically closed field $\K$ (e.g. $\C$) is $\K$ itself. In 1877, \href{https://mathshistory.st-andrews.ac.uk/Biographies/Frobenius/}{Ferdinand Georg Frobenius} (1849--1917) proved that the only finite-dimensional division algebras over $\R$ are $\R$, $\C$ and $\mathbb{H}$. The Gelfand-Mazur theorem was first stated in 1938 by \href{https://mathshistory.st-andrews.ac.uk/Biographies/Mazur/}{Stanisław Mazur} (1905--1981) \cite{mazur}, and proved in 1941 by Gelfand \cite{gelfand}. The Gelfand-Mazur-Kaplansky theorem was proved in 1949 by \href{https://mathshistory.st-andrews.ac.uk/Biographies/Kaplansky/}{Irving Kaplansky} (1917--2006) \cite{kaplansky}. These theorems further demonstrate the power of norms (and thus analysis) in extending results from finite-dimensional algebras to infinite-dimensional algebras.\\

\ind Both of these theorems also hold under (certain) weaker assumptions, see \cite{bhatt}.

\begin{lemma}
Suppose $\As$ is a unital algebra and $a,b\in\As$. Then $\1-ab\in\As\inv$ if and only if $\1-ba\in\As\inv$.
\end{lemma}
\begin{proof}
Suppose $\1-ab\in\As\inv$ and define $c=(\1-ab)^{-1}$. Then:
\begin{align*}
    (\1+bca)(\1-ba)&=\1-ba+bca-bcaba=\1-ba+bc(\1-ab)a=\1-ba+ba=\1 \\
    (\1-ba)(\1+bca)&=\1-ba+bca-babca=\1-ba+b(\1-ab)ca=\1-ba+ba=\1
\end{align*}
Thus $1-ba\in\As\inv$ and $(1-ba)^{-1}=1+bca$. The converse follows by switching $a$ and $b$.
\end{proof}

\begin{lemma}
Suppose $\As$ is a unital algebra and $a,b\in\As$. Then $\sigma(ab)\exc\{0\}=\sigma(ba)\exc\{0\}$.\\
In other words, $ab$ and $ba$ have the same spectrum, except possibly for $0$.
\end{lemma}
\begin{proof}
Suppose $\lam\in\sigma(ab)$ and $\lam\ne0$. Then $ab-\lam\1\notin\As\inv$, so $\1-\frac{1}{\lam}ab=-\frac{1}{\lam}(ab-\lam\1)\notin\As\inv$. Applying the previous lemma to $a$ and $\frac{1}{\lam}b$ shows that $1-\frac{1}{\lam}ba$ is also not invertible, so neither is $ba-\lam\1=-\lam\left(\1-\frac{1}{\lam}ba\right)$. Thus $\lam\in\sigma(ba)$. The converse follows by switching $a$ and $b$.
\end{proof}

\begin{theorem}\label{ab-ba/=1}
Suppose $\As$ is a unital complex Banach algebra and $\lam\ne0$. Then there are no elements $a,b\in\As$ such that $ab-ba=\lam\1$.
\end{theorem}
\begin{proof}
By \Cref{ftba}, $\sigma(ab)\ne\emp$. Suppose $\mu\in\sigma(ab)$. Then $ab-\mu\1=ba-(\mu-\lam)\1\notin\As\inv$, so $\mu-\lam\in\sigma(ba)$. Since $\lam\ne0$, at least one of $\mu$ and $\mu-\lam$ must be nonzero. Suppose $\mu\ne0$ (the other case is similar). By the previous lemma, we have $\mu\in\sigma(ba)$, so $ba-\mu\1=ab-(\mu+\lam)\1\notin\As\inv$, and so $\mu+\lam\in\sigma(ab)$. Repeating this process yields $\mu+n\lam\in\sigma(ab)$ for all $n\in\N$. This is a contradiction, since by \Cref{ftba}, $\sigma(ab)$ must be bounded.
\end{proof}

This theorem has profound implications in quantum mechanics, which we will discuss in \Cref{rebuildingqm}.

\subsection{Characters}\label{characters}
\ind We will now develop the notion of \emph{characters}. Characters are a central and fundamental concept in representation theory, and a topic of active research in abstract algebra. Here, we will only need them to formulate and prove the \hyperref[gelfandrep]{Gelfand representation theorem}, which will eventually lead us to the \hyperref[gelfandnaimark]{Gelfand-Naimark theorem}.

\begin{definition}
Suppose $\As$ is a Banach algebra. A \textbf{character} of $\As$ is a nonzero homomorphism $\chi:\As\to\F$.\\
The set of all characters of $\As$ is denoted by $\Sigma(\As)$ or $\Phi_\As$.
\end{definition}

\begin{example}
Suppose $\As=\F^n$ (with componentwise operations). Then for each $k=1,2,...,n$, the map $\chi:\As\to\F$, $\chi(a_1,a_2,...,a_n)=a_k$ is a character of $\As$.
\end{example}

\begin{warning}
Some sources use $\Sigma(\As)$ to refer to the state space of a \csta. This is unrelated to the Gelfand spectrum.
\end{warning}

\begin{proposition}\label{char1inv}
Suppose $\As$ is a unital Banach algebra and $\chi\in\Sigma(\As)$. Then:
\begin{enumerate}
    \item $\chi(\1)=1$
    \item If $a\in\As\inv$, then $\chi(a)\ne0$ and $\chi(a^{-1})=\frac{1}{\chi(a)}$.
\end{enumerate}
\end{proposition}
\begin{proof}
\begin{enumerate}
    \item For all $a\in\As$, we have $\chi(a)=\chi(a\1)=\chi(a)\chi(1)$. Since $\chi\ne0$, there exists $a\in\As$ such that $\chi(a)\ne0$. Thus $\chi(\1)=\frac{\chi(a)}{\chi(a)}=1$.
    \item By (1), we have $1=\chi(\1)=\chi(aa^{-1})=\chi(a)\chi(a^{-1})$. Thus $\chi(a)\ne0$ and $\chi(a^{-1})=\frac{1}{\chi(a)}$.\qedhere
\end{enumerate}
\end{proof}

\ind You might expect the notion of characters to lead to a rich framework within the theory of Banach algebras. In some sense, it does, as we will see shortly. However, we should not get ahead of ourselves. We first need to rule out the majority of cases where characters do not tell us anything, simply because there are none of them. In most cases where $\As$ is non-abelian, it has no characters. This is true even in the simplest examples one could think of:

\begin{example}
$M_n(\F)$ has no characters if $n\ge 2$. To see this, suppose $E_{kl}$ is the $n\times n$ matrix whose $(k,l)$-entry is $1$ and all other entries are $0$. Then $E_{k,l}E_{l,k}=E_{k,k}$, and $E_{k,l}^2=0$ if $k\ne l$. Suppose $\chi$ is a character of $M_n(\F)$. Then for all $k\ne l$, we have $\chi(E_{k,l})^2=\chi(E_{k,l}^2)=\chi(0)=0$, and so $\chi(E_{k,l})=0$. As for $E_{k,k}$, we have $\chi(E_{k,k})=\chi(E_{k,l}E_{l,k})=\chi(E_{k,l})\chi(E_{l,k})=0\cdot0=0$. Thus $\chi(E_{k,l})=0$ for all $k,l\in\{1,2,...,n\}$. Since every $A\in M_n(\F)$ is a linear combination of these matrices, we have $\chi(A)=0$ for all $A\in M_n(\F)$, and so $\chi=0$, a contradiction. Thus $M_n(\F)$ has no characters.
\end{example}

\ind Most other non-abelian Banach algebras, such as $\Bs(X)$ (where $X$ is a Banach space and $\dim(X)\ge2$), have the same problem as they contain a copy of $M_2(\F)$. As such, we will continue our discussion of characters for \emph{abelian} Banach algebras only.

\begin{theorem}\label{charbd}
Suppose $\As$ is an abelian Banach algebra and $\chi\in\Sigma(\As)$. Then $\chi$ is bounded and $\norm{\chi}=1$.
\end{theorem}
\begin{proof}
Assume without loss of generality that $\As$ is unital (otherwise, replace $\As$ with its unitization $\As_1$ defined in \Cref{embedding} and extend $\chi$ to $\As_1$ by defining $\chi(\1_{\As_1})=1$). Suppose $a\in\As$ and define $\lam=\chi(a)$. If $\abs{\lam}>\norm{a}$, then $\norm{\frac{1}{\lam}a}<1$, so $1-\frac{1}{\lam}a\in\As\inv$. Define $b=\left(\1-\frac{1}{\lam}a\right)^{-1}$. Then $\1=b(\1-\frac{1}{\lam}a)=b-\frac{1}{\lam}ba$. This yields $1=\chi(\1)=\chi(b-\frac{1}{\lam}ba)=\chi(b)-\frac{\chi(b)\chi(a)}{\lam}=\chi(b)-\chi(b)=0$, a contradiction. Thus $\abs{\lam}=\abs{\chi(a)}\le\norm{a}$ for all $a\in\As$, so $\chi$ is bounded and $\norm{\chi}\le1$. Finally, since $\chi(\1)=1$, we have $\norm{\chi}=1$.
\end{proof}

\begin{theorem}\label{kermaxideal}
Suppose $\As$ is a unital abelian complex Banach algebra and $\chi\in\Sigma(\As)$. Then $\ker(\chi)$ is a maximal ideal of $\As$. Moreover, if $\Ms$ is a maximal ideal of $\As$, then there is exactly one character $\chi\in\Sigma(\As)$ such that $\ker(\chi)=\Ms$.
\end{theorem}
In other words, if $\text{M}(\As)$ is the set of all maximal ideals of $\As$, then the map $\kappa:\Sigma(\As)\to\text{M}(\As)$ given by $\kappa(\chi)=\ker(\chi)$ is a bijection.

\begin{center}
\begin{tikzpicture}
\draw[-,very thick,red,fill=red!10] (0,0) ellipse (2 and 0.8);
\draw[-,very thick,blue,fill=blue!10] (6,0) ellipse (2 and 0.8);
\node[yshift=4pt,align=center,font=\large\sffamily] at (0,0) {$\Sigma(\As)$\\Characters of $\As$};
\node[yshift=4pt,align=center,font=\large\sffamily] at (6,0) {$\operatorname{M}(\As)$\\Maximal ideals of $\As$};
\draw[->,thick] (2,0.5) arc (120:60:2) node[midway,above]{$\chi\mapsto\ker(\chi)$};
\draw[->,thick] (4,-0.5) arc (300:240:2) node[midway,below]{$\Ms\mapsto g\circ\pi$};
\end{tikzpicture}\\
$\Sigma(\As)$ and $\operatorname{M}(\As)$ are really two sides of the same coin.
\end{center}

\begin{proof}
Suppose $\chi\in\Sigma(\As)$. Then $\ker(h)$ is a non-trivial ideal of $\As$, and since $\im(\chi)=\C$, we have $\As/\ker(\chi)\cong\C$. Thus $\ker(\chi)$ is a maximal ideal of $\As$.\\

Suppose $\Ms$ is a maximal ideal of $\As$. By \Cref{maxclosed}, $\Ms$ is closed, so by \Cref{quotientbalg}, $\As/\Ms$ is a unital Banach algebra. Suppose $\pi:\As\to\As/\Ms$ is the quotient map. If $a\in\As$ and $\pi(a)$ is not invertible in $\As/\Ms$, then $\pi(\As a)=\pi(a)[\As/\Ms]$ is a proper ideal of $\As/\Ms$. Define $\Is=\{b\in\As\mid\pi(b)\in\pi(\As a)\}=\pi^{-1}(\pi(\As a))$. Then $\Is$ is a proper ideal of $\As$ and $\Is\sups\Ms$. Since $\Ms$ is maximal, we have $\Is=\Ms$. Thus $\pi(a\As)\sub\pi(\Is)=\pi(\Ms)=\{0\}$, and so $\pi(a)=0$. In other words, if $\pi(a)$ is not invertible, it must be zero. By the Gelfand-Mazur theorem (\Cref{gelfmazur}), $\As/\Ms\cong\C$, i.e. $\As/\Ms=\{\lam u+\Ms\mid\lam\in\C\}$ for some $u\in\As$. Define $g:\As/\Ms\to\C$ by $g(\lam u+\Ms)=\lam$ and define $\chi:\As\to\C$ by $h=g\circ\pi$. Then $\chi$ is a homomorphism (as it is a composition of homomorphisms) and $\ker(\chi)=\Ms$.\\

Now suppose $\chi_1,\chi_2\in\Sigma(\As)$ are nonzero homomorphisms and $\ker(\chi_1)=\ker(\chi_2)$. Since they are linear, we have $\alpha\chi_1=\beta\chi_2$ for some $\alpha,\beta\in\C$. This yields $\alpha=\alpha\chi_1(\1)=\beta\chi_2(\1)=\beta$, so $\alpha=\beta$, and so $\chi_1=\chi_2$.
\end{proof}

\ind The last theorem only holds for \emph{complex} unital abelian Banach algebras, i.e. it requires $\F=\C$. Specifically, we used the Gelfand-Mazur theorem to conclude that $\As/\Ms\cong\C$ and thus has complex dimension $1$. In fact, if $\F=\R$, one can show that $\As/\Ms\cong\R$ or $\C$ and so has real dimension $1$ or $2$. This allows for a slight modification of the last theorem to hold, where we pass to equivalence classes of characters under complex conjugation. See \cite{ingelstam} for details.\\

\ind For the next definition and theorem, we will need the concepts of nets and the weak* topology. See \Cref{weak*top} for more details.

\begin{definition}\label{gelfspec}
Suppose $\As$ is an abelian Banach algebra. The \textbf{Gelfand spectrum} (or \textit{topological spectrum}) of $\As$ is the set $\Sigma(\As)$ of all characters of $\As$, with the weak* topology inherited from $\As^*$.
\end{definition}
\begin{remark}
Since every character of $\As$ is a bounded linear functional, we have $\Sigma(\As)\subset\As^*$.
\end{remark}

\begin{theorem}\label{gelfspecch}
Suppose $\As$ is an abelian Banach algebra. Then $\Sigma(\As)$ is a locally compact Hausdorff space. If $\As$ is unital, then $\Sigma(\As)$ is compact.
\end{theorem}
\begin{proof}
Since $\As^*$ is Hausdorff, so is $\Sigma(\As)$ (since every subset of a Hausdorff space is Hausdorff). Suppose $B$ is the closed unit ball in $\As^*$. By the \hyperref[banachalaoglu]{Banach-Alaoglu theorem}, $B$ is weak*-compact. We will denote the set of \emph{all} homomorphisms $\chi:\As\to\F$ (including the zero homomorphism) by $\overline{\Sigma}(\As)$. Then $\overline{\Sigma}(\As)\sub B$. We will show that $\overline{\Sigma}(\As)$ is weak*-closed.\\

Suppose $(\chi_\alpha)$ is a net in $\overline{\Sigma}(\As)$ that weak*-converges to some $\chi\in B$. Then for all $a,b\in\As$, we have $\chi(ab)=\lim_\alpha \chi_\alpha(ab)=\lim_\alpha \chi_\alpha(a)\chi_\alpha(b)=\chi(a)\chi(b)$. Similarly, $\chi(\lam a+b)=\lam\chi(a)+\chi(b)$ for all $a,b\in\As$ and all $\lam\in\F$. Thus $\chi$ is a homomorphism, i.e. $\chi\in\overline{\Sigma}(\As)$. Thus $\overline{\Sigma}(\As)$ is weak*-closed, and so it is weak*-compact (since every closed subset of a compact set is compact).\\

Since $\As^*$ is Hausdorff, the singleton $\{0\}$ is closed in $\As^*$ (and thus in $\overline{\Sigma}(\As)$), and so its complement $\Sigma(\As)=\overline{\Sigma}(\As)\exc\{0\}$ is open in $\overline{\Sigma}(\As)$. Thus $\Sigma(\As)$ is locally compact (since every open subset of a compact set is locally compact).\\

Now suppose $\As$ is unital. Then for any net $(\chi_\alpha)$ in $\Sigma(\As)$ that converges to $\chi\in\Sigma(\As)$, we have $\chi(\1)=\lim_\alpha \chi_\alpha(\1)=\lim_\alpha 1=1\ne0$, so $\chi\ne0$, and so $\chi\in\Sigma(\As)$. Thus $\Sigma(\As)$ is a weak*-closed subset of $B$, and so it is weak*-compact.
\end{proof}
\begin{remarks}\leavevmode
\begin{enumerate}
    \item This proof also shows that $\overline{\Sigma}(\As)$ is the one-point compactification of $\Sigma(\As)$ (since it is the closure of $\Sigma(\As)$ in $\As^*$ and contains exactly one more point). See \cite[Pages 117-118]{singh} for more details.
    \item If $\As$ is not unital, then the set $\overline{\Sigma}(\As)$ of \textit{all} homomorphisms $\chi:\As\to\F$ is still weak*-closed (and thus weak*-compact), but there may be nets of nonzero homomorphisms that converge to the zero homomorphism, as the next example shows. If $\As$ is unital, this is not possible, as they all have to satisfy $\chi(\1)=1$.
\end{enumerate}
\end{remarks}

\begin{example}
Suppose $\As=C_0(\R)$ and define the sequence $(\phi_n)$ by $\phi_n:C_0(\R)\to\C$, $\phi_n(f)=f(n)$. Then each $\phi_n$ is a nonzero homomorphism. However, for all $f\in C_0(\R)$, we have $\limn\phi_n(f)=\limn f(n)=0$, so $(\phi_n)$ weak*-converges to $0$.
\end{example}

\begin{theorem}\label{gelfspec=spec}
Suppose $\As$ is a unital abelian complex Banach algebra. Then for all $a\in\As$, we have:
\begin{equation*}
    \sigma(a)=\left\{\chi(a)\mid\chi\in\Sigma(\As)\right\}
\end{equation*}
In other words, the set of values attained by the characters $\chi\in\Sigma(\As)$ at $a$ is exactly $\sigma(a)$.
\end{theorem}

This explains the name ``Gelfand spectrum'' and the notation $\Sigma(\As)$, since for any $a\in\As$, the elements of $\Sigma(\As)$ can be evaluated at $a$ to recover its spectrum $\sigma(a)$.

\begin{proof}
Suppose $\chi\in\Sigma(\As)$ and define $\lam=\chi(a)$. Then $\chi(a-\lam\1)=\chi(a)-\lam\chi(\1)=\chi(a)-\chi(a)=0$, so $a-\lam\1\in\ker(\chi)$. Thus $a-\lam\1\notin\As\inv$, and so $\lam\in\sigma(a)$.\\

Now suppose $\lam\in\sigma(a)$. Then $a-\lam\1$ is not invertible. Thus $(a-\lam\1)\As$ is a proper ideal of $\As$. Suppose $\Ms$ is a maximal ideal of $\As$ containing $(a-\lam\1)\As$. By \Cref{kermaxideal}, there is a character $\chi\in\Sigma(\As)$ such that $\ker(\chi)=\Ms$. This yields $0=\chi(a-\lam\1)=\chi(a)-\lam\chi(1)=\chi(a)-\lam$. Thus $\chi(a)=\lam\in\sigma(a)$. Thus $\sigma(a)=\{\chi(a)\mid\chi\in\Sigma(\As)\}$.
\end{proof}

\begin{definition}
Suppose $\As$ is an abelian Banach algebra and $a\in\As$. The \textbf{Gelfand transform} of $a$ is the function $\ahat:\Sigma(\As)\to\F$ given by $\ahat(h)=h(a)$.\\
The \textbf{Gelfand transform} of $\As$ is the function $\gamma:\As\to C_0(\Sigma(\As))$ given by $\gamma(a)=\ahat$.
\end{definition}
\begin{remarks}\leavevmode
\begin{enumerate}
    \item We do not know \emph{a priori} that $\ahat\in C_0(\Sigma(\As))$, but we will show this in the proof of the next theorem.
    \item Some sources refer to the Gelfand transform (both of an individual element and of the whole algebra) as the \emph{Gelfand representation}.
\end{enumerate}
\end{remarks}

\begin{theorem}[Gelfand Representation Theorem]\label{gelfandrep}
Suppose $\As$ is an abelian Banach algebra. Then for all $a\in\As$, we have $\ahat\in C_0(\Sigma(\As))$. Moreover, the Gelfand transform of $\As$ is a continuous homomorphism and $\norm{\gamma}\le1$, and for all $a\in\As$, we have $\norm{\ahat}_\infty=r(a)$. Furthermore, if $\As$ is unital and complex, then $\ker(\gamma)$ is the intersection of all maximal ideals\footnotemark\ of $\As$.
\end{theorem}
\vspace*{-10pt}
\footnotetext{This is also known as the \emph{Jacobson radical} of $\As$, though this terminology is rarely used outside of ring theory.}
\begin{proof}
Suppose $\chi\in\Sigma(\As)$ and $(\chi_\alpha)$ is a net in $\Sigma(\As)$ such that $\chi_\alpha\to\chi$ in $\Sigma(\As)$. Then $\chi_\alpha\to\chi$ weak* in $\As^*$. Thus for all $a\in\As$, we have $\ahat(\chi_\alpha)=\chi_\alpha(a)\to\chi(a)=\ahat(\chi)$. Thus $\ahat$ is continuous. Moreover, we can extend $\ahat$ to a continuous function on $\overline{\Sigma}(\As)$ by defining\footnotemark\ $\ahat(0)=0$. Thus $\ahat$ vanishes at infinity, and so $\ahat\in C_0(\Sigma(\As))$.\\

We now show that $\gamma$ is a homomorphism. Suppose $a,b\in\As$ and $\lam\in\F$. Then for all $\chi\in\Sigma(\As)$, we have:
\begin{equation*}
    \gamma(\lam a+b)(\chi)=\widehat{\lam a+b}(\chi)=\chi(\lam a+b)=\lam\chi(a)+\chi(b)=\lam\ahat(\chi)+\bhat(\chi)=\lam\gamma(a)(\chi)+\gamma(b)(\chi)
\end{equation*}
Thus $\gamma(\lam a+b)=\lam\gamma(a)+\gamma(b)$, and so $\gamma$ is linear. Now suppose $a,b\in\As$. Then for all $\chi\in\Sigma(\As)$, we have:
\begin{equation*}
    \gamma(ab)(\chi)=\widehat{ab}(\chi)=\chi(ab)=\chi(a)\chi(b)=\ahat(\chi)\bhat(\chi)=\gamma(a)(\chi)\gamma(b)(\chi)
\end{equation*}
Thus $\gamma(ab)=\gamma(a)\gamma(b)$, and so $\gamma$ is a homomorphism.\\

For all $a\in\As$, we have $\abs{\ahat(\chi)}=\abs{\chi(a)}\le\norm{a}$ (since $\norm{\chi}=1$), so $\norm{\gamma(a)}_\infty=\norm{\ahat}_\infty\le\norm{a}$. Thus $\gamma$ is bounded and $\norm{\gamma}\le1$. Also, for all $a\in\As$, we have:
\begin{equation*}
    \norm{\ahat}_\infty=\sup_{\chi\in\Sigma(\As)} \abs{\chi(a)}=\sup_{\lam\in\sigma(a)} \abs{\lam}=r(a)
\end{equation*}

Finally, if $\As$ is unital and complex, we have:
\begin{align*}
    a\in\ker(\gamma)
    &\iff\gamma(a)=\ahat\equiv0 \\
    &\iff h(a)=0 \text{ for all } h\in\Sigma(\As) \\
    &\iff a\in\ker(h) \text{ for all } h\in\Sigma(\As) \\
    &\iff a\in\Ms \text{ for every maximal ideal } \Ms \text{ of } \As \tag{by \Cref{kermaxideal}} \\
    &\iff a\in\capp_{\substack{\Ms\text{ is a maximal}\\\text{ideal of } \As}} \Ms
\end{align*}
Thus $\ker(\gamma)$ is the intersection of all maximal ideals of $\As$.
\end{proof}
\vspace*{-10pt}
\footnotetext{In other words, we extend $\ahat$ from $\Sigma(\As)$ (which corresponds to $\As$) to its closure $\overline{\Sigma}(\As)$ (which corresponds to the unitization $\As_1$ of $\As$, see \Cref{embedding}). Since the new `point' is the zero homomorphism, the only suitable value of $\ahat(0)$ is $0$.}
\begin{remark}
If $\As$ is unital, it also follows directly that $\norm{\gamma}=1$, since $\gamma(\1)(\chi)=\widehat{\1}(\chi)=\chi(\1)=1$ for all $\chi\in\Sigma(\As)$.
\end{remark}

\ind The correspondence between characters and maximal ideals in \Cref{kermaxideal} is used to show that $\ker(\gamma)$ is the intersection of all maximal ideals of $\As$. This is the only part of the proof that requires $\As$ to be unital and complex. However, even if $\As$ is real or non-unital, the Gelfand transform still gives a meaningful representation of $\Sigma(\As)$ into $C_0(\Sigma(\As))$. See \cite{albiac} for related results.

\begin{example}
Suppose $\As=L^1(\R)$ (with the convolution product). The characters of $L^1(\R)$ are point evaluations of the Fourier transform, i.e. maps that send $f\in L^1(\R)$ to $\fh(\eta)=\int_{-\infty}^\infty e^{-2\pi i\eta x}f(x)\dx$ for fixed $\eta\in\R$. Thus each $\chi\in\Sigma(L^1(\R))$ corresponds to exactly one $\eta\in\R$, and the Gelfand transform of $f\in L^1(\R)$ is $\fh$. In other words, the Gelfand transform of $L^1(\R)$ is simply the Fourier transform.
\end{example}

\subsection{The Exponential}\label{exponential}
\ind In calculus, one learns about the all-important power series formula for the natural exponential: $e^x=\sumn[0]\frac{1}{n!}x^n$. In linear algebra (or perhaps a first course on differential equations), one learns how this formula can be adapted to define the \emph{matrix exponential}:
\begin{equation*}
    e^A=\sumn[0]\frac{1}{n!}A^n
\end{equation*}
\ind The first remarkable fact about this is that it is always well-defined (as long as we use the convention $A^0=I$, the identity matrix). The second is how it can be used to solve systems of linear differential equations (or linear difference equations, or several other types of equations). But this is merely a special case of a \emph{much} more general type of exponential, which we will discuss here.

\begin{definition}
Suppose $\As$ is a unital Banach algebra and $a\in\As$. The \textbf{exponential} of $a$, denoted by $\exp(a)$ or $e^a$, is given by:
\begin{equation*}
    \exp(a)=\sumn[0]\frac{1}{n!}a^n
\end{equation*}
Where we define $a^0=\1$, the identity of $\As$.
\end{definition}
\begin{remark}
This is well-defined for all $a\in\As$, as it is a power series in $a$ and its radius of convergence is $\infty$.
\end{remark}

\begin{example}
Suppose $\As=M_2(\C)$ and $A=\mmat{1}{5}{0}{2}$. Then for all $n\in\N_0$, we have $A^n=\mmat{1}{5(2^n-1)}{0}{2^n}$. This yields:
\begin{equation*}
    \exp(A)
    =\sumn[0]\frac{1}{n!}A^n
    =\sumn[0]\frac{1}{n!}\mmat{1}{5(2^n-1)}{0}{2^n}
    =\mqty(\sumn[0]\frac{1}{n!}&5\sumn[0]\frac{2^n-1}{n!}\\[4pt]0&\sumn[0]\frac{2^n}{n!})
    =\mmat{e}{5(e^2-e)}{0}{e^2}
\end{equation*}
\end{example}

\begin{proposition}
Suppose $\As$ is a unital Banach algebra and $a,b\in\As$. Then:
\begin{enumerate}
    \item $\norm{\exp(a)}\le e^\norm{a}$
    \item If $a$ and $b$ commute, then $\exp(a+b)=\exp(a)\exp(b)$.
    \item $\exp(a)\in\As\inv$ and $(\exp(a))^{-1}=\exp(-a)$.
\end{enumerate}
\end{proposition}
\begin{proof}\leavevmode
\begin{enumerate}
    \item $\Disp\norm{\exp(a)}=\norm{\sumn[0]\frac{1}{n!}a^n}\le\sumn[0]\frac{1}{n!}\norm{a^n}\le\sumn[0]\frac{1}{n!}\norm{a}^n=e^\norm{a}$
    \item Since $a$ and $b$ commute, we can use the binomial theorem to rewrite $\exp(a+b)$ as follows:
    \begin{align*}
        \exp(a+b)=\sumn[0]\frac{1}{n!}(a+b)^n=\sumn[0]\frac{1}{n!}\sum_{k=0}^n \binom{n}{k}a^{n-k}b^k=\sumn[0]\sum_{k=0}^n \frac{1}{(n-k)!}a^{n-k}\frac{1}{k!}b^k
    \end{align*}
    Substituting $r=n-k$ and $s=k$, we get:
    \begin{align*}
        \exp(a+b)=\sum_{r=0}^\infty \sum_{s=0}^\infty \frac{1}{r!}a^r\frac{1}{s!}b^s=\left(\sum_{r=0}^\infty \frac{1}{r!}a^r\right)\left(\sum_{s=0}^\infty\frac{1}{s!}b^s\right)=\exp(a)\exp(b)
    \end{align*}
    \item Since $a$ and $-a$ commute, by (2), we have $\exp(a-a)=\exp(a)\exp(-a)$ and $\exp(a-a)=\exp(-a)\exp(a)$. Since $\exp(a-a)=\exp(0)=\1$, we have that $\exp(a)\in\As\inv$ and $\exp(a)^{-1}=\exp(-a)$.\qedhere
\end{enumerate}
\end{proof}

Equation (2) above does not hold in general, as the following example shows.
\begin{example}
Suppose $A,B\in M_2(\C)$, $A=\mmat{1}{0}{0}{0}$ and $B=\mmat{0}{1}{0}{0}$. Then $\exp(A)=\mmat{e}{0}{0}{1}$, $\exp(B)=\mmat{1}{1}{0}{1}$ and $A+B=\mmat{1}{1}{0}{0}$. This yields $\exp(A+B)=\mmat{e}{e-1}{0}{1}$, which is not equal to $\exp(A)\exp(B)=\mmat{e}{e}{0}{1}$ or $\exp(B)\exp(A)=\mmat{e}{1}{0}{1}$. This is possible because $A$ and $B$ do not commute.
\end{example}

\ind While $\exp(a+b)\ne\exp(a)\exp(b)$ in general, there are several formulas that can be used to mitigate this problem, such as the Lie product formula (see \cite[Corollary 14.4, Page 208]{kemp}) and the Baker-Campbell-Hausdorff formula (see \cite[Theorem 14.21, Page 218]{kemp}).\\

\ind This idea of extending the exponential function to Banach algebras using its series definition can be applied with other analytic functions, such as $\sin(x)=\sumn[0] \frac{(-1)^n}{(2n+1)!}x^{2n+1}$. This leads to the theory of functional calculus, which is discussed in \cite[Chapter 12]{ward}, \cite[Chapter 10]{grudin} and \cite{verstraten}.

{
\definecolor{col1}{HTML}{98FCDB}
\definecolor{col2}{HTML}{9CDAF8}
\begin{tcolorbox}[interior style={left color=col1!50,right color=col2!50},frame style={left color=col1!70!black,right color=col2!70!black},fontupper=\sffamily,width=0.75\linewidth,center,top=2mm,bottom=2mm,left=2mm,right=2mm,arc=4mm,before skip=10pt,after skip=10pt,enhanced]
The rest of this chapter is a side quest, intended to lead up to the GKZ theorem. It will not be required for our later discussion of \cstas\ and quantum mechanics. You could skip to the next chapter if you \emph{really} want, but where's the fun in that?
\end{tcolorbox}
}

\begin{lemma}
Suppose $\As$ is a unital complex Banach algebra and $a,b\in\As$. Also suppose there is a constant $M>0$ such that $\norm{\exp(\lam a)b\exp(-\lam a)}\le M$ for all $\lam\in\C$. Then $ab=ba$.
\end{lemma}
Basically, given any $a,b\in\As$, either $b$ commutes with $a$, or its norm can be made arbitrarily large by conjugating it with $\exp(\lam a)$ for some sufficiently large $\lam$.
\begin{proof}
Define $f:\C\to\As$ by $f(\lam)=\exp(\lam a)b\exp(-\lam a)$. Since the exponential function is entire, $f$ is a product of entire functions, and so it is entire. By assumption, it is also bounded, so by \Cref{liouvillenvs}, it is constant. Thus $f(\lam)=f(0)=b$ for all $\lam\in\C$. Expanding $f$ as a power series in $\lam$, we get:
\begin{align*}
    f(\lam)&=(1+\lam a+O(\lam^2))b(1-\lam a+O(\lam^2))=(b+\lam ab+O(\lam^2))(1-\lam a+O(\lam^2))=b+\lam ab-\lam ba+O(\lam^2)\\
    &=b+\lam(ab-ba)+O(\lam^2)
\end{align*}
Thus $ab=ba$.
\end{proof}
\begin{remark}
The higher-order terms in the above expansion do not tell us anything, as they all automatically vanish if $ab=ba$.
\end{remark}

The following example shows that this lemma does not hold for \textit{real} Banach algebras.
\begin{example}
Suppose $\As=M_2(\R)$ and define $A=\mmat{0}{-1}{1}{0}$ and $B=\mmat{0}{1}{0}{0}$. Then $AB=\mmat{0}{0}{0}{1}$ but $BA=\mmat{1}{0}{0}{0}$, so $A$ and $B$ do not commute. However, for all $\lam\in\R$, we have $\exp(\lam A)=\mmat{\cos(\lam)}{-\sin(\lam)}{\sin(\lam)}{\cos(\lam)}$, and so $\norm{\exp(\lam A)}=1$. Thus $\norm{\exp(\lam A)B\exp(-\lam A)}\le\norm{\exp(\lam A)}\norm{B}\norm{\exp(-\lam A)}=1\cdot 1\cdot 1=1$ for all $\lam\in\R$, even though $A$ and $B$ do not commute. Note that if we allow $\lam\in\C$, then for $\lam=it$ where $t\in\R$, we have $\norm{\exp(\lam A)B\exp(-\lam A)}=\cosh(2t)$. Since $\cosh(2t)\to\infty$ as $t\to\pm\infty$, this is indeed unbounded.
\end{example}

\begin{corollary}
Suppose $\As$ is a unital non-abelian complex algebra. Then there is NO norm on $\As$ that makes it a unital complex Banach algebra and is preserved by similarity transformations, i.e. transformations of the form $b\mapsto aba^{-1}$.
\end{corollary}
\begin{proof}
Suppose $\norm{\cdot}$ is such a norm and $a,b\in\As$. Since $\norm{\cdot}$ is preserved by similarity transformations, we have $\norm{\exp(\lam a)b\exp(-\lam a)}=\norm{\exp(\lam a)b(\exp(\lam a))^{-1}}=\norm{b}$ for all $\lam\in\C$. By the previous lemma, we have $ab=ba$, and so $\As$ is abelian.
\end{proof}

\begin{example}
$M_n(\C)$ is a unital complex algebra, and it is non-abelian if $n\ge 2$. Thus there is no matrix norm that is preserved by similarity transformations. The standard operator norm is preserved by \emph{unitary} similarity transformations, but not by arbitrary ones.
\end{example}

To prove the next lemma, we will need the following result from complex analysis:

\begin{theorem}[Borel-Carathéodory theorem]\label{borelcaratheodory}
Suppose $g:\Dbar(0,R)\to\C$ is holomorphic on $\Dbar(0,R)$, the closed disk of radius $R$ centered at $0$. Then for any $r<R$, we have:
\begin{equation*}
    \sup_{\abs{z}\le r} \abs{g(z)}
    \le\frac{2r}{R-r}\sup_{\abs{z}\le R} \Re(g(z))+\frac{R+r}{R-r}\abs{g(0)}
\end{equation*}
\end{theorem}
See \cite[Theorem 8.3.3, Pages 258-259]{taylor} for a proof.

\begin{lemma}\label{gkzentire}
Suppose $f:\C\to\C$ is an entire function such that $f(0)=1$, $f'(0)=0$ and there is a constant $M>0$ such that $0<\abs{f(z)}<e^{M\abs{z}}$ for all $z\in\C$. Then $f(z)=1$ for all $z\in\C$.
\end{lemma}
\begin{proof}
Since $f$ is entire and never zero, there is an entire function $g$ such that $f(z)=e^{g(z)}$. By assumption, we have $1=f(0)=e^{g(0)}$, so we can assume $g(0)=0$ (by adding an integer multiple of $2\pi i$ to $g$). We also have $0=f'(0)=g'(0)e^{g(0)}$, so $g'(0)=0$. Note that $e^{\Re(g(z))}=\abs{f(z)}<e^{M\abs{z}}$, so $\Re(g(z))<M\abs{z}$ for all $z\in\C$. Suppose $0<r<R$. By the \hyperref[borelcaratheodory]{Borel-Carathéodory theorem}, for all $\abs{z}\le r$, we have:
\begin{equation*}
    \abs{g(z)}
    \le\frac{2r}{R-r}\sup_{\abs{z}\le R} \Re(g(z))+\frac{R+r}{R-r}\abs{g(0)}
    \le\frac{2r}{R-r}MR+0
    =\frac{2MrR}{R-r}
\end{equation*}
Setting $R=2r$ yields $\abs{g(z)}\le 4mr$ for all $\abs{z}\le r$. In other words, $\abs{g(z)}$ can be bounded above on $\Dbar(0,r)$ by $4mr$, which is linear in $r$. Thus $g$ is either constant or linear in $z$. Since $g(0)=g'(0)=0$, it follows that $g\equiv0$. Thus $f\equiv1$.
\end{proof}

\begin{lemma}\label{gkzequiv}
Suppose $\As$ is a unital algebra, $\phi:\As\to\F$ is linear and $\phi(\1)=1$. Then the following are equivalent:
\begin{enumerate}
    \item If $a\in\As$ and $\phi(a)=0$, then $\phi(a^2)=0$.
    \item If $a\in\As$, then $\phi(a^2)=\phi(a)^2$.
    \item If $a,b\in\As$ and $\phi(a)=0$, then $\phi(ab)=0$.
    \item If $a,b\in\As$, then $\phi(ab)=\phi(a)\phi(b)$.
\end{enumerate}
\end{lemma}
\begin{proof}\leavevmode
\begin{enumerate}[leftmargin=.55in]
    \item[($1\Rightarrow2$)] Suppose $a\in\As$. Since $\phi$ is linear and $\phi(\1)=1$, we have $\phi(a-\phi(a)\1)=\phi(a)-\phi(a)\phi(\1)=\phi(a)-\phi(a)=0$. By (1), we have $\phi((a-\phi(a)\1)^2)=0$. Expanding this, we get:
    \begin{align*}
        0&=\phi\left((a-\phi(a)\1)^2\right)=\phi\left(a^2-2\phi(a)a+\phi(a)^2\1\right)=\phi(a^2)-2\phi(a)\phi(a)+\phi(a)^2\phi(\1)\\
        &=\phi(a^2)-2\phi(a)^2+\phi(a)^2=\phi(a^2)-\phi(a)^2
    \end{align*}
    Thus $\phi(a^2)=\phi(a)^2$.
    \item[($2\Rightarrow3$)] We first show that if $\phi(a)=0$, then $\phi(ab+ba)=0$. Replacing $a$ with $a+b$ in (2) yields $\phi((a+b)^2)=\phi(a+b)^2$. Expanding both sides, we get:
    \begin{align*}
        \phi(a)^2+\phi(ab+ba)+\phi(b)^2=\phi(a)^2+2\phi(a)\phi(b)+\phi(b)^2 \quad\therefore\quad \phi(ab+ba)=2\phi(a)\phi(b)
    \end{align*}
    Thus, if $\phi(a)=0$, then $\phi(ab+ba)=0$. By (2), this implies $\phi((ab+ba)^2)=0^2=0$. Since this holds for all $b\in\As$, we can replace $b$ with $bab$ to get $\phi(a(bab)+(bab)a)=0$. Note that $(ab+ba)^2+(ab-ba)^2=2(a(bab)+(bab)a)$. This yields:
    \begin{align*}
        \phi((ab+ba)^2)+\phi((ab-ba)^2)=2\phi(a(bab)+(bab)a) \quad\therefore\quad \phi((ab-ba)^2)=0
    \end{align*}
    By (2), we have $\phi(ab-ba)^2=0$, and so $\phi(ab-ba)=0$. Finally, we have $ab=\frac{1}{2}((ab+ba)+(ab-ba))$, and so $\phi(ab)=\frac{1}{2}(\phi(ab+ba)+\phi(ab-ba))=\frac{1}{2}(0+0)=0$.
    \item[($3\Rightarrow4$)] Suppose $a,b\in\As$. As in the proof of $1\Rightarrow2$, we have $\phi(a-\phi(a)\1)=0$. By (3), we have $\phi((a-\phi(a)\1)b)=0$. Expanding the left side, we get:
    \begin{align*}
        0=\phi((a-\phi(a)\1)b)=\phi(ab-\phi(a)b)=\phi(ab)-\phi(a)\phi(b)
    \end{align*}
    Thus $\phi(ab)=\phi(a)\phi(b)$.
    \item[($4\Rightarrow1$)] This follows trivially by setting $a=b$.\qedhere
\end{enumerate}
\end{proof}
\begin{remark}
The condition $\phi(a)=0\Rightarrow\phi(ba)=0$ is also equivalent to the above conditions, which can be seen by using $ba=\frac{1}{2}((ab+ba)-(ab-ba))$ at the end of the proof of $2\Rightarrow3$.
\end{remark}

\ind We are now ready to prove the Gleason-Kahane-Żelazko (GKZ) theorem, a remarkable theorem about characters of unital complex Banach algebras. We know from \Cref{char1inv} that for a linear functional on such an algebra to be a character, at the very least, it must map the identity element to $1$ and invertible elements to nonzero numbers. The GKZ theorem states that this is enough: \emph{every} linear functional with these two properties is a character!

\begin{center}
\begin{minipage}{16em}\begin{center}
\personbox{Fuchsia!50}{\href{https://en.wikipedia.org/wiki/Andrew_M._Gleason}{Andrew M. Gleason}\\(1921--2008)}{\includegraphics[width=1.3in]{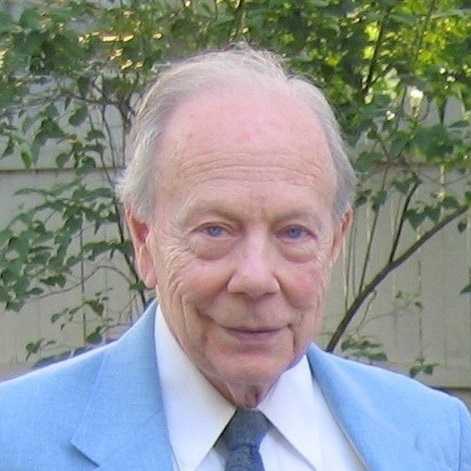}}\label{xgleason} 
\end{center}\end{minipage}
\begin{minipage}{16em}\begin{center}
\personbox{Emerald!50}{\href{https://mathshistory.st-andrews.ac.uk/Biographies/Kahane/}{Jean-Pierre Kahane}\\(1926--2017)}{\includegraphics[width=1.3in]{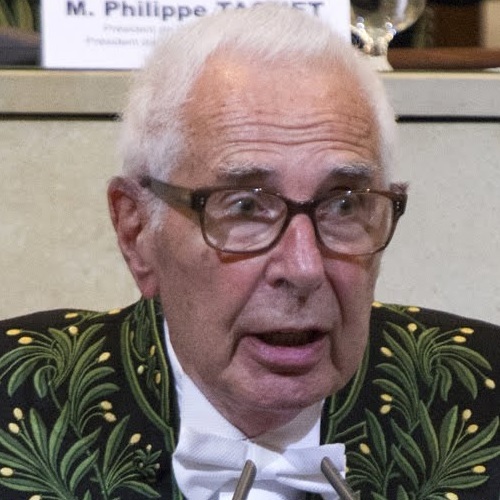}}\label{xkahane} 
\end{center}\end{minipage}
\begin{minipage}{16em}\begin{center}
\personbox{RoyalBlue!50}{\href{https://www.mathgenealogy.org/id.php?id=86239}{Wiesław Żelazko}\\(1933--)}{\includegraphics[width=1.3in]{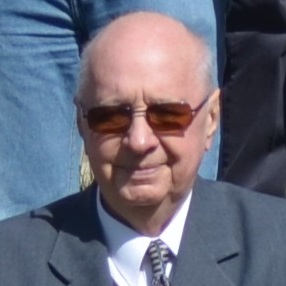}}\label{xzelazko} 
\end{center}\end{minipage}
\end{center}

\ind The GKZ theorem was first proved for unital \textit{abelian} complex Banach algebras in 1967 by Gleason\footnote{In addition to being a mathematician, Gleason was also a codebreaker in the US Navy in World War II and the Korean War. \href{https://www.britannica.com/biography/Leslie-Lamport}{Leslie Lamport} (1941--), the inventor of \LaTeX, was a second-generation student of Gleason.} \cite{gleason}, and then in 1968 by Kahane and Żelazko \cite{kahane}. The general (non-abelian) case was proved by Żelazko\footnote{Żelazko was a student of Mazur (known for the \hyperref[gelfmazur]{Gelfand-Mazur theorem}).} in 1968 \cite{zelazko}.

\begin{theorem}[Gleason-Kahane-Żelazko Theorem or GKZ Theorem]\label{gkz}
Suppose $\As$ is a unital complex Banach algebra and $\phi:\As\to\C$ is linear. Then $\phi\in\Sigma(\As)$ if and only if $\phi(\1)=1$ and $\phi(a)\ne0$ for all $a\in\As\inv$.
\end{theorem}
\begin{proof}
($\Rightarrow$) This follows directly from \Cref{char1inv}.\\
($\Leftarrow$) Suppose $\phi(\1)=1$ and $\phi(a)\ne0$ for all $a\in\As\inv$. We first show that $\phi$ is bounded. Suppose $a\in\As$. If $\phi(a)=0$, then $\abs{\phi(a)}\le\norm{a}$ holds trivially. If $\phi(a)\ne0$, then $\phi\left(\1-\frac{a}{\phi(a)}\right)=\phi(\1)-\frac{\phi(a)}{\phi(a)}=1-1=0$, so $\1-\frac{a}{\phi(a)}\in\ker(\phi)$. Thus $\1-\frac{a}{\phi(a)}\notin\As\inv$, so by \Cref{neumann}, we must have $\norm{\frac{a}{\phi(a)}}\ge1$, i.e. $\norm{\phi(a)}\le\norm{a}$. Thus $\phi$ is bounded.\\

Now suppose $a\in\As$ and $\phi(a)=0$. We will show that $\phi(a^2)=0$. Assume without loss of generality that $\norm{a}\le1$ (otherwise, replace $a$ with $\frac{a}{\norm{a}}$). Define $f:\C\to\C$ by $f(z)=\sumn[0] \frac{\phi(a^n)}{n!}z^n$. Since $\abs{\phi(a^n)}\le\norm{a^n}\le\norm{a}^n\le1$ for all $n\in\N_0$, it follows that $f$ is entire and $\abs{f(z)}\le e^{\norm{a}\abs{z}}$. We also have $f(0)=\phi(\1)=1$ and $f'(0)=\phi(a)=0$. Since $\phi$ is bounded (and thus continuous), we have:
\begin{equation*}
    f(z)=\sumn[0] \frac{\phi((za)^n)}{n!}=\phi\left(\sumn[0] \frac{(za)^n}{n!}\right)=\phi(\exp(za))
\end{equation*}
Since $\exp(za)\in\As\inv$, we have $f(z)\ne0$ for all $z\in\C$, i.e. $\abs{f(z)}>0$ for all $z\in\C$. By \Cref{gkzentire}, we have $f(z)=1$ for all $z\in\C$. Thus all terms in the definition of $f$ with $n\ge1$ vanish. In particular, $\phi(a^2)=0$.\\

We have just shown that if $a\in\As$ and $\phi(a)=0$, then $\phi(a^2)=0$. By \Cref{gkzequiv}, we have $\phi(ab)=\phi(a)\phi(b)$ for all $a,b\in\As$. Thus $\phi$ is a homomorphism (and it is nonzero by assumption), and so $\phi\in\Sigma(\As)$.
\end{proof}

\begin{warning}
The GKZ theorem is NOT true for general linear operators $\phi:\As\to\Bs$, i.e. if $\phi:\As\to\Bs$ is linear and maps invertible elements to invertible elements, it may not be multiplicative.
\end{warning}

\begin{example}
Suppose $\As=M_n(\C)$ and $\phi:\As\to\As$ is given by $\phi(A)=A^\text{T}$. Then $\phi$ is linear and maps invertible matrices to invertible matrices. However, it is not multiplicative as $\phi(AB)=(AB)^\text{T}=B^\text{T}A^\text{T}\ne A^\text{T}B^\text{T}$.
\end{example}

\begin{warning}
The GKZ theorem is NOT true for \emph{real} Banach algebras, as the next example shows.
\end{warning}

\begin{example}
Suppose $\As=C([0,1],\R)$ and define $\phi:\As\to\R$ by $\phi(f)=\int_0^1 f(x)\dx$. Then $\phi(\1)=1$. Also, if $f\in\As\inv$, it is continuous and does not vanish, so it is always positive or always negative, and so its integral cannot be zero. Thus $\phi(f)\ne0$ for every $f\in\As\inv$. However, $\phi$ is not multiplicative. On the other hand, if $\As=C([0,1],\C)$, then $f(x)=e^{2\pi ix}$ is invertible, but $\phi(f)=0$.
\end{example}

\ind There is a modification of the GKZ theorem for real Banach algebras due to S. H. Kulkarni (1984) \cite{kulkarni}: $\phi\in\Sigma(\As)$ if and only if $\phi(\1)=1$ and $\phi(a)^2+\phi(b)^2\ne0$ for all $a,b\in\As$ such that $ab=ba$ and $a^2+b^2\in\As\inv$.\\

\ind There is also a generalization to topological vector spaces, due to A. Golbaharan (2020) \cite{golbaharan}.

\conclbox{Chapter 2 Conclusion} 

\ind The results we have proved thus far show the power of combining algebra and analysis: Using the \hyperref[neumann]{Neumann series} and tools from complex analysis (specifically Liouville's theorem and Laurent's theorem), we proved the \hyperref[ftba]{fundamental theorem of Banach algebras} and the \hyperref[specradius]{spectral radius formula}. We then went on to prove the \hyperref[gelfmazur]{Gelfand-Mazur theorem} and established the correspondence between characters and maximal ideals in \Cref{kermaxideal}. Finally, we proved the \hyperref[gelfandrep]{Gelfand representation theorem}, thereby linking all abelian Banach algebras to the example of $C_0(X)$ through the Gelfand transform.\\

\ind In general, however, the Gelfand transform may not be injective or surjective. Indeed, there are abelian Banach algebras that cannot be realized as $C_0(X)$, such as those with the trivial multiplication operation $ab=0$ for all $a,b\in\As$. In the next chapter, we will explore a special class of Banach algebras, known as \cstas, that are much more well-behaved. In particular, we will show that for abelian \cstas, the Gelfand transform is bijective, and as such, all abelian \cstas\ can be realized as $C_0(X)$.

\section{$C^*$-algebras}\label[chapter]{cstas}
\ind Having developed the rich and insightful theory of Banach algebras in the previous chapter, we are now ready to move on to \cstas. Essentially, \cstas\ are Banach algebras with an additional operation * that works like the conjugate transpose of a matrix or the conjugate of a complex-valued function, that satisfies the all-important \emph{\cst-axiom}:
\begin{equation*}
    \norm{a^*a}=\norm{a}^2
\end{equation*}

\vspace*{3mm}

\begin{minipage}{\dimexpr\linewidth-1.5in} 
\ind \cstas\ were first defined in 1946 by Rickart \cite{rickart1946}. He originally called them \emph{$B^*$-algebras} (the $B$ stands for `Banach'), and the condition $\norm{a^*a}=\norm{a}^2$ was called the \emph{$B^*$-axiom}. In 1947, Segal \cite{segal1947} defined a \csta\ as a closed *-subalgebra of $\Bs(\Hs)$, where $\Hs$ is a Hilbert space. As we will see, these definitions turn out to be essentially equivalent.\\

\ind The $C$ in \csta\ stands for ``closed'', not ``complex''. It is perfectly fine to consider \cstas\ over $\R$, but in applications (primarily in quantum mechanics), all of our vector spaces are over $\C$, so we are naturally more interested in \cstas\ over $\C$. We should also emphasize that despite the presence of * everywhere, \cstas\ are \emph{not} required to be inner product spaces.
\end{minipage}
\begin{minipage}{1.5in}\begin{center} 
\personbox{Green!20}{\href{https://en.wikipedia.org/wiki/Charles_Earl_Rickart}{Charles Earl Rickart}\\(1913--2002)}{\includegraphics[width=1.3in]{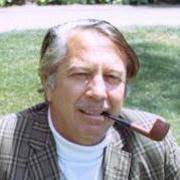}}\label{xrickart} 
\end{center}\end{minipage}

\vspace*{3mm}

\ind The main theorems we will prove in this chapter are the \hyperref[gelfandnaimark]{Gelfand-Naimark theorem} and the \hyperref[gns]{Gelfand-Naimark-Segal (GNS) theorem}. The former is for abelian \cstas\ and uses the example of $C_0(X)$, while the latter is for arbitrary \cstas\ and uses the example of $\Bs(\Hs)$.

\begin{warning}
Some sources refer to both of the aforementioned theorems as the Gelfand-Naimark theorem.
\end{warning}

\ind For further reading on \cstas, see \cite{muscat}, \cite{davidson}, \cite{strung}, \cite{conway2}, \cite{murphy} and \cite{bratteli}.

\subsection{$C^*$-algebras}
\begin{definition}
A \textbf{*-algebra} (or \textit{star algebra} or \textit{involutive algebra}) is an associative algebra $\As$, together with a map $^*:\As\to\As$, such that:
\begin{enumerate}
    \item If $a\in\As$, then $(a^*)^*=a$. \hfill (* is an involution)
    \item If $a,b\in\As$, then $(a+b)^*=a^*+b^*$. \hfill (* is conjugate linear)
    \item If $a\in\As$ and $\lam\in\F$, then $(\lam a)^*=\overline{\lam}a^*$.
    \item If $a,b\in\As$, then $(ab)^*=b^*a^*$. \hfill (* is anti-multiplicative)
\end{enumerate}
\end{definition}

\begin{proposition}
Suppose $\As$ is a unital *-algebra. Then $\1^*=\1$.
\end{proposition}
\begin{proof}
By definition, we have $\1a=a\1=a$ for all $a\in\As$. Applying * to all three expressions yields $a^*\1^*=\1^*a^*=a^*$. Thus $\1^*$ is also an identity of $\As$. Since the identity of $\As$ is unique, we have $\1^*=\1$.
\end{proof}

\begin{definition}
A \textbf{normed *-algebra} is a normed algebra with a map $^*:\As\to\As$ that satisfies the above properties. A \textbf{Banach *-algebra} is a normed *-algebra that is complete with respect to the induced metric.
\end{definition}

The norm provides \emph{analytic} structure to the algebra, while the * operation provides \emph{algebraic} structure. In general, there is no link between them. The necessary link comes now, with the following definition:

\begin{definition}
A \textbf{\csta} is a Banach *-algebra $\As$ such that $\norm{a^*a}=\norm{a}^2$ for all $a\in\As$.
\end{definition}

The above property (the \emph{\cst-axiom}) links the norm and the * operation together. It is sort of a compatibility requirement, and it has profound implications on both the algebraic and analytic structures of $\As$, as we will see later\footnote{Think, if you will, of two people at a party who are perfect for each other. Alone, they are sad and depressed, not having any fun. But once they start chatting, sparks fly, and the story of a lifetime begins. This is the beauty of the \cst-axiom.}.

\begin{remark}
It is enough to require that $\norm{a^*a}\ge\norm{a}^2$ for all $a\in\As$, since this would imply that $\norm{a}^2\le\norm{a^*a}\le\norm{a^*}\norm{a}$, so $\norm{a}\le\norm{a^*}$. Switching $a$ and $a^*$ yields $\norm{a^*}\le\norm{a}$, so $\norm{a^*}=\norm{a}$, and so $\norm{a^*a}\le\norm{a}^2$.
\end{remark}

\begin{examples}\leavevmode
\begin{enumerate}
    \item $\F$ is a \csta, with the absolute value norm and the complex conjugate as the adjoint.
    \item $M_n(\F)$ is a \csta, with the operator norm and the conjugate transpose as the adjoint.
    \item Suppose $X$ is a locally compact Hausdorff space. Then $C_0(X)$ is a \csta, with the supremum norm and the adjoint given by complex conjugation, i.e. for each $f\in C_0(X)$, the function $f^*\in C_0(X)$ is given by $f^*(x)=\overline{f(x)}$ for all $x\in X$.
    \item Suppose $(X,\As,\mu)$ is a measure space. Then $L^\infty(X,\As,\mu)$ is a \csta, with the essential supremum norm and the adjoint given by complex conjugation.
    \item Suppose $\Hs$ is a Hilbert space. Then $\Bs(\Hs)$ is a \csta, with the operator norm and the Hilbert space adjoint, i.e. for each $A\in\Bs(\Hs)$, the operator $A^*\in\Bs(\Hs)$ is given by $\ip{Ax,y}=\ip{x,A^*y}$ for all $x,y\in\Hs$.
\end{enumerate}
\end{examples}

{
\definecolor{col}{HTML}{1ACE90}
\begin{tcolorbox}[title=The Road to \cstas, fonttitle=\sffamily\bfseries\selectfont,colback=white,colframe=col!80,coltitle=black,halign title=center,top=1mm,bottom=1mm,left=1mm,right=1mm,before skip=10pt,after skip=10pt,enhanced]
\begin{center}
\begin{tikzpicture}[font=\sffamily]
\def\hs{3.75} 
\def\vs{1.5} 
\definecolor{col1}{HTML}{F41208} 
\definecolor{col2}{HTML}{0812F4} 
\node[align=center] (00) at (0,0) {Vector space};
\node[align=center] (01) at (\hs,0) {Algebra};
\node[align=center] (02) at (2*\hs,0) {*-algebra};
\node[align=center] (10) at (0,-\vs) {Normed vector space};
\node[align=center] (11) at (\hs,-\vs) {Normed algebra};
\node[align=center] (12) at (2*\hs,-\vs) {Normed *-algebra};
\node[align=center] (20) at (0,-2*\vs) {Banach space};
\node[align=center] (21) at (\hs,-2*\vs) {Banach algebra};
\node[align=center] (22) at (2*\hs,-2*\vs) {Banach *-algebra};
\node[align=center] (33) at (3*\hs,-3*\vs) {\csta};
\foreach \a/\b in {1/0,1/1,2/0,2/1}
{\draw[->,thick,col1] (\a\b)--(\a\the\numexpr\b+1\relax);
\draw[->,thick,col2] (\b\a)--(\the\numexpr\b+1\relax\a);}
\draw[->,thick,col1] (00)--(01) node[midway,above] {multiplication};
\draw[->,thick,col1] (01)--(02) node[midway,above] {*-operation};
\draw[->,thick,col2] (00)--(10) node[midway,left] {norm};
\draw[->,thick,col2] (10)--(20) node[midway,left] {completeness};
\draw[->,thick,col1] (22) edge[bend left=18] (33);
\draw[->,thick,col2] (22) edge[bend right=18] (33);
\node[col1!50!col2,xshift=3pt] at (2.5*\hs,-2.5*\vs) {\cst-axiom};
\draw[->,very thick,col1!50!black] (0.2*\hs,0.5*\vs)--(1.4*\hs,0.5*\vs) node[midway,above] {Increasing \textit{algebraic} structure};
\draw[->,very thick,col2!50!black] (-0.6*\hs,-0.2*\vs)--(-0.6*\hs,-1.4*\vs) node[midway,left,align=right,text width=2cm] {Increasing\\\textit{analytic}\\structure};
\end{tikzpicture}\\
\vspace{2mm}
\sffamily \cstas\ have \emph{far} more algebraic and analytic structure than all of the constructs we have introduced earlier.
\end{center}
\end{tcolorbox}
}

\ind To get the ball rolling, we will now prove a few results about \cstas\ that will routinely come up later in this chapter. It is helpful to observe how the \cst-axiom shows up in the proofs. Once we get to the later sections of this chapter, we will usually apply results like these instead of using the \cst-axiom explicitly.

\begin{proposition}\label{norma*=norma}
Suppose $\As$ is a \csta. Then for all $a\in\As$, we have $\norm{a^*}=\norm{a}$.\\
In other words, the * operation preserves norms.
\end{proposition}
\begin{proof}
If $a=0$, this holds trivially as both sides are zero. Suppose $a\ne0$. Then we have $\norm{a}^2=\norm{a^*a}\le\norm{a^*}\norm{a}$. Dividing both sides by $\norm{a}$, we get $\norm{a}\le\norm{a^*}$. Switching $a$ and $a^*$ yields $\norm{a^*}\le\norm{a}$, and so $\norm{a^*}=\norm{a}$.
\end{proof}

\begin{warning}
The * operation is NOT an isometry if $\F=\C$, as it is conjugate linear rather than linear.
\end{warning}

\begin{corollary}
Suppose $\As$ is a \csta. Then the * operation on $\As$ is a homeomorphism.
\end{corollary}
\begin{proof}
Suppose $\ep>0$ and set $\delta=\ep$. By the previous proposition, for any $a,b\in\As$, we have $\norm{a^*-b^*}=\norm{(a-b)^*}=\norm{a-b}$. Thus if $\norm{a-b}<\delta$, we have $\norm{a^*-b^*}<\delta=\ep$, and so * is continuous. Since * is an involution, it is its own inverse, and so it is a homeomorphism.
\end{proof}

The next proposition shows that the norm in a \csta\ behaves somewhat like the operator norm:

\begin{proposition}\label{cstaopnorm}
Suppose $\As$ is a \csta. Then for all $a\in\As$, we have:
\begin{equation*}
    \norm{a}=\sup_{x\in\As,\norm{x}\le1} \norm{ax}=\sup_{x\in\As,\norm{x}\le1} \norm{xa}
\end{equation*}
\end{proposition}
\begin{proof}
Define $M=\sup_{x\in\As,\norm{x}\le1} \norm{ax}$. Since $\norm{ax}\le\norm{a}\norm{x}=\norm{a}$ for all $x\in\As$, $\norm{x}\le1$, we have $M\le\norm{a}$. Setting $x=\frac{a^*}{\norm{a}}$ (which has norm $1$ by \Cref{norma*=norma}) yields $\norm{ax}=\norm{\frac{aa^*}{\norm{a}}}=\frac{\norm{a}^2}{\norm{a}}=\norm{a}$. Thus $M=\norm{a}$. A similar proof shows that $\norm{a}=\sup_{x\in\As,\norm{x}\le1} \norm{xa}$.
\end{proof}
\begin{remark}
An equivalent way to write this is $\Disp\norm{a}=\sup_{x\in\As,x\ne0} \frac{\norm{ax}}{\norm{x}}=\sup_{x\in\As,x\ne0} \frac{\norm{xa}}{\norm{x}}$ (provided $\As\ne\{0\}$).
\end{remark}

\xnewpage

In \Cref{unitalnormedalgebra}, we remarked that the condition $\norm{\1}=1$ is not vacuous and thus a necessary part of the definition. However, if $\As$ is a \csta, this is indeed vacuous, as the \cst-axiom prevents any other possibility.

\begin{corollary}
Suppose $\As$ is a \csta\ with an identity. Then $\norm{\1}=1$.
\end{corollary}
This follows directly from the previous proposition by setting $a=\1$.

\begin{theorem}
Suppose $\As$ is a \csta. Then $\As$ is isometrically isomorphic to a subalgebra of $\Bs(\As)$.
\end{theorem}
\begin{proof}
For each $a\in\As$, define $f_a:\As\to\As$ by $f_a(x)=ax$. Clearly $f_a$ is linear. By \Cref{cstaopnorm}, we have $\norm{f_a}=\sup_{x\in\As,\norm{x}\le1} \norm{ax}=\norm{a}$. Thus $f_a$ is bounded and $\norm{f_a}=\norm{a}$.\\

Now define $\rho:\As\to\Bs(\As)$ by $\rho(a)=f_a$. Clearly $\rho$ is linear, and as we have just shown, $\norm{\rho(a)}=\norm{f_a}=\norm{a}$, so $\rho$ is an isometry. It is also a homomorphism as $\rho(ab)=f_{ab}=f_a\circ f_b=\rho(a)\circ\rho(b)$. Thus $\rho:\As\to\Bs(\As)$ is an isometric homomorphism, and so it is an isometric isomorphism from $\As$ to $\im(\rho)\sub\Bs(\As)$.
\end{proof}
\begin{remark}
This theorem holds for all Banach algebras $\As$ in which $\norm{a}=\sup_{x\in\As,\norm{x}\le1} \norm{ax}$ for all $a\in\As$. This is reminiscent of Cayley's theorem, which states that every group $G$ is isomorphic (as a group) to a subgroup of the corresponding symmetric group $\operatorname{Sym}(G)$.
\end{remark}

\begin{definition}
Suppose $\As$ is a \csta\ and $a\in\As$.
\begin{itemize}
    \item $a$ is \textbf{Hermitian} (or \textit{self-adjoint}) if $a^*=a$.
    \item $a$ is \textbf{unitary} if $aa^*=a^*a=\1$ (this only applies if $\As$ is unital).
    \item $a$ is \textbf{normal} if $aa^*=a^*a$.
\end{itemize}
\end{definition}
\begin{remark}
It follows directly that every Hermitian or unitary element is normal.
\end{remark}

Just like for matrices, every element $a\in\As$ can be expressed uniquely as $a=b+ic$, where $b,c\in\As$ are Hermitian. Explicitly, these are given by $b=\frac{a+a^*}{2}$ and $c=\frac{a-a^*}{2i}$.

\begin{theorem}\label{sprad=norm}
Suppose $\As$ is a \csta\ and $a\in\As$ is normal. Then $r(a)=\norm{a}$.
\end{theorem}
\begin{proof}
Since $a$ is normal, we have:
\begin{equation*}
    (a^*a)^*(a^*a)=(a^*(a^*)^*)(a^*a)=(a^*a)^2=a^*aa^*a=a^*a^*aa=(a^*)^2a^2=(a^2)^*a^2
\end{equation*}
Thus $\norm{a^2}^2=\norm{(a^2)^*a^2}=\norm{(a^*a)^*(a^*a)}=\norm{a^*a}^2=\norm{a}^4$. Taking the square root of both sides yields $\norm{a^2}=\norm{a}^2$. It follows by induction that $\norm{a^{2^m}}=\norm{a}^{2^m}$ for all $m\in\N$. By the spectral radius formula (\Cref{specradius}), we have:
\begin{equation*}
    r(a)=\limn\norm{a^n}^{1/n}=\limm\norm\big{a^{2^m}}^{2^{-m}}=\limm\big(\norm{a}^{2^m}\big)^{2^{-m}}=\norm{a} \qedhere
\end{equation*}
\end{proof}
\begin{remark}
An important special case of this theorem is when $\As$ is abelian. In this case, every $a\in\As$ is normal, and so $r(a)=\norm{a}$ for \emph{all} $a\in\As$.
\end{remark}

\begin{corollary}\label{onlynorm}
For any *-algebra $\As$, there is at most one norm on $\As$ that makes it a \csta. This norm is given by $\norm{\cdot}:\As\to[0,\infty)$, $\norm{a}=\sqrt{r(a^*a)}$.
\end{corollary}
\begin{proof}
Suppose $\norm{\cdot}$ is a norm on $\As$ that makes it a \csta. Since $a^*a$ is Hermitian (and thus normal), by the previous theorem, we have $r(a^*a)=\norm{a^*a}=\norm{a}^2$. Thus $\norm{a}=\sqrt{r(a^*a)}$. Since the spectral radius does not depend on the norm, this formula uniquely determines the norm $\norm{\cdot}$.
\end{proof}

\begin{example}
Suppose $\As=M_n(\C)$ with the usual operations and the * operation given by the conjugate transpose. Then the only norm that makes it a \csta\ is the spectral norm. For any matrix $A\in M_n(\C)$, its singular values are the square roots of the eigenvalues of $A^*A$. Thus $\norm{A}=\sqrt{r(A^*A)}$ is also the largest singular value of $A$.
\end{example}

\ind The last corollary shows that the norm of a \csta\ is uniquely determined. On the other hand, the * operation is also uniquely determined: For any Banach algebra $\As$, there is at most one * operation on $\As$ that makes it a \csta. See \cite[Theorem 7.9, Page 59]{wilde} for a proof.

\xnewpage

\begin{definition}
Suppose $\As$ is a *-algebra. A \textbf{*-subalgebra} of $\As$ is a subalgebra $\Bs$ of $\As$ such that $a^*\in\Bs$ for all $a\in\Bs$.\\
In other words, a *-subalgebra of $\As$ is a subalgebra of $\As$ that is also closed under the * operation.
\end{definition}
\begin{remark}
Some sources simply call this a \emph{subalgebra}. To avoid confusion, we will not do this here.
\end{remark}

It follows directly from the definitions that a subset of a \csta\ $\As$ is a \csta\ if and only if it is a closed *-subalgebra of $\As$.

\begin{definition}
Suppose $\As$ and $\Bs$ are \cstas. A \textbf{*-homomorphism} from $\As$ to $\Bs$ is an algebra homomorphism $\phi:\As\to\Bs$ such that $\phi(a^*)=\phi(a)^*$ for all $a\in\As$.\\
A \textbf{*-isomorphism} is a *-homomorphism that is also a bijection.\\
In other words, a *-homomorphism is an algebra homomorphism that preserves the structure of the adjoint operation.
\end{definition}

\begin{proposition}
Suppose $\As$ and $\Bs$ are \cstas\ and $\phi:\As\to\Bs$ is a *-homomorphism. Then $\im(\phi)$ is a *-subalgebra of $\Bs$.
\end{proposition}
\begin{proof}
Since $\phi$ is a homomorphism, by \Cref{imskeri}, $\im(\phi)$ is a subalgebra of $\Bs$. Suppose $b\in\im(\phi)$. Then we have $b=\phi(a)$ for some $a\in\phi$. Thus $\phi(a^*)=\phi(a)^*=b^*\in\im(\phi)$, and so $\im(\phi)$ is a *-subalgebra of $\Bs$.
\end{proof}
\begin{remark}
It can be shown that $\im(\phi)$ is also closed in $\Bs$ (and thus also a \csta), see \cite[Theorem 1.3.2, Page 13]{arveson1}.
\end{remark}

\begin{lemma}\label{specsub1}
Suppose $\As$ and $\Bs$ are unital algebras and $\phi:\As\to\Bs$ is a unital homomorphism. Then for all $a\in\As$, we have $\sigma(\phi(a))\sub\sigma(a)$.\\
In other words, unital homomorphisms cannot add new elements to the spectrum.
\end{lemma}
\begin{proof}
Suppose $\lam\notin\sigma(a)$. Then $a-\lam\1_\As\in\As\inv$. By \Cref{unithominv}, we have $\phi(a-\lam\1_\As)=\phi(a)-\lam\1_\Bs\in\Bs\inv$, and so $\lam\notin\sigma(\phi(a))$. Thus $\sigma(\phi(a))\sub\sigma(a)$.
\end{proof}

\begin{theorem}\label{*homobdd}
Suppose $\As$ and $\Bs$ are \cstas\ and $\phi:\As\to\Bs$ is a *-homomorphism. Then $\phi$ is bounded and $\norm{\phi}\le1$.
\end{theorem}
\begin{proof}
Assume without loss of generality that $\As$ and $\Bs$ are unital and $\phi$ is a unital *-homomorphism (otherwise, replace them with their unitizations $\As_1$ and $\Bs_1$ defined in \Cref{embedding} and extend $\phi$ to $\As_1$ by defining $\phi(\1_{\As_1})=\1_{\Bs_1}$). Suppose $a\in\As$ is Hermitian. Then $\phi(a)=\phi(a^*)=\phi(a)^*$, so $\phi(a)$ is also Hermitian. Thus $r(a)=\norm{a}$ and $r(\phi(a))=\norm{\phi(a)}$. By \Cref{specsub1}, we have $\sigma(\phi(a))\sub\sigma(a)$, so $r(\phi(a))\le r(a)$. Thus $\norm{\phi(a)}\le\norm{a}$.\\

Now suppose $a\in\As$. Since $a^*a$ is Hermitian, so is $\phi(a^*a)$, and so $\norm{\phi(a)}^2=\norm{\phi(a)^*\phi(a)}=\norm{\phi(a^*a)}\le\norm{a^*a}=\norm{a}^2$. Thus $\norm{\phi(a)}\le\norm{a}$, and so $\phi$ is bounded and $\norm{\phi}\le1$.
\end{proof}

\ind The last theorem should surprise (and hopefully delight) you. A *-homomorphism is a purely \emph{algebraic} notion. Why on Earth should *-homomorphisms automatically be continuous? This is yet another example of the power of the \cst-axiom. Once you demand enough structure for a map between \cstas\ (in this case, that the map is a *-homomorphism), you get even \emph{more} structure (in this case, continuity) for free!\\

\ind In fact, $\phi$ does not even need to be a *-homomorphism. As long as it preserves addition, vector multiplication and the adjoint (but not necessarily scalar multiplication), it is automatically continuous, and $\norm{\phi(a)}\le\norm{a}$ for all $a\in\As$. See \cite{tomforde} for a proof.

\begin{proposition}\label{cstaspec}
Suppose $\As$ is an abelian \csta, $a\in\As$ and $\chi\in\Sigma(\As)$. Then:
\begin{enumerate}
    \item If $a$ is Hermitian, then $\chi(a)\in\R$.
    \item $\chi(a^*)=\overline{\chi(a)}$
    \item $\chi(a^*a)\ge0$
    \item If $\As$ is unital and $a$ is unitary, then $\abs{\chi(a)}=1$.
\end{enumerate}
\end{proposition}
\begin{proof}
We will assume $\F=\C$ to prove (1) and (2), since they hold trivially if $\F=\R$.
\begin{enumerate}
    \item Suppose $c\in\As$ such that $\chi(c)=1$, and $t\in\R$. Then we have:
    \begin{equation*}
        \abs{\chi(a+itc)}^2\le\norm{a+itc}^2=\norm{(a+itc)(a-itc)}=\norm{a^2+t^2c^2}\le\norm{a}^2+\abs{t}^2\norm{c}^2=\norm{a}^2+\abs{t}^2
    \end{equation*}
    Suppose $\chi(a)=\alpha+i\beta$, where $\alpha,\beta\in\R$. Then:
    \begin{equation*}
        \norm{a}^2+t^2\ge\abs{\chi(a+itc)}^2=\abs{\alpha+i(\beta+t)c}^2=\alpha^2+\beta^2+2\beta t+t^2
    \end{equation*}
    Thus $\norm{a}^2\ge\alpha^2+\beta^2+2\beta t$. This cannot hold for arbitrarily large $t\in\R$ unless $\beta=0$. Thus $\beta=0$, and so $\chi(a)\in\R$.
    \item Suppose $a=b+ic$, where $b$ and $c$ are Hermitian. By (1), we have $\chi(b),\chi(c)\in\R$. This yields:
    \begin{equation*}
        \chi(a^*)=\chi(b-ic)=\chi(b)-i\chi(c)=\overline{\chi(b)+i\chi(c)}=\overline{\chi(b+ic)}=\overline{\chi(a)}
    \end{equation*}
    \item By (2), we have $\chi(a^*a)=\chi(a^*)\chi(a)=\overline{\chi(a)}\chi(a)=\abs{\chi(a)}^2\ge0$.
    \item By \Cref{char1inv} and \Cref{charbd}, we have $\norm{\chi}=1$ and $\chi(\1)=1$. By (2), we have $\abs{\chi(a)}^2=\chi(a^*)\chi(a)=\chi(a^*a)=\chi(\1)=1$, so $\abs{\chi(a)}=1$. \qedhere
\end{enumerate}
\end{proof}
\begin{remark}
(2) shows that every homomorphism $\chi:\As\to\C$ is a *-homomorphism.
\end{remark}

\begin{lemma}\label{specsub2}
Suppose $\As$ is a unital algebra and $\Bs$ is a unital subalgebra of $\As$. Then for all $a\in\Bs$, we have $\sigma_\As(a)\sub\sigma_\Bs(a)$.\\
In other words, the spectrum of $a$ in the larger algebra is contained in the spectrum of $a$ in the smaller algebra.
\end{lemma}
This follows directly from \Cref{specsub1} by setting $\phi:\Bs\to\As$, $\phi(a)=a$ (the inclusion map, which is a unital homomorphism).

\begin{theorem}\label{hermunitspec}
Suppose $\As$ is a \csta\ and $a\in\As$.
\begin{enumerate}
    \item If $a$ is Hermitian, then $\sigma(a)\sub\R$.
    \item If $a=c^*c$ for some $c\in\As$, then $\sigma(a)\sub[0,\infty)$.
    \item If $a$ is unitary, then $\sigma(a)\sub\partial D(0,1)$ (the unit circle $\{z\in\C\mid\abs{z}=1\}$).
\end{enumerate}
\end{theorem}
\begin{proof}
Assume without loss of generality that $\As$ is unital (otherwise, replace $\As$ with its unitization $\As_1$ defined in \Cref{embedding}). We will prove these results by first passing to a subalgebra of $\As$ that is abelian, which will allow us to use characters to help us.\\

Define $\Ss$ as the closed unital *-subalgebra of $\As$ generated by $a$ (or equivalently, the closure of the smallest subalgebra of $\As$ containing $\1$, $a$ and $a^*$). Then every element of $\As$ can be expressed as a (finite or infinite) power series in $a$ and $a^*$. Since $a$ is normal, all such elements commute, and so $\Ss$ is abelian.\\

We now consider the set $\{\chi(a)\mid\chi\in\Sigma(\Ss)\}$, i.e. the set of values obtained by evaluating all characters of $\Ss$ at $a$. By \Cref{gelfspec=spec}, we have $\{\chi(a)\mid\chi\in\Sigma(\Ss)\}=\sigma_\Ss(a)$, and by \Cref{specsub2}, we have $\sigma_\As(a)\sub\sigma_\Ss(a)$. Thus $\sigma_\As(a)\sub\{\chi(a)\mid\chi\in\Sigma(\Ss)\}$.\\

Finally, we use \Cref{cstaspec} to conclude all three results at once:
\begin{enumerate}
    \item If $a$ is Hermitian, then $\chi(a)\in\R$ for all $\chi\in\Sigma(\Ss)$, and so $\sigma_\As(a)\sub\R$.
    \item If $a=c^*c$ for some $c\in\As$, then $\chi(a)=\chi(c^*c)\ge0$ for all $\chi\in\Sigma(\Ss)$, and so $\sigma_\As(a)\sub[0,\infty)$.
    \item If $a$ is unitary, then $\abs{\chi(a)}=1$ for all $\chi\in\Sigma(\Ss)$, and so $\sigma_\As(a)\sub\partial D(0,1)$. \qedhere
\end{enumerate}
\end{proof}

\begin{theorem}\label{hermunit}
Suppose $\As$ is a unital Banach *-algebra and $a\in\As$. Then:
\begin{enumerate}
    \item $\exp(a^*)=\exp(a)^*$
    \item If $a$ is Hermitian, then $\exp(ia)$ is unitary.
\end{enumerate}
\end{theorem}
\begin{proof}
\begin{enumerate}
    \item $\Disp\exp(a^*)=\sumn[0]\frac{1}{n!}(a^*)^n=\sumn[0]\frac{1}{n!}(a^n)^*=\left(\sumn[0]\frac{1}{n!}a^n\right)^*=\exp(a)^*$
    \item Replacing $a$ with $ia$ in (1) yields $\exp((ia)^*)=\exp(ia)^*$, i.e. $\exp(ia)^*=\exp(-ia^*)$. Since $a$ is Hermitian, we have $\exp(-ia^*)=\exp(-ia)=(\exp(ia))^{-1}$. Thus $\exp(ia)^*=\exp(ia)^{-1}$, and so $\exp(ia)$ is unitary.\qedhere
\end{enumerate}
\end{proof}

\begin{example}
$A=\mmat{0}{\pi}{\pi}{0}$ is Hermitian, and $\exp(iA)=\exp\mmat{0}{i\pi}{i\pi}{0}=\mmat{-1}{0}{0}{-1}$ is unitary.
\end{example}

\xnewpage

\ind We are almost ready to prove the \hyperref[gelfandnaimark]{Gelfand-Naimark theorem}. The last ingredient we will need is the following result, which is a slight generalization of the Stone-Weierstrass theorem.

\begin{theorem}\label{swlc}
Suppose $X$ is a locally compact Hausdorff space and $\Rs$ is a *-subalgebra of $C_0(X)$ such that:
\begin{itemize}
    \item $\Rs$ \emph{separates points}: For any $x,y\in X$, $x\ne y$, there exists $f\in\Rs$ such that $f(x)\ne f(y)$.
    \item $\Rs$ \emph{vanishes nowhere}: For any $x\in X$, there exists $f\in\Rs$ such that $f(x)\ne0$.
\end{itemize}
Then $\Rs$ is dense in $C_0(X)$.
\end{theorem}
See \cite{stonewlc} for a proof.

\begin{theorem}
Suppose $\As$ is an abelian \csta. Then the Gelfand transform $\gamma$ of $\As$ is an isometric *-isomorphism from $\As$ onto $C_0(\Sigma(\As))$.
\end{theorem}
\begin{proof}
By the Gelfand representation theorem (\Cref{gelfandrep}), the Gelfand transform $\gamma$ of $\As$ is a continuous *-homomorphism from $\As$ to $C_0(\Sigma(\As))$, and for all $a\in\As$, we have $\norm{\ahat}_\infty=r(a)$. Since $\As$ is abelian, every $a\in\As$ is normal, so by \Cref{sprad=norm}, we have $r(a)=\norm{a}$. Thus $\norm{\ahat}_\infty=r(a)=\norm{a}$ for all $a\in\As$, and so $\gamma$ is an isometry.\\

For all $a\in\As$ and all $\chi\in\Sigma(\As)$, we have $\widehat{a^*}(\chi)=\chi(a^*)=\overline{\chi(a)}=\overline{\ahat(\chi)}=\ahat^*(\chi)$. Thus $\widehat{a^*}=\ahat^*$, i.e. $\gamma(a^*)=\gamma(a)^*$, and so $\gamma$ is a *-homomorphism.\\

We now show that $\gamma$ is onto. Define $\Rs=\im(\gamma)$. Since $\gamma$ is a *-homomorphism, $\Rs$ is a *-subalgebra of $C_0(\Sigma(\As))$.
\begin{itemize}
    \item Suppose $\chi_1,\chi_2\in\sigma(\As)$ and $\chi_1\ne\chi_2$. Then there exists $a\in\As$ such that $\chi_1(a)\ne\chi_2(a)$, i.e. $\ahat(\chi_1)\ne\ahat(\chi_2)$. Thus $\Rs$ separates points.
    \item Now suppose $\chi\in\Sigma(\As)$. Then $\chi\ne0$ by assumption, so there exists $b\in\As$ such that $\chi(b)\ne0$, i.e. $\bhat(\chi)\ne0$. Thus $\Rs$ vanishes nowhere.
\end{itemize}
By \Cref{swlc}, $\Rs$ is dense in $C_0(\Sigma(\As))$. Finally, by \Cref{isomimclosed}, $\Rs$ is closed in $C_0(\Sigma(\As))$, so it must be equal to $C_0(\Sigma(\As))$. Thus $\gamma$ is onto, and so it is an isometric *-isomorphism.
\end{proof}

\ind In the above proof, we did not need to show directly that $\gamma$ is injective, this follows automatically (using \Cref{isomcontinj}) from the fact that it is an isometry. Consequently, the Gelfand transform of any abelian \csta\ is injective. Furthermore, if $\As$ is unital and complex, we can deduce that the intersection of all maximal ideals of $\As$ is the zero ideal\footnote{In the language of ring theory, we would say every unital abelian complex \csta\ is \emph{Jacobson semisimple}.} $\{0\}$ (since it is $\ker(\gamma)$ and $\gamma$ is injective).

\begin{theorem}[Gelfand-Naimark Theorem]\label{gelfandnaimark}
Every abelian \csta\ is isometrically *-isomorphic to $C_0(X)$, where $X$ is some locally compact Hausdorff space.
\end{theorem}
\begin{proof}
Suppose $\As$ is an abelian \csta. Define $X=\Sigma(\As)$ and $\gamma:\As\to C_0(X)$ as the Gelfand transform of $\As$. By the previous theorem, $\gamma$ is an isometric *-isomorphism.
\end{proof}

\begin{example}
Suppose $X$ is a locally compact Hausdorff space. We already know that $C_0(X)$ is an abelian \csta. Let's see what happens when we construct the Gelfand transform of $C_0(X)$. The characters of $C_0(X)$ are the evaluation maps $\chi_z:C_0(X)\to\F$, $\chi_z(f)=f(z)$ (where $z\in X$). Each character $\chi_z$ naturally corresponds to exactly one $z\in X$. Also, for each $f\in C_0(X)$, the Gelfand transform of $f$ is given by $\fh:\Sigma(C_0(X))\to\F$, $\fh(\chi_z)=\chi_z(f)=f(z)$. Thus $C_0(\Sigma(C_0(X)))\cong C_0(X)$. In other words, if we perform this construction on $C_0(X)$, we get back $C_0(X)$.
\end{example}

\subsection{Positive Elements and Positive Linear Functionals}
\ind We now introduce the notions of \emph{positive elements} and \emph{positive linear functionals}. Both of these are generalizations of non-negative real numbers, positive semidefinite matrices (Hermitian matrices whose eigenvalues are all non-negative), and non-negative functions (real-valued functions that only assume non-negative values). It may be more appropriate to call these \textit{non-negative} elements and linear functionals, since strictly speaking, the aforementioned values are allowed to be zero. Nonetheless, we call them ``positive'' as it is simpler.

\begin{definition}
Suppose $\As$ is a \csta. An element $a\in\As$ is \textbf{positive} if $a^*=a$ and $\sigma(a)\sub[0,\infty)$.\\
In other words, $a$ is positive if it is Hermitian and its spectrum consists only of non-negative real numbers.
\end{definition}

\begin{examples}\leavevmode
\begin{enumerate}
    \item A matrix $A\in M_n(\F)$ is positive if and only if it is positive semidefinite.
    \item Suppose $X$ is a locally compact Hausdorff space. An element $f\in C_0(X)$ is positive if and only if $f(x)\ge0$ for all $x\in X$.
    \item Suppose $\Hs$ is a Hilbert space. An element $T\in\Bs(\Hs)$ is positive if and only if $\ip{Tx,x}\ge0$ for all $x\in\Hs$.
    \item Suppose $(X,\As,\mu)$ is a measure space. An element $f\in L^\infty(X,\As,\mu)$ is positive if and only if $f(x)\ge0$ a.e. on $X$.
\end{enumerate}
\end{examples}

\begin{proposition}
Suppose $\As$ is a \csta\ and $a\in\As$. Then the following are equivalent:
\begin{enumerate}
    \item $a$ is positive.
    \item $a=b^2$ for some Hermitian $b\in\As$.
    \item $a=c^*c$ for some $c\in\As$.
\end{enumerate}
\end{proposition}
\begin{proof}
\begin{enumerate}[leftmargin=.55in]
    \item[($1\Rightarrow2$)] Note that $\sigma(a)\sub[0,\norm{a}]$ and the square root function $\Phi:[0,\norm{a}]\to\R$, $\Phi(x)=\sqrt{x}$ is continuous on $[0,\norm{a}]$. Thus $b=\Phi(a)$ is well-defined. Since $\Phi$ is real-valued, $b$ is Hermitian, and $b^2=\Phi(a)^2=a$.
    \item[($2\Rightarrow3$)] Since $b$ is Hermitian, we have $a=b^2=b^*b$, so we can set $c=b$.
    \item[($3\Rightarrow1$)] Clearly $a$ is Hermitian, as $a^*=(c^*c)^*=c^*(c^*)^*=c^*c=a$. By \Cref{hermunitspec}, we have $\sigma(a)\sub[0,\infty)$, and so $a$ is positive. \qedhere
\end{enumerate}
\end{proof}

\begin{example}
The matrix $A=\mmat{25}{40}{40}{65}$ is positive (it is Hermitian and its spectrum is $\{0.28,89.72\}\subset[0,\infty)$). It can be expressed as $A=B^2$, where $B=\mmat{3}{4}{4}{7}$ is Hermitian (and thus $A=B^*B$ as well). This is not the only way to express $A$ as $C^*C$, e.g. we could also use $C=\mmat{5}{8}{0}{1}$.
\end{example}
\begin{example}
The function $f:[-1,1]\to\C$, $f(x)=x^4$ is positive (since $f(x)\ge0$ for all $x\in[-1,1]$). It can be expressed as $f=g^2$, where $g(x)=x^2$ (and $g$ is Hermitian as it is real-valued). We could also express $f$ as $h^*h$, where $h=ix\abs{x}$.
\end{example}

\begin{definition}
Suppose $\As$ is a \csta. A linear functional $f:\As\to\F$ is \textbf{positive} if $f(a^*a)\ge0$ for all $a\in\As$.\\
A \textbf{state} on $\As$ is a positive linear functional $f:\As\to\F$ such that $\norm{f}=1$.
\end{definition}
\begin{remark}
A state is simply a normalized positive linear functional, i.e. if $f\ne0$ is a positive linear functional on $\As$, then $\frac{f}{\norm{f}}$ is a state on $\As$. You might be wondering if this requires us to assume $f$ is bounded, but as we will show in \Cref{+lf}, this follows automatically from being positive.
\end{remark}

\begin{example}
Suppose $f:M_2(\C)\to\C$, $f\mmat{a}{b}{c}{d}=a$. This is clearly linear, and it is positive as every matrix of the form $A^*A$ is positive semidefinite, so its upper-left entry must be non-negative.
\end{example}

We first present two important examples of states that we will need later:

\begin{proposition}\label{ipx=state}
Suppose $\Hs$ is a Hilbert space and $x\in\Hs$, $\norm{x}=1$. Define $f:\Bs(\Hs)\to\F$ by $f(A)=\ip{Ax,x}$. Then $f$ is a state on $\Bs(\Hs)$.
\end{proposition}
\begin{proof}
For any $A,B\in\Bs(\Hs)$ and any $\lam\in\F$, we have:
\begin{equation*}
    f(\lam A+B)=\ip{(\lam A+B)x,x}=\ip{\lam Ax+Bx,x}=\lam\ip{Ax,x}+\ip{Bx,x}=\lam f(A)+f(B)
\end{equation*}
Thus $f$ is linear. Also, for any $A\in\Bs(\Hs)$, we have $f(A^*A)=\ip{A^*A,x}=\ip{Ax,Ax}=\norm{Ax}^2\ge0$. Thus $f$ is positive. Moreover, we have:
\begin{equation*}
    \abs{f(A)}=\abs{\ip{Ax,x}}\le\norm{Ax}\norm{x}=\norm{A}\norm{x}^2=\norm{A}
\end{equation*}
Thus $f$ is bounded and $\norm{f}\le1$. Finally, if $A=I$ (the identity operator on $\Hs$), we have:
\begin{equation*}
    f(I)=\ip{Ix,x}=\ip{x,x}=\norm{x}^2=1
\end{equation*}
Thus $\norm{f}=1$, and so $f$ is a state on $\Bs(\Hs)$.
\end{proof}
\begin{remark}
The vector $x$ is not unique, for example we can use $\alpha x$ for any $\alpha\in\F$, $\abs{\alpha}=1$ (or for general $\alpha\in\F$, the corresponding $f$ will be a positive linear functional with norm $\norm{f}=\abs{\alpha}^2$).
\end{remark}

\xnewpage

\begin{proposition}\label{char=state}
Suppose $\As$ is an abelian \csta. Then every character of $\As$ is a state on $\As$.
\end{proposition}
\begin{proof}
Suppose $\chi\in\Sigma(\As)$. By definition, $\chi$ is a linear functional on $\As$. By \Cref{cstaspec}, we have $\chi(a^*a)\ge0$ for all $a\in\As$, so $\chi$ is positive. By \Cref{charbd}, $\chi$ is bounded and $\norm{\chi}=1$. Thus $\chi$ is a state on $\As$.
\end{proof}

To prove the next proposition, we will need the following modification of the Hahn-Banach theorem:

\begin{theorem}[Hahn-Banach Theorem for Positive Linear Functionals on \cstas]\label{hbtcsta}
Suppose $\As$ is a \csta, $\Bs$ is a closed subalgebra of $\As$ and $f:\Bs\to\F$ is a positive linear functional. Then $f$ has an extension to a positive linear functional on $\As$ with the same norm, i.e. there is a positive linear functional $F:\As\to\F$ such that $F(a)=f(a)$ for all $a\in\Bs$ and $\norm{F}=\norm{f}$.
\end{theorem}
This is a special case of a much more general theorem. See \cite{muhamadiev} for a proof and related results.

\begin{proposition}\label{stateexists}
Suppose $\As$ is a \csta\ and $a\in\As$ is normal. If $\lam\in\sigma(a)$, then there is a state $f:\As\to\F$ such that $f(a)=\lam$.
\end{proposition}
\begin{proof}
Define $\Ss$ as in \Cref{hermunitspec} (again, $\Ss$ is abelian as $a$ is normal). By \Cref{gelfspec=spec}, we have $\sigma_\Ss(a)=\{\chi(a)\mid\chi\in\Sigma(\Ss)\}$, and by \Cref{specsub2}, we have $\sigma_\As(a)\sub\sigma_\Ss(a)$. Thus if $\lam\in\sigma_\As(a)$, there is a character $\chi\in\Sigma(\Ss)$ such that $\chi(a)=\lam$. By \Cref{char=state}, $\chi$ is a state on $\Ss$, so by \Cref{hbtcsta}, $\chi$ can be extended to a positive linear functional $F$ on $\As$ with $\norm{F}=\norm{\chi}=1$.
\end{proof}

\begin{corollary}\label{stateexists+}
Suppose $\As$ is a \csta\ and $a\in\As$ is positive. Then there is a state $f:\As\to\F$ such that $f(a)=\norm{a}$.
\end{corollary}
\begin{proof}
Since $a$ is positive, it is Hermitian (and thus normal), so by \Cref{sprad=norm}, we have $r(a)=\norm{a}$. By \Cref{ftba}, $\sigma(a)$ is compact, so there exists $\lam\in\sigma(a)$ such that $\abs{\lam}=r(a)$. Since $\sigma(a)\sub[0,\infty)$, this implies that $\lam=r(a)$, and so $\lam\in\sigma(a)$. By the previous proposition, there is a state $f:\As\to\F$ such that $f(a)=\lam=\norm{a}$.
\end{proof}

\ind We will now prove a version of the Cauchy-Schwarz inequality for positive linear functionals. Usually in functional analysis, one would first prove that a map is an inner product and then deduce that it satisfies the Cauchy-Schwarz inequality. Here, we are going the other way: We first prove the Cauchy-Schwarz inequality holds for some ``inner product'' (which we will define in \Cref{gnsconst}), and later use this to prove it is indeed an inner product.

\begin{theorem}[Cauchy-Schwarz Inequality for Positive Linear Functionals]\label{cauchycsta}
Suppose $\As$ is a \csta\ and $f:\As\to\C$ is a positive linear functional. Then for all $a,b\in\As$, we have:
\begin{equation*}
    \abs{f(b^*a)}^2\le f(a^*a)f(b^*b)
\end{equation*}
\end{theorem}
\begin{proof}
Since $f$ is positive, for all $a,b\in\As$ and all $\lam,\mu\in\F$, we have $f((\lam a+\mu b)^*(\lam a+\mu b))\ge0$. Expanding the left side, we get:
\begin{align*}
    f\left((\lam a+\mu b)^*(\lam a+\mu b)\right)
    &=f\left((\overline{\lam}a^*+\overline{\mu}b^*)(\lam a+\mu b)\right)
    =f\left(\overline{\lam}\lam a^*a+\overline{\lam}\mu a^*b+\overline{\mu}\lam b^*a+\overline{\mu}\mu b^*b\right) \\
    &=\overline{\lam}\lam f(a^*a)+\overline{\lam}\mu f(a^*b)+\overline{\mu}\lam f(b^*a)+\overline{\mu}\mu f(b^*b)
\end{align*}
We now set $\lam=\overline{f(b^*a)}$ and $\mu=-f(a^*a)$ (note that \textcolor{Green!80!black}{$\mu\in\R$}) to get:
\begin{align*}
    0&\le f(b^*a)\overline{f(b^*a)}f(a^*a)-f(b^*a)f(a^*a)\textcolor{red!70!black}{f(a^*b)}-\textcolor{Green!80!black}{\overline{f(a^*a)}}\overline{f(b^*a)}f(b^*a)+\textcolor{Green!80!black}{\overline{f(a^*a)}}f(a^*a)f(b^*b) \\
    &=f(b^*a)\overline{f(b^*a)}f(a^*a)-f(b^*a)f(a^*a)\textcolor{red!70!black}{\overline{f(b^*a)}}-\textcolor{Green!80!black}{f(a^*a)}\overline{f(b^*a)}f(b^*a)+\textcolor{Green!80!black}{f(a^*a)}f(a^*a)f(b^*b) \\
    &=f(a^*a)\left(f(b^*a)\overline{f(b^*a)}-f(b^*a)\overline{f(b^*a)}-\overline{f(b^*a)}f(b^*a)+f(a^*a)f(b^*b)\right) \\
    &=f(a^*a)\left(\abs{f(b^*a)}^2-\abs{f(b^*a)}^2-\abs{f(b^*a)}^2+f(a^*a)f(b^*b)\right) \\
    &=f(a^*a)\left(f(a^*a)f(b^*b)-\abs{f(b^*a)}^2\right)
\end{align*}
Since $f$ is positive, we have $f(a^*a)\ge0$. Thus the second factor $f(a^*a)f(b^*b)-\abs{f(b^*a)}^2$ must also be non-negative, i.e. $\abs{f(b^*a)}^2\le f(a^*a)f(b^*b)$.
\end{proof}

\begin{theorem}\label{+lf}
Suppose $\As$ is a unital \csta\ and $f:\As\to\F$ is a linear functional. Then $f$ is positive if and only if it is bounded and $\norm{f}=f(\1)$.
\end{theorem}
\begin{proof}
($\Rightarrow$) Suppose $a\in\As$ is Hermitian and $\norm{a}<1$. We will show that there exists $b\in\As$ such that $\1-a=b^*b$. Note that the function $\phi:D(0,1)\to\C$, $\phi(z)=\sqrt{1-z}$ is holomorphic on $D(0,1)$, so it has a Taylor series $\sumn[0] r_nz^n$ that converges absolutely for all $\abs{z}<1$ (explicitly, $r_n=\frac{1}{4^n(1-2n)}\binom{2n}{n}$). Define $b=\sumn r_na^n$. Then $b\in\As$ and $b^2=1-a$. Also, since all $r_n\in\R$, we have $b^*=b$, and so $b^*b=b^2=\1-a$. This yields $f(\1)-f(a)=f(\1-a)=f(b^*b)\ge0$. Thus $f(a)\in\R$ and $f(\1)\ge f(a)$. Replacing $a$ with $-a$ yields $f(\1)\ge f(-a)=-f(a)$, and so $f(\1)\ge\abs{f(a)}$.\\

Now suppose $a\in\As$ and $\norm{a}<1$. Then $a^*a$ is Hermitian and $\norm{a^*a}=\norm{a}^2<1$. Thus $f(a^*a)\le f(\1)$. By \Cref{cauchycsta}, we have:
\begin{equation*}
    \abs{f(a)}^2=\abs{f(a^*)}^2=\abs{f(a^*\1)}^2\le f(\1^*\1)f(a^*a)=f(\1)f(a^*a)\le f(\1)^2
\end{equation*}
All in all, we have $\abs{f(a)}\le f(\1)$ for all $a\in\As$, $\norm{a}<1$. Thus $f$ is bounded and $\norm{f}\le f(\1)$. Finally, since $\norm{\1}=1$, we have $\norm{f}=f(\1)$.\\

($\Leftarrow$) Suppose $f$ is bounded and $\norm{f}=f(\1)$. We want to show that $f(c^*c)\ge0$ for all $c\in\As$. Assume without loss of generality that $\norm{c}\le1$ (otherwise, replace $c$ with $\frac{c}{\norm{c}}$). Then $\abs{f(\1)-f(c^*c)}=\abs{f(\1-c^*c)}\le\norm{f}\norm{\1-c^*c}=f(\1)\norm{\1-c^*c}\le f(\1)$
Thus $f(c^*c)\ge0$ for all $c\in\As$, $\norm{c}\le1$, and so $f$ is positive.
\end{proof}
\begin{remark}
If $\As$ is not unital, it still follows that every positive linear functional on $\As$ is bounded, which can be shown by unitizing $\As$ (see \Cref{embedding}).
\end{remark}

\subsection{The GNS Construction}\label{gnsconst}
\ind We are now ready to define the Gelfand-Naimark-Segal construction (or \emph{GNS construction}), which will led us to the \hyperref[gns]{GNS theorem}.

\begin{center}
\begin{minipage}{16em}\begin{center}
\personbox{RubineRed!50}{\href{https://mathshistory.st-andrews.ac.uk/Biographies/Gelfand/}{Israel M. Gelfand}\\(1913--2009)}{\includegraphics[width=1.3in]{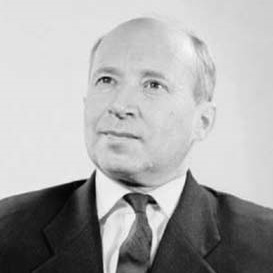}}\label{xgelfand} 
\end{center}\end{minipage}
\begin{minipage}{16em}\begin{center}
\personbox{PineGreen!50}{\href{https://mathshistory.st-andrews.ac.uk/Biographies/Naimark/}{Mark A. Naimark}\\(1909--1978)}{\includegraphics[width=1.3in]{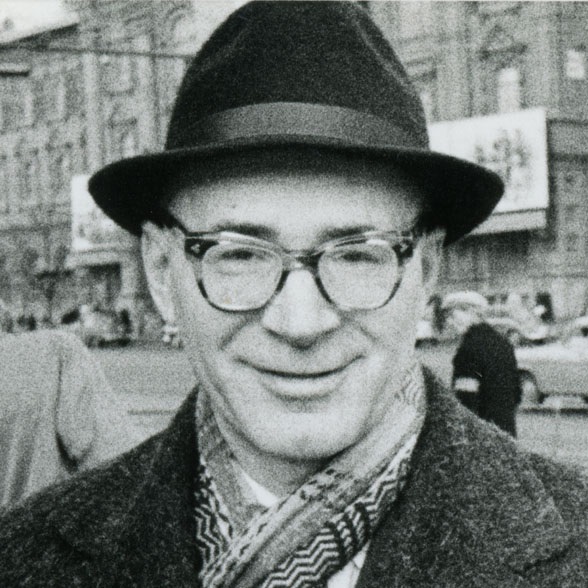}}\label{xnaimark} 
\end{center}\end{minipage}
\begin{minipage}{16em}\begin{center}
\personbox{Cerulean!50}{\href{https://mathshistory.st-andrews.ac.uk/Biographies/Segal/}{Irving E. Segal}\\(1918--1998)}{\includegraphics[width=1.3in]{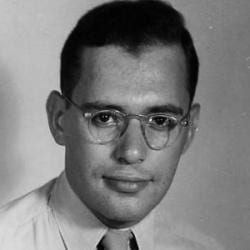}}\label{xsegal} 
\end{center}\end{minipage}
\end{center}

\begin{definition}
Suppose $\As$ is a \csta. A \textbf{representation} (or $^*$-representation) of $\As$ is a Hilbert space $\Hs$, together with a *-homomorphism $\pi:\As\to\Bs(\Hs)$.\\
A representation of $\As$ is \textbf{faithful} if the map $\pi$ is injective.
\end{definition}
\begin{remark}
Formally, the representation is the \emph{pair} $(\Hs,\pi)$, but we usually refer to $\pi$ itself as the representation, as is customary in representation theory.
\end{remark}

{
\definecolor{col}{HTML}{12C7AF}
\begin{tcolorbox}[title=The GNS Construction Explained, fonttitle=\sffamily\bfseries\selectfont,colback=white,colframe=col!60,coltitle=black,halign title=center,top=1mm,bottom=1mm,left=1mm,right=1mm,before skip=10pt,after skip=10pt,enhanced]
The steps of the GNS construction can be outlined as follows:
\begin{enumerate}
    \item Construct an inner product space from the positive linear functional $f$ and complete it to get a Hilbert space $\Hs_f$.
    \item For each element of the \csta\ $\As$, define a bounded operator on $\Hs_f$ by left-multiplication.
    \item Define a map $\pi_f$ sending each element of $\As$ to the corresponding operator.
\end{enumerate}
We now have a Hilbert space $\Hs_f$ and a *-homomorphism $\pi_f$. This is the GNS representation of $(\As,f)$.
\end{tcolorbox}
}

\ind It is important to note that the results in the rest of this section are constructed from an arbitrary \csta\ $\As$ as well as a positive linear functional $f:\As\to\C$. Thus whatever we construct in the end will (in general) depend on the choice of $f$.

\begin{theorem}
Suppose $\As$ is a \csta\ and $f$ is a positive linear functional on $\As$. Define $\Ns_f$ and $\ip{\cdot,\cdot}:(\As/\Ns_f)\times(\As/\Ns_f)\to\F$ as follows:
\begin{align*}
    \Ns_f=\left\{a\in\As\mid f(a^*a)=0\right\} && \ip{a+\Ns_f,b+\Ns_f}=f(b^*a)
\end{align*}
Then $\Ns_f$ is a closed left ideal of $\As$, and $\ip{\cdot,\cdot}$ is an inner product on $\As/\Ns_f$.
\end{theorem}
\begin{proof}
By \Cref{cauchycsta}, if $f(a^*a)=0$, then for all $b\in\As$, we have $\abs{f(b^*a)}^2\le f(a^*a)f(b^*b)=0$, and so $f(b^*a)=0$. Thus $\Ns_f$ is equivalently given by $\Ns_f=\{a\in\As\mid f(b^*a)=0\text{ for all }b\in\As\}$.\\

We first show that $\Ns_f$ is a closed left ideal of $\As$. Suppose $a\in\Ns_f$ and $b\in\As$. Then $f((ba)^*(ba))=f(a^*b^*ba)=f((b^*ba)^*a)=0$, so $ba\in\Ns_f$. Thus $\Ns_f$ is a left ideal of $\As$. It is also closed as it is the pre-image of $\{0\}$ under the continuous map $a\mapsto f(a^*a)$.\\

We now show that $\ip{\cdot,\cdot}$ is well-defined. Suppose $a,b,c,d\in\As$ such that $c\in a+\Ns_f$ and $d\in b+\Ns_f$, i.e. $a$ and $c$ are in the same equivalence class under $\Ns_f$, as are $b$ and $d$. Then we have $c=a+m$ and $d=b+n$ for some $m,n\in\Ns_f$. This yields:
\begin{align*}
    f(d^*c)
    &=f((b+n)^*(a+m))
    =f((b^*+n^*)(a+m))
    =f(b^*a+b^*m+n^*a+n^*m) \\
    &=f(b^*a)+f(b^*m)+f(n^*a)+f(n^*m)
    =f(b^*a)+f(b^*m)+\overline{f(a^*n)}+f(n^*m) \\
    &=f(b^*a)+0+\overline{0}+0
    =f(b^*a)
\end{align*}
Thus $\ip{\cdot,\cdot}$ is well-defined.\\

We now show that $\ip{\cdot,\cdot}$ is an inner product on $\As/\Ns_f$. We will denote equivalence classes $d+\Ns_f\in\As/\Ns_f$ by $[d]$.
\begin{enumerate}
    \item $\ip{a,a}=f(a^*a)\ge0$ for all $a\in\As$, and $\ip{a,a}=0\iff f(a^*a)=0\iff a\in\Ns_f\iff[a]=[0]$.
    \item $\ip{[a],[b]}=f(b^*a)=\overline{f(a^*b)}=\overline{\ip{[b],[a]}}$
    \item $\ip{\lam[a]+[b],[c]}=\ip{[\lam a+b],c}=f(c^*(\lam a+b))=f(\lam c^*a+c^*b)=\lam f(c^*a)+f(c^*b)=\lam\ip{[a],[c]}+\ip{[b],[c]}$
\end{enumerate}
Thus $\ip{\cdot,\cdot}$ is an inner product on $\As/\Ns_f$.
\end{proof}

\ind Now that we have constructed an inner product space $\As/\Ns_f$, we can start to devise a representation of $\As$. First, we complete the inner product space $\As/\Ns_f$ to get a Hilbert space, call it $\Hs$ (see \Cref{embedding} for details). We now need to define a *-homomorphism $\pi:\As\to\Bs(\Hs)$. To do this, we need the following result from functional analysis:

\begin{lemma}\label{bdopext}
Suppose $\Hs$ is a Hilbert space, $\Vs$ is a dense subspace of $\Hs$ and $S:\Vs\to\Hs$ is a bounded linear operator. Then there is exactly one $T\in\Bs(\Hs)$ such that $T|_\Vs=S$ and $\norm{T}=\norm{S}$.
\end{lemma}
In other words, every bounded linear operator defined on a dense subspace of $\Hs$ can be extended uniquely to a bounded linear operator on all of $\Hs$ without changing its operator norm. See \cite[Proposition 2.59, Page 58]{ward} or \cite[Theorem 9.28, Pages 168-169]{haase} for a proof.

\begin{definition}
Suppose $\As$ is a \csta\ and $f:\As\to\C$ is a positive linear functional. The \textbf{Gelfand-Naimark-Segal representation} (or \textit{GNS representation}) of $f$ is the representation $(\Hs_f,\pi_f)$, where $\Hs_f$ is the completion of the inner product space $\As/\Ns_f$, and $\pi_f:\As\to\Bs(\Hs_f)$ is given by:
\begin{align*}
    \pi_f(a):\Hs_f\to\Hs_f && \pi_f(a)(b+\Ns_f)=ab+\Ns_f
\end{align*}
And $\pi_f(a)$ is extended to all of $\Hs_f$ by \Cref{bdopext}.
\end{definition}
\begin{remark}
Formally, $(\Hs_f,\pi_f)$ is the GNS representation of the \emph{pair} $(\As,f)$, but we usually refer to it as the GNS representation of $f$, as the \csta\ $\As$ is usually clear from the context (and can be deduced from $f$, being its domain).
\end{remark}

\begin{proposition}
The GNS representation of $f$ is a representation of $\As$.
\end{proposition}
\begin{proof}
By construction, $\Hs_f$ is a Hilbert space. We will denote equivalence classes $d+\Ns_f\in\As/\Ns_f$ by $[d]$. Suppose $a,b,c\in\As$ and $\lam\in\F$. Then we have:
\begin{itemize}
    \item $\pi_f(\lam a+b)([c])=[(\lam a+b)c]=[\lam ac+bc]=\lam[ac]+[bc]=\lam\pi_f(a)([c])+\pi_f(b)([c])$
    \item $\pi_f(ab)([c])=[(ab)c]=[abc]=[a][bc]=[a]\pi_f(b)([c])=\pi_f(a)\pi_f(b)([c])$
    \item $\pi_f(a^*)[c]=[a^*c]=[a]^*[c]=\pi_f(a)^*([c])$
\end{itemize}
Thus $\pi_f$ is a *-homomorphism, and so the GNS representation of $f$ is a representation of $\As$.
\end{proof}
\begin{remark}
It follows directly from \Cref{*homobdd} that $\pi_f$ is bounded and $\norm{\pi_f}\le1$.
\end{remark}

\ind If $\pi_f$ was an isometry, we could immediately conclude the \hyperref[gns]{GNS theorem}, as we already have all the other ingredients for the proof. Clearly $\pi_f$ is not an isometry in general, as we could choose $f$ to be the zero functional, which would yield $\Ns_f=\As$, $\Hs_f=0$ and so $\pi_f(a)=0$ for all $a\in\As$. One might hope that we can choose a particular ``non-degenerate'' positive linear functional $f$ that makes $\pi_f$ an isometry. Unfortunately, this is not always possible.\\

\ind What we have to do instead is construct many representations from different positive linear functionals and ``patch'' them together to create a representation with no degeneracy\footnote{Think of a video game where you have to move an aircraft with a joystick. If you only have one joystick (say, that moves forward and backward), you can only move the aircraft forward or backward. If you have three joysticks, one forward/backward, one left/right and one up/down, you can move the aircraft anywhere you like. This is what we are doing with the direct sum representation: A single joystick (positive linear functional) is not enough, so we combine many joysticks to gain full control of the \csta.}. This is done via the \emph{direct sum representation}:

\begin{definition}
Suppose $\As$ is a \csta\ and $\{(\Hs_\alpha,\pi_\alpha)\}_{\alpha\in I}$ is a collection of representations of $\As$. The \textbf{direct sum} of $\{(\Hs_\alpha,\pi_\alpha)\}_{\alpha\in I}$ is given by:
\begin{equation*}
    \left(\bigoplus_{\alpha\in I} \Hs_\alpha,\bigoplus_{\alpha\in I} \pi_\alpha\right)
\end{equation*}
Where $\bigoplus_{\alpha\in I} \Hs_\alpha$ is the direct sum of the Hilbert spaces $\Hs_\alpha$, and $\bigoplus_{\alpha\in I} \pi_\alpha:\As\to\Bs\left(\bigoplus_{\alpha\in I} \Hs_\alpha\right)$ is the direct sum of the maps $\pi_\alpha:\As\to\Bs(\Hs_\alpha)$.
\end{definition}
\begin{remark}
The indexing set $I$ can be \emph{any} set, it can even be uncountable.
\end{remark}

\ind To explain the above definition a bit more, the direct sum (of the \emph{Hilbert spaces}) $\bigoplus_{\alpha\in I} \Hs_\alpha$ is the space of all elements $\{x_\alpha\}_{\alpha\in I}$ of the Cartesian product $\prod_{\alpha\in I} \Hs_\alpha$ such that $\sum_{\alpha\in I} \norm{x_\alpha}_{\Hs_\alpha}^2<\infty$, equipped with the inner product given by\footnote{You might be wondering if this is well-defined. The sum $\sum_{\alpha\in I} \norm{x_\alpha}_{\Hs_\alpha}^2$ is defined as the supremum of all partial sums of finitely many terms. If this is finite, there can only be countably many terms that are nonzero. Consequently, the sum $\sum_{\alpha\in I} \ip{x_\alpha,y_\alpha}$ is a countable sum of complex numbers, and by the Cauchy-Schwarz inequality, it is absolutely convergent. Thus the inner product on $\bigoplus_{\alpha\in I} \Hs_\alpha$ is well-defined, i.e. it is independent of the order in which the terms are added.} $\ip{\{x_\alpha\}_{\alpha\in I},\{y_\alpha\}_{\alpha\in I}}=\sum_{\alpha\in I} \ip{x_\alpha,y_\alpha}$. This is indeed a Hilbert space, see \cite[Page 24]{conway2}.\\

\ind The direct sum (of the \emph{maps}) $\bigoplus_{\alpha\in I} \pi_\alpha:\As\to\Bs\left(\bigoplus_{\alpha\in I} \Hs_\alpha\right)$ is the map that sends each $a\in\As$ to $(\bigoplus_{\alpha\in I} \pi_\alpha)(a)\in\Bs(\bigoplus_{\alpha\in I} \Hs_\alpha)$, which is the operator that sends $\{x_\alpha\}_{\alpha\in I}$ to $\{\pi_\alpha(x_\alpha)\}_{\alpha\in I}$. If $I$ is finite, this operator can be viewed as a block-diagonal matrix. The fact that $\bigoplus_{\alpha\in I} \pi_\alpha$ is a *-homomorphism follows easily from the definitions.

\begin{example}
Suppose $\As=\C$, $\Hs_1=\C$, $\Hs_2=\C^2$ and $\pi_1:\C\to M_1(\C)$ and $\pi_2:\C\to M_2(\C)$ are given by:
\begin{align*}
    \pi_1(z)=4z\mqty(1) && \pi_2(z)=z\mmat{0}{1}{2}{0}
\end{align*}
Where $(1)$ is the $1\times 1$ identity matrix. Then $\Hs_1\oplus\Hs_2=\C\times\C^2$, which can be identified with $\C^3$ in the natural way, and:
\begin{align*}
    \pi_1\oplus\pi_2:\C\to M_3(\C) \, , \qquad (\pi_1\oplus\pi_2)(z)=z\mqty(4&0&0\\0&0&1\\0&2&0)
\end{align*}
\end{example}

\begin{definition}
Suppose $\As$ is a \csta. The \textbf{universal representation} of $\As$ is given by:
\begin{equation*}
    (\mathbfcal{H},\boldsymbol\pi)=\left(\bigoplus_{f\in\text{S}(\As)} \Hs_f,\bigoplus_{f\in\text{S}(\As)} \pi_f\right)
\end{equation*}
Where $\text{S}(\As)$ is the set of all states on $\As$, and for each $f\in\text{S}(\As)$, $(\Hs_f,\pi_f)$ is the GNS representation of $f$.
\end{definition}

\begin{proposition}
The universal representation of $\As$ is faithful, and the map $\boldsymbol\pi$ is an isometry.
\end{proposition}
\begin{proof}
Suppose $a\in\As$. Assume without loss of generality that $\norm{a}\le1$ (otherwise, replace $a$ with $\frac{a}{\norm{a}}$). We will show that $\norm{\boldsymbol\pi(a)}_{\mathbfcal{H}}\ge\norm{a}_\As$. Since $(aa^*)^2=(aa^*)^*(aa^*)$ is positive, by \Cref{stateexists+}, there is a state $f\in\text{S}(\As)$ such that $f((aa^*)^2)=\norm{(aa^*)^2}$. This yields:
\begin{align*}
    \norm{[aa^*]}_{\Hs_f}^2
    &=\ip{[aa^*],[aa^*]}_{\Hs_f}
    =f((aa^*)^*(aa^*))
    =f\left((aa^*)^2\right)
    =\norm{(aa^*)^2}_\As
    =\norm{(aa^*)^*(aa^*)}_\As
    =\norm{aa^*}_\As^2
\end{align*}
Thus $\norm{[aa^*]}_{\Hs_f}=\norm{aa^*}_\As$. Since $\norm{a}\le1$, we have:
\begin{align*}
    \norm{\boldsymbol\pi(a)}_{\mathbfcal{H}}
    &\ge\norm{\pi_f(a)}_{\Hs_f}
    \ge\abs{\pi_f(a)[a^*]}_{\Hs_f}
    =\norm{[aa^*]}_{\Hs_f}
    =\norm{aa^*}_\As
    =\norm{a}_\As^2
\end{align*}
Thus $\norm{\boldsymbol\pi(a)}_{\mathbfcal{H}}\ge\norm{a}_\As$. Since $\boldsymbol\pi$ is a *-homomorphism, we also have $\norm{\boldsymbol\pi(a)}_{\mathbfcal{H}}\le\norm{a}_\As$, and so $\norm{\boldsymbol\pi(a)}_{\mathbfcal{H}}=\norm{a}_\As$. Thus $\boldsymbol\pi$ is an isometry. By \Cref{isomcontinj}, $\boldsymbol\pi$ is injective, and so $(\mathbfcal{H},\boldsymbol\pi)$ is faithful.
\end{proof}

\begin{theorem}[Gelfand-Naimark-Segal Theorem or GNS Theorem]\label{gns}
Every \csta\ is isometrically *-isomorphic to a closed *-subalgebra of $\Bs(\Hs)$, where $\Hs$ is some Hilbert space.
\end{theorem}
\begin{proof}
Suppose $\As$ is a \csta. Define $\mathbfcal{H}$ and $\boldsymbol\pi$ by the universal representation of $\As$. By the previous proposition, $\boldsymbol\pi$ is an isometric *-homomorphism, so its image $\im(\boldsymbol\pi)$ is a closed *-subalgebra of $\Bs(\mathbfcal{H})$. Thus $\boldsymbol\pi$ is an isometric *-isomorphism from $\As$ onto $\im(\boldsymbol\pi)$.
\end{proof}

\begin{example}
Suppose $\Hs$ is a Hilbert space. We already know that $\Bs(\Hs)$ is a \csta. Let's see what happens when we construct the universal representation of $\Bs(\Hs)$. By \Cref{ipx=state}, we can construct states on $\Bs(\Hs)$ by $f_x(A)=\ip{Ax,x}$, where $x\in\Hs$, $\norm{x}=1$. Each state $f_x$ naturally corresponds to at least one $x\in\Hs$. Thus the universal representation of $\Bs(\Hs)$ embeds it into a (possibly larger) \csta\ $\Bs(\mathbfcal{H})$.
\end{example}

\ind The way we concluded the \hyperref[gns]{GNS theorem} from the universal representation of a \csta\ is not the usual way this is done in the literature. Most sources, such as \cite{conway1}, \cite{conway2}, \cite{arveson1}, \cite{berkeleycsta} and \cite{verstraten}, use a different approach involving \emph{cyclic vectors} and \emph{cyclic representations}.

\conclbox{Chapter 3 Conclusion} 

\ind Through the \hyperref[gelfandnaimark]{Gelfand-Naimark theorem} and the \hyperref[gns]{GNS theorem}, we have classified \emph{all} \cstas: They can all be realized as a closed *-subalgebra of $\Bs(\Hs)$ for some Hilbert space $\Hs$. This shows that Segal's definition of a \csta\ is essentially equivalent to Rickart's definition: $B^*$-algebras and \cstas\ are one and the same. In addition, every \emph{abelian} \csta\ $\As$ can be realized as $C_0(X)$ for some locally compact Hausdorff space $X$, and in this case, we also have that $X$ is homeomorphic to $\Sigma(\As)$ (strictly speaking, we only have that $C_0(X)$ is isometrically isomorphic to $C_0(\Sigma(\As))$, but this implies that $X$ is homeomorphic to $\Sigma(\As)$ by the Banach-Stone theorem, see \cite[Chapter VI, Theorem 2.1, Page 172]{conway2}).\\

\ind There is much more to be said about \cstas, particularly the representation theory that arises from them. Since the GNS representation depends on our choice of positive linear functional(s), we can have many isometric *-isomorphisms from the same \csta\ into algebras of bounded operators on different Hilbert spaces. Some of these will have more desirable properties than others. This is further discussed in \cite{strung}, \cite{berkeleycsta} and \cite{thill}.\\

\ind There are also special classes of \cstas\ with additional properties, such as von Neumann algebras (also known as $W^*$-algebras), approximately finite-dimensional (AF) algebras and uniformly hyperfinite (UHF) algebras. See \cite{strung}, \cite{davidson}, \cite{conway1}, \cite{murphy} and \cite{thill} for more information on these.

\section{Rebuilding Quantum Mechanics}\label[chapter]{rebuildingqm}
\ind Now that we have developed all the necessary mathematical tools, we are ready see how they tie in to the theory of quantum mechanics. Our first goal is to draw the link between the two formulations of the Dirac-von Neumann axioms that we saw in \Cref{introduction}.\\

\ind The link between bounded operators on a Hilbert space and elements of a \csta\ was established in the \hyperref[gns]{GNS theorem}. This, together with \Cref{ipx=state}, shows that we can view states in the Hilbert space formulation (unit vectors in $\Hs$) as states in the \csta\ formulation (unit positive linear functionals on $\As$). This settles the question of compatibility of the two formulations.\\

\ind We will now consider the ``particle in a box'' model. Suppose we have a quantum particle of mass $m$ in a (one-dimensional) box $[0,L]$. The state space here is $\Hs=L^2([0,L])$. We can find the energy eigenstates, as one would in a first course on quantum mechanics, by solving the one-dimensional time-independent Schrödinger equation:
\begin{align*}
    -\frac{\hbar^2}{2m}\pdv[2]{\psi(x)}{x}+V(x)\psi(x)=0 &&
    V(x)=\begin{cases} 0 & 0\le x\le L \\ \infty & \text{otherwise} \end{cases}
\end{align*}
This yields the following energy eigenstates $\psi_n$ with corresponding energies $E_n$:

\begin{center}
\begin{tikzpicture}[scale=1,every node/.style={transform shape}]
\clip (-6,-0.2) rectangle (8,3.3);
\fill[gray!10,even odd rule] (-1,-0.2) rectangle (3,3.3) (0,0) rectangle (2,3.3);
\draw[<->,very thick] (0,3.3) node[below left=3pt]{$V$}--(0,0) node[left]{$0$}--(2,0)--(2,3.3);
\draw[-,thick,red] (0,0.2) sin (1,0.5) cos (2,0.2) node[right]{$\psi_1$};
\draw[-,thick,orange] (0,1) sin (0.5,1.3) cos (1,1) sin (1.5,0.7) cos (2,1) node[right]{$\psi_2$};
\draw[-,thick,Green] (0,1.8) sin (1/3,2.1) cos (2/3,1.8) sin (1,1.5) cos (4/3,1.8) sin (5/3,2.1) cos (2,1.8) node[right]{$\psi_3$};
\draw[-,thick,blue] (0,2.6) sin (0.25,2.9) cos (0.5,2.6) sin (0.75,2.3) cos (1,2.6) sin (1.25,2.9) cos (1.5,2.6) sin (1.75,2.3) cos (2,2.6) node[right]{$\psi_4$};
\node at (-4,1.5) {$\Disp \psi_n(x)=\sqrt{\frac{2}{L}}\sin\left(\frac{n\pi x}{L}\right)$};
\node at (5.5,1.5) {$\Disp E_n=\frac{n^2\pi^2\hbar^2}{2mL^2}$};
\end{tikzpicture}
\end{center}

\ind In this system, the position operator $\xh$ is \emph{bounded}, since for any state $\psi\in\Hs$, $\norm{\psi}_\Hs=1$, we have:
\begin{equation*}
    \norm{\xh\psi}_\Hs^2
    =\norm{x\psi(x)}_\Hs^2
    =\int_0^L \abs{x\psi(x)}^2\dx
    \le\int_0^L L^2\abs{\psi(x)}^2\dx
    =L^2
\end{equation*}
Intuitively, it \emph{must} be bounded, since the particle is confined to the region $[0,L]$.\\

\ind In the \csta\ formulation, the position operator $\xh$ will be a self-adjoint element of our \csta\ $\As$, which is itself a closed *-subalgebra of $\Bs(\Hs)$. How big this \csta\ is depends on how many observables we want to include: The more quantities you wish to measure, the larger your space needs to be. This is true even in the Hilbert space formulation.\\

\ind What about the states $\psi_n$? As we saw in the Dirac-von Neumann axioms, the states are now given by positive linear functionals $\omega_n:\As\to\C$. Again, what these are exactly will depend on the \csta\ $\As$, but we can work out what they will evaluate to for the observables we wish to consider. For example, evaluating the position operator on the eigenstates $\psi_n$, we get:
\begin{equation*}
    \ip{\xh\psi_n,\psi_n}
    =\int_0^L x\psi_n(x)\overline{\psi_n(x)}\dx
    =\int_0^L x\abs{\psi_n(x)}^2\dx
    =\frac{2}{L}\int_0^L x\sin\left(\frac{n\pi x}{L}\right)^2\dx
    =\frac{2}{L}\frac{L^2}{4}=\frac{L}{2}
\end{equation*}
In other words, the expected value of $\xh$ in any eigenstate is $\frac{L}{2}$. Thus in the \csta\ formulation, we have $\omega_n(\xh)=\frac{L}{2}$ for all $n\in\N$. This makes sense, as $\xh$ measures the position of the particle in a homogeneous box $[0,L]$, so on average, it should be $\frac{L}{2}$. For mixed states, the expected value of $\xh$ will depend on time, as the time evolution of each eigenstate is different (it is proportional to the energy $E_n$).\\

\ind We now consider another observable. Define $A:\Hs\to\Hs$ as follows:
\begin{equation*}
    A\psi(x)=-2\cos\left(\frac{2\pi x}{L}\right)\psi(x)
\end{equation*}
This is a self-adjoint operator in $\Bs(\Hs)$, so it is indeed an observable. The expected value of $A$ when the system is in the eigenstate $\psi_n$ is:
\begin{align*}
    \ip{A\psi_n,\psi_n}
    &=\int_0^L -2\cos\left(\frac{2\pi x}{L}\right)\psi_n(x)\overline{\psi_n(x)}\dx
    =-2\int_0^L \cos\left(\frac{2\pi x}{L}\right)\abs{\psi_n(x)}^2\dx \\
    &=-\frac{4}{L}\int_0^L \cos\left(\frac{2\pi x}{L}\right)\sin\left(\frac{n\pi x}{L}\right)^2\dx
    =\begin{cases} 1 & n=1 \\ 0 & n\ne1 \end{cases}
\end{align*}
This observable has an expected value of $1$ in the ground state $\psi_1$ and $0$ in every other eigenstate. Again, for mixed states, the expected value will depend on time. Evaluating these expected values for all the observables we wish to consider lets us construct the states $\omega_n:\As\to\C$ corresponding to the eigenstates $\psi_n$.\\

\ind The distinction between pure states and mixed states can be formalized, and is in fact already built into the \csta\ formulation: The state space of a quantum system is a \emph{convex cone} in the dual space $\As^*$ of the \csta\ $\As$, and the pure states of the system correspond to the \emph{extreme points} of this cone. The notions of convex cones and extreme points are part of a wider theory of \emph{topological vector spaces}. See \cite[Chapter 8]{ward} for more details.

\subsection{Unbounded Operators}
\ind There is one final issue that we have swept under the rug until now. While the \hyperref[gns]{GNS theorem} asserts that every \csta\ can be embedded into some $\Bs(\Hs)$ (the space of \emph{bounded} operators on a Hilbert space), it does not tell us anything about how to deal with \emph{unbounded} operators. This issue cannot be ignored, as many of the observables we deal with in quantum mechanics are unbounded. Fortunately though, most of the analysis still follows through. We will outline one instance of this here, but for a rigorous treatment, see \cite[Chapter VIII]{reedsimon}.\\

\ind As an example, we will look at the \emph{position} and \emph{momentum} operators (in one spatial dimension, though the general case follows similarly):
\begin{align*}
    \xh\psi(x)=x\psi(x) && \ph\psi(x)=-i\hbar\pdv{\psi(x)}{x}
\end{align*}
These are both self-adjoint\footnote{Strictly speaking, they are not self-adjoint, as they are not defined on all of $L^2(\R)$. They are not even densely defined. However, they are \emph{essentially self-adjoint}, i.e. they can be extended to closed self-adjoint operators, see \cite[Chapter 13]{ward}.} operators on $L^2(\R)$. They satisfy the \emph{canonical commutation relation}:
\begin{equation*}
    \xh\ph-\ph\xh=i\hbar\1
\end{equation*}
\ind At first glance, this seems to contradict \Cref{ab-ba/=1}. However, since both $\xh$ and $\ph$ are unbounded, they are not elements of the Banach algebra in question (in this case $\Bs(L^2(\R))$), and so they are able to satisfy this relation. This has the (unfortunate but interesting) consequence that quantum mechanics is an \emph{incomplete} theory of physics, as it does not allow us to simultaneously predict the position and momentum of any particle. This is formally manifested in the \emph{Heisenberg uncertainty principle}:
\begin{equation*}
    \Delta\xh\Delta\ph\ge\frac{\hbar}{2}
\end{equation*}\\

\ind We know from \Cref{hermunit} that self-adjoint \emph{bounded} operators on a Hilbert space can be transformed into unitary operators via the complex exponential. Although $\xh$ and $\ph$ are unbounded, we can try taking complex exponentials of them, just like we did in \Cref{exponential}:
\begin{align*}
    e^{ik\xh}\psi(x)
    &=\sumn[0]\frac{1}{n!}(ik\xh)^n\psi(x)
    =\sumn[0]\frac{1}{n!}(ik)^nx^n\psi(x)
    =\sumn[0]\frac{1}{n!}(ikx)^n\psi(x)
    =e^{ikx}\psi(x) \\
    e^{im\ph}\psi(x)
    &=\sumn[0]\frac{1}{n!}\left(i\left(-im\hbar\pdv{x}\right)\right)^n\psi(x)
    =\sumn[0]\frac{1}{n!}\left(-m\hbar\pdv{x}\right)^n\psi(x)
    =\sumn[0]\frac{1}{n!}(-m\hbar)^n\pdv[n]{\psi(x)}{x}
    =\psi(x-m\hbar)
\end{align*}
In the end, we see that $e^{ik\xh}$ is nothing but the phase shift operator, and $e^{im\ph}$ is nothing but the spatial translation operator. These are perfectly well-defined (and also unitary) operators on $L^2(\R)$!\\

\ind While our computation just now was not rigorous, it shows that unbounded operators can still meaningfully generate unitary operators via the complex exponential\footnote{There is a similar phenomenon in probability theory: Non-integrable random variables may still have a well-defined characteristic function, and this allows the central limit theorem to be applied in more general settings.}. This idea can be formalized through Stone's theorem on one-parameter unitary groups. See \cite[Section VIII.4]{reedsimon} and \cite{moller} for more details.

\xnewpage

\subsection{Concluding Remarks}
\ind So what have we done? We introduced the notion of \cstas\ in \Cref{cstas} and linked them to algebras of bounded operators on Hilbert spaces via the \hyperref[gnsconst]{GNS construction}. This shows that the \hyperref[diracvn]{Dirac-von Neumann axioms} are valid in the \csta\ formulation, and they are compatible with everything we already know about quantum mechanics from the Hilbert space formulation.\\

\ind But why is the \csta\ formulation any better? With simple quantum systems, such as a particle in a finite potential well, a particle in a box or a simple harmonic oscillator (these are among what physicists like to call `toy examples'), the Hilbert space formulation works just fine, and there is no reason to bring in \cstas. But for more complicated quantum systems, especially those with an infinite number of particles, the \csta\ formulation becomes essential.\\

\ind Suppose we want to study a quantum system consisting of many particles by analyzing them individually. We can create a state space for each particle and `patch' these spaces together (via a direct sum representation) to get the state space for the entire system. This sounds good in principle, but it has a few issues:
\begin{itemize}
    \item The states of the particles may depend on one another for physical reasons (e.g. quantum entanglement, Pauli exclusion principle).
    \item If the number of particles is infinite, we need to be careful when `patching' their state spaces together as this may not yield a Hilbert space of states (in the same way that an infinite sum of bounded functions or operators may not be bounded).
    \item If the particles are interacting with one another (e.g. radioactive decay, pair production), then the state spaces of each particle are insufficient to model the system, as the particles are not conserved throughout the experiment (particles may be created or destroyed).
\end{itemize}
If any of these issues is present in our system (which, physically speaking, is very likely), the Hilbert space formulation will be extremely cumbersome, or even inadequate, to properly analyze it.\\

\ind This is where the \csta\ formulation comes in: It solves \emph{all} of these issues at once. We can build the dependence between the states of particles into the linear functionals we use to construct our states, we can take direct sums of infinitely many representations just as easily as we do with finitely many, and we can account for particle interactions simply by introducing more observables (formally, we might have a larger \csta\ to deal with, but this is not intrinsically any more difficult as we have the \hyperref[gns]{GNS theorem} to help us).\\

\ind Besides showing how powerful \cstas\ are as a mathematical tool for quantum mechanics, the mere fact that it allows us to deal with systems of large (finite or infinite) numbers of particles is far-reaching in itself. Importantly, it allows us to study the collective behavior of particles in a macroscopic system, e.g. the kinetic theory of gases, even at the quantum level. This naturally leads on to the theory of \emph{quantum statistical mechanics}, which is discussed in \cite{bratteli} and \cite{bogolubov}.\\

\ind But that is another story......


\appendix

\section{Embedding}\label{embedding}
\ind In \Cref{banachalgs}, we introduced normed algebras in their most general form, but we were quick to specialize to unital complex Banach algebras. At that point, this seemed like a matter of necessity. Indeed, the \hyperref[ftba]{fundamental theorem of Banach algebras}, the \hyperref[specradius]{spectral radius formula} and the \hyperref[gkz]{GKZ theorem} require all three assumptions: \emph{unital}, \emph{complex} and \emph{Banach}. As far as the theory goes, however, we do not actually lose any generality by making these assumptions. This is because every normed algebra ``lives inside'' a unital complex Banach algebra. More formally, it can be embedded into a unital complex Banach algebra while preserving its algebra operations and its norm. We will now outline how to do this:


\begin{center}
\begin{tikzpicture}[scale=1,font=\sffamily]
\definecolor{col1}{HTML}{A8C4FF}
\definecolor{col2}{HTML}{4028A0}
\newtcbox{\stbox}[1][]{colback=#1!10,colframe=#1,coltitle=black,hbox,halign=center,rounded corners,top=2pt,bottom=2pt,left=2pt,right=2pt}
\node (a) at (-6,2) {\stbox[col1]{Normed algebra}};
\node (b) at (-3.5,0) {\stbox[col1!66.7!col2]{Banach algebra}};
\node (c) at (2.5,0) {\stbox[col1!33.3!col2]{Unital Banach algebra}};
\node (d) at (6,2) {\stbox[col2]{Unital complex Banach algebra}};
\draw[->,thick] (a)--node[midway,left,yshift=-1mm]{Completion} (b);
\draw[->,thick] (b)--node[midway,above]{Unitization} (c);
\draw[->,thick] (c)--node[midway,right,yshift=-1.5mm]{Complexification} (d);
\end{tikzpicture}
\end{center}


\begin{cthm}[Completion]
Suppose $\As$ is a normed algebra. Define $\mathfrak{C}$ as the set of all Cauchy sequences in $\As$, and define the equivalence relation $\sim$ on $\mathfrak{C}$ by $(a_n)\sim(b_n)\Rightarrow\limn\norm{a_n-b_n}=0$. Finally, define $\widehat{\As}=\mathfrak{C}/\sim$ with the following operations:
\begin{align*}
    [(a_n)]+[(b_n)]=[(a_n+b_n)] && \lam[(a_n)]=[(\lam a_n)] && [(a_n)][(b_n)]=[(a_nb_n)] && \norm{[(a_n)]}=\limn\norm{a_n}
\end{align*}
Then $\widehat{\As}$ is a Banach algebra, and the function $f:\As\to\widehat{\As}$, $f(a)=[(a)]$ ($f$ maps $a$ to the equivalence class of the constant sequence $(a,a,a,...)$) is an isometry. $\widehat{\As}$ is known as the \textbf{completion} of $\As$.
\end{cthm}
See \cite[Proposition 7.17, Page 107]{muscat} or \cite[Theorem 2.32]{ward} for a proof of the corresponding theorem for normed vector spaces, and \cite[\S1, Proposition 12]{bonsall} for a proof of the extension to normed algebras. If $\As$ is an inner product space (not necessarily an algebra), we can also extend the inner product to the completion by defining $\ip{[(a_n)],[(b_n)]}=\limn\ip{a_n,b_n}$.

\begin{cthm}[Unitization]
Suppose $\As$ is a non-unital Banach algebra. Define $\As_1=\As\times\F$ with the following operations:
\begin{align*}
    (a,x)+(b,y)&=(a+b,x+y) & \lam(a,x)&=(\lam a,\lam x) \\ (a,x)(b,y)&=(ab+xb+ya,xy) & \norm{(a,x)}_1&=\norm{a}+\abs{x}
\end{align*}
Then $\As_1$ is a unital Banach algebra with identity $(0,1)$. Also, the function $f:\As\to\As_1$, $f(a)=(a,0)$ is an isometric homomorphism. $\As_1$ is known as the \textbf{unitization} of $\As$.
\end{cthm}
See \cite[Lemma 1.4]{wilde} or \cite[Proposition I.1.3, Pages 2-3]{davidson} for a proof.

\begin{cthm}[Complexification]
Suppose $\As$ is a unital real Banach algebra. Define $\As_\C=\As\times\As$ with the following operations:
\begin{align*}
    (a,b)+(c,d)=(a+c,b+d) && (\mu+i\nu)(a,b)=(\mu a-\nu b,\mu b+\nu a) && (a,b)(c,d)=(ac-bd,ad+bc)
\end{align*}
For each $a\in\As$, define $T_a\in\Bs(\As)$ by $T_a(b)=ab$ (this is known as the \textit{left regular representation} of $\As$). Define $T_a'\in\Bs(\As_\C)$ by $S_a(b,c)=(T_a(b),T_a(c))=(ab,ac)$. Finally, define the norm $\norm{\cdot}_\C$ on $\As_\C$ by $\norm{(a,b)}_\C=\norm{S_a+iS_b}$. Then $\As_\C$ is a unital complex Banach algebra with identity $(\1,0)$, and the function $f:\As\to\As_\C$, $f(a)=(a,0)$ is an isometric homomorphism. $\As_\C$ is known as the \textbf{complexification} of $\As$.
\end{cthm}
See \cite[Theorem 1.3.1]{rickart} for a proof.

\begin{example}
Suppose $X$ is a locally compact (but not compact) Hausdorff space and $\As=C_c(X,\R)$, the set of all continuous functions $f:X\to\R$ with compact support. This normed algebra is not unital, complex or Banach. We can give it all three of these properties by successively applying the above procedures:
\begin{enumerate}
    \item The completion of $C_c(X,\R)$ is $C_0(X,\R)$, the set of all continuous functions $f:X\to\R$ that vanish at infinity.
    \item The unitization of $C_0(X,\R)$ is $C_0(X,\R)\oplus\R1$, the set of all functions $f:X\to\R$ that can be expressed as the sum of a function in $C_0(X,\R)$ and a (real) constant function.
    \item The complexification of $C_0(X,\R)\oplus\R1$ is $C_0(X,\C)\oplus\C1$, the set of all functions $f:X\to\C$ that can be expressed as the sum of a function in $C_0(X,\C)$ and a (complex) constant function.
\end{enumerate}
As intended, $C_0(X,\C)\oplus\C1$ is a unital complex Banach algebra.
\end{example}


\section{The Weak* Topology}\label{weak*top}
\ind In \Cref{gelfspec} and \Cref{gelfspecch}, we used the notions of nets and the weak* topology. Here, we will elaborate a little on these.

\begin{definition}
Suppose $X$ is a normed vector space. For each $x\in X$, the \textbf{evaluation map} of $x$ is given by $\xh:X^*\to\F$, $\xh(f)=f(x)$.\\
The \textbf{weak* topology} on $X^*$ is the smallest topology on $X^*$ that makes all of the evaluation maps continuous.
\end{definition}
\begin{remarks}\leavevmode
\begin{enumerate}
    \item A sequence $(f_n)$ in $X^*$ converges to $f\in X^*$ in the weak* topology if and only if it converges pointwise to $f$, i.e. $f_n(x)\to f(x)$ (in $\F$) for all $x\in X$.
    \item The weak* topology is smaller than the weak topology, which is smaller than the strong topology.
    \item If $X$ is a reflexive Banach space (in particular, if $X$ is a Hilbert space), then the weak* topology is identical to the weak topology.
\end{enumerate}
\end{remarks}

\begin{definition}
A \textbf{directed set} is a non-empty set $P$, together with a binary relation $\preceq$ on $P$, such that:
\begin{enumerate}
    \item If $a\in P$, then $a\preceq a$. \hfill (Reflexivity)
    \item If $a,b,c\in P$, $a\preceq b$ and $b\preceq c$, then $a\preceq c$. \hfill (Transitivity)
    \item If $a,b\in P$, then there exists $c\in P$ such that $a\preceq c$ and $b\preceq c$. \hfill (Upper bound property)
\end{enumerate}
A \textbf{net} in a set $X$ is a function $f:X\to P$, where $P$ is any directed set.
\end{definition}
\begin{examples}\leavevmode
\begin{enumerate}
    \item $\N$ is a directed set with the relation $\le$. The resulting nets are simply sequences, e.g. $a_n=2n$.
    \item $\R$ is a directed set with the relation $\le$. The resulting nets are functions on $\R$, e.g. $f(x)=x^2-1$.
\end{enumerate}
\end{examples}

\ind You might be wondering why we used nets instead of sequences to prove \Cref{gelfspecch}. This is a necessary adjustment, due to the following result:
\begin{theorem}
Suppose $X$ is an infinite-dimensional normed vector space. Then the weak* topology on $X^*$ is not first countable.
\end{theorem}
See \cite[Theorem 6.26, Pages 237-238]{aliprantis} for a proof.\\

\ind Since the weak* topology on $X^*$ is not first countable, sequences are not general enough to describe its topological behavior, and so we have to resort to nets. See \cite{singh} for a more detailed discussion of nets and first countability.

\begin{theorem}[Banach-Alaoglu Theorem]\label{banachalaoglu}
Suppose $X$ is a normed vector space. Then the closed unit ball in $X^*$ is weak*-compact.
\end{theorem}
See \cite[Theorem 1.23, Page 9]{douglas} or \cite[Chapter V, Theorem 3.1]{conway2} for a proof.
\begin{warning}
This does NOT imply that $X^*$ is weak*-locally compact! While every closed ball centered at any $f\in X^*$ is weak*-compact, it is not a weak*-neighborhood of $f$. In fact, it is not even a weak neighborhood of $f$, see \cite[Corollary 6.27, Page 238]{aliprantis}.
\end{warning}

\newpage

\clearpage 
\phantomsection 
\addcontentsline{toc}{section}{References}
\printbibliography

\vspace{5mm}


{
\definecolor{col}{HTML}{5221C8}
\begin{tcolorbox}[title=Photo Credits, fonttitle=\Large\sffamily\bfseries\selectfont,colback=col!10,colframe=col!40,coltitle=black,halign title=center,top=1mm,bottom=1mm,left=1mm,right=1mm,drop fuzzy shadow,before skip=10pt,after skip=10pt,enhanced]
\begin{itemize}[align=left,leftmargin=0.4in,font=\sffamily\selectfont]
    \item[{\hyperref[xdirac]{Dirac}:}] \url{https://www.wsj.com/articles/the-pleasure-and-pain-of-scientific-predictions-11597935357}
    \item[{\hyperref[xvonneumann]{von Neumann}:}] \url{https://en.wikipedia.org/wiki/John_von_Neumann}
    \item[{\hyperref[xgleason]{Gleason}:}] \url{https://www.geni.com/people/Andrew-Gleason/6000000030388180603}
    \item[{\hyperref[xkahane]{Kahane}:}] \url{http://www.bibmath.net/bios/index.php?action=affiche&quoi=kahane}
    \item[{\hyperref[xzelazko]{Żelazko}:}] \url{https://www.wikidata.org/wiki/Q657990}
    \item[{\hyperref[xrickart]{Rickart}:}] \url{https://peoplepill.com/people/charles-earl-rickart/}
    \item[{\hyperref[xgelfand]{Gelfand}:}] \url{https://en.wikipedia.org/wiki/File:IM_Gelfand.jpg}
    \item[{\hyperref[xnaimark]{Naimark}:}] \url{https://www.maa.org/book/export/html/117813}
    \item[{\hyperref[xsegal]{Segal}:}] \url{https://www.gf.org/fellows/all-fellows/irving-e-segal/}
\end{itemize}
\end{tcolorbox}
}

\end{document}